\documentclass[12pt
]{article}

\usepackage[all]{xy}
\usepackage[latin1]{inputenc}
\usepackage{amsfonts}
\usepackage{amsmath}
\usepackage{amssymb}
\usepackage{amscd}
\usepackage{amsthm}
\usepackage{indentfirst}
\usepackage[hmargin=2.98cm
,vmargin=4.0cm
]{geometry}

\usepackage{tikz-cd}

\usepackage[colorlinks=true,linkcolor=blue,
  citecolor=blue,
  urlcolor=blue,
  filecolor=blue]{hyperref}

\newtheorem{thm}{Theorem}[section]
\newtheorem{lemma}[thm]{Lemma}

\newtheorem{prop}[thm]{Proposition}
\newtheorem{coroll}[thm]{Corollary}
\newtheorem{claim}[thm]{Claim}

\theoremstyle{definition}

\newtheorem{defin}[thm]{Definition}
\newtheorem{rem}[thm]{Remark}

\newtheorem{exam}[thm]{Example}

\newtheorem{notation}[thm]{Notation}
\newtheorem*{acknow}{Acknowledgements}

\newcommand{\R}{{\mathbb{R}}}

\newcommand{\Q}{{\mathbb{Q}}}

\newcommand{\Z}{{\mathbb{Z}}}
\newcommand{\N}{{\mathbb{N}}}
\newcommand{\C}{{\mathbb{C}}}

\newcommand{\CP}{{\mathbb{C}P}}

\newcommand{\cA}{{\mathcal{A}}}
\newcommand{\cB}{{\mathcal{B}}}
\newcommand{\cC}{{\mathcal{C}}}

\newcommand{\cF}{{\mathcal{F}}}
\newcommand{\cG}{{\mathcal{G}}}
\newcommand{\cH}{{\mathcal{H}}}
\newcommand{\cI}{{\mathcal{I}}}
\newcommand{\cJ}{{\mathcal{J}}}
\newcommand{\cK}{{\mathcal{K}}}

\newcommand{\cM}{{\mathcal{M}}}

\newcommand{\cP}{{\mathcal{P}}}

\newcommand{\cR}{{\mathcal{R}}}

\newcommand{\fc}{\colon}

\newcommand{\ol}{\overline}
\newcommand{\wt}{\widetilde}
\newcommand{\wh}{\widehat}

\newcommand{\tb}{\textbf}

\newcommand{\bs}{\boldsymbol}

\DeclareMathOperator{\sgn}{sgn}

\DeclareMathOperator{\im}{im}
\DeclareMathOperator{\id}{id}
\DeclareMathOperator{\area}{area}
\DeclareMathOperator{\Vol}{Vol}

\DeclareMathOperator{\Slit}{Slit}

\DeclareMathOperator{\val}{val}
\DeclareMathOperator{\pr}{pr}

\DeclareMathOperator{\ev}{ev}

\DeclareMathOperator{\ini}{ini}
\DeclareMathOperator{\ter}{ter}

\DeclareMathOperator{\cone}{cone}
\DeclareMathOperator{\co}{co}

\DeclareMathOperator{\supp}{supp}

\DeclareMathOperator{\Int}{Int}

\DeclareMathOperator{\coker}{coker}

\DeclareMathOperator{\res}{res}

\begin{document}
	
	\title{Symplectophagy}
	\author{Adi Dickstein\footnote{Partially supported by the Fondation Courtois, the Jenny Panitch Beckow Memorial Scholarship, Canada, 2026, and by the Israel Science Foundation (grant no.\ 1102/20)}, Yaniv Ganor, and Frol Zapolsky}
	
	\date{}
	
	\setcounter{tocdepth}{2}
	
	\renewcommand{\labelenumi}{(\roman{enumi})}

	\maketitle
	\begin{abstract}
      We introduce the notion of a symplectophagy, which is a symplectic isotopy satisfying conditions reminiscent of the symplectic camel. We present new symplectic rigidity results which take the form of obstructions to the existence of symplectophagies in certain situations. The main tool is the relative symplectic cohomology recently introduced by Varolg\"une\c s, as well as its extension to pairs, due to the authors and L.\ Polterovich, for which we establish an invariance property, which is the chief technical novelty of the present paper. An important role is also played by the recent calculations of relative symplectic cohomology in complex projective spaces due to Dickstein--Ganor.
	\end{abstract}

    \tableofcontents

\section{Introduction and main results}\label{s:intro}

A central theme in symplectic topology is the dichotomy between rigidity and flexibility, see the introduction to \cite{MS_1995_Intro}. Flexibility refers to phenomena in symplectic geometry governed by algebraic topology, whereas rigidity is the existence of additional obstructions having, in a sense, a purely symplectic character. One of the first examples of rigidity was Gromov's celebrated nonsqueezing theorem \cite{Gromov_1985_J-hol_curves}, which elucidated a profound difference between symplectic and volume-preserving geometry. Nowadays symplectic rigidity pervades the field. One of the early examples beyond restrictions on symplectic embeddings was the symplectic camel problem, which deals with a certain kind of symplectic isotopies \cite{Eliashberg_Gromov_1991_convex_symp_mfds}, \cite{McDuff_Traynor_1993_camel}, \cite{Viterbo_1992_GenFun}.

In the present paper we define a new class of symplectic isotopies called \emph{symplectophagies}, which map one subset of a symplectic manifold into another under certain restrictions, and thus resemble the symplectic camel. Our main results are obstructions to the existence of symplectophagies in certain situations. These obstructions come from the theory of relative symplectic cohomology for compact subsets of a symplectic manifold, recently introduced by Varolg\"une\c s \cite{Varolgunes_2021_MV_and_relSH}, as well as an enhancement thereof to pairs of compact subsets, defined in \cite{DGPZ_2024_Symp_top_and_IVMs} by the authors together with L.\ Polterovich. The information we extract from these theories which allows us to prove the rigidity results comes from the \emph{torsion} of the relative symplectic cohomology modules \emph{over the Novikov ring}, and takes the form of capacity-like quantities.

Another crucial ingredient is the recent calculation by Dickstein--Ganor of the relative symplectic cohomology in complex projective spaces for toric balls, their complements and boundaries, together with appropriate restriction maps. The main technical novelty in the present paper, which is necessary for some of our results, is the invariance property for the relative symplectic cohomology of pairs, Theorem \ref{thm: invariance}.

\subsection{Definition and examples}\label{ss:def_examples}

Throughout the paper, by an \tb{isotopy} we always mean one that is parametrized by $[0,1]$ and that starts at the identity at time $0$, unless explicitly stated otherwise.

\begin{defin}\label{def:symplectophagy}
    Let $(M,\omega)$ be a symplectic manifold, and let $A, B \subset M$. A Hamiltonian isotopy $(\phi_t)_{t\in[0,1]}$ is called a \textbf{weak symplectophagy of $B$ into $A$} if:
    \begin{itemize}
        \item $\phi_1(B)\subset A$,
        \item for every $t\in [0,1]$, we have $\phi_t(B)\subset A\cup B$.
    \end{itemize}
    If in addition for every $t\in [0,1]$, we have $\phi_t(A\cap B)\subset A$ then $(\phi_t)_{t\in[0,1]}$ is called a \textbf{genuine symplectophagy} of $B$ into $A$.
\end{defin}

\begin{rem}\label{rem:symplectophagy_name_cpt_suppt}\phantom{M}
  \begin{itemize}
    \item The suffix ``\emph{-phagy}'' comes from Greek and refers to eating or consumption. A symplectophagy of $B$ into $A$ can be viewed as a process of consumption, at the end of which $B$ is entirely contained in $A$. In the case of a genuine symplectophagy the piece $A\cap B$, having already been swallowed by $A$, can never leave it.
    \item We do not impose conditions on the support of the Hamiltonian isotopy. However, note that if $(\phi_t)_t$ is a weak (respectively, genuine) symplectophagy of $B$ into $A$ and $B$ is compact, then there is a compactly supported Hamiltonian isotopy which is a weak (respectively, genuine) symplectophagy of $B$ into $A$. This can be achieved by suitably cutting off the generating Hamiltonian.
  \end{itemize}

\end{rem}

Our main goal in this paper is to prove rigidity results in the form of obstructions to symplectophagy in some situations. There are certain \emph{obvious} obstructions to symplectophagy, among which:
\begin{itemize}
  \item \emph{Symplectic embeddability}: if there is no symplectic embedding of $B$ into $A$, for instance if $\Vol(B)>\Vol(A)$, then, of course, there can be no symplectophagy of $B$ into $A$.
  \item \emph{Topological obstructions}: as we will see, sometimes there is no isotopy, not necessarily symplectic, satisfying the above containment requirements, for purely topological reasons, in which case there is no symplectophagy either.
\end{itemize}

Our rigidity results are formulated in Section \ref{ss:main_results}. Beforehand, however, we would like to consider some illustrative examples, which show that the notion of symplectophagy is related to existing phenomena, on the one hand, and that at times it can be quite delicate, on the other hand. In all the examples we only consider the case where the symplectic volume of $B$ does not exceed that of $A$, see Remark \ref{rem:symplectophagy_name_cpt_suppt}.
\begin{exam}\label{ex:two_disks}
  Let $A,B\subset \R^2$ be two smooth closed disks with transversely intersecting boundaries (Figure \ref{fig:two_disks}). If $A\cap B=\varnothing$, there is no weak symplectophagy of $B$ into $A$, because an isotopy mapping $B$ into $A\cup B$ must keep $B$ inside $B$. If $A\cap B\neq \varnothing$, then there is a weak symplectophagy of $B$ into $A$, because $B$ can be ``dragged'' into $A$ through one of the components of $A\cap B$. If $A\cap B$ has a single connected component, then a similar argument shows that there exists a genuine symplectophagy of $B$ into $A$. If $A\cap B$ has more than one connected component, then there is no genuine symplectophagy, due to a purely topological obstruction, see Example \ref{ex:top_obstructions} below.
  \begin{figure}
  \begin{minipage}[]{.45\textwidth}
    \centering
    \includegraphics[width=0.6\linewidth]{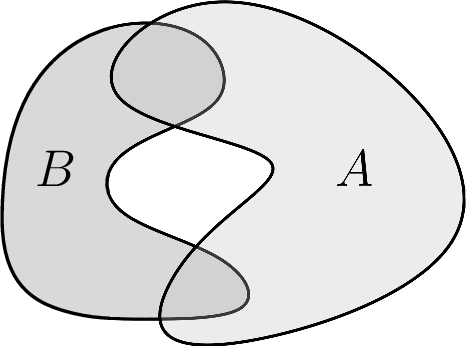}

    (a)
  \end{minipage}
  \begin{minipage}[]{.45\textwidth}
    \centering
    \includegraphics[width=0.5\linewidth]{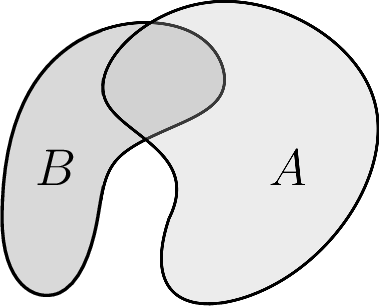}

    (b)
  \end{minipage}
  \caption{(a) two intersection components; (b) one intersection component}
  \label{fig:two_disks}
  \end{figure}
\end{exam}

\begin{exam}\label{ex:disk_annulus}
  Let $D\subset \R^2$ be a smooth closed disk and $A$ be a smooth closed annulus of positive area such that $\partial A \cap D = \partial D$, see Figure \ref{fig:disk_annulus_sphere}(a). Then there is no weak symplectophagy of $A$ into $D$. The obstruction here is a combination of topological and area considerations. Indeed, any symplectomorphism of $\R^2$ mapping $A$ into $D$ must map the whole disk bounded by the outer boundary circle of $A$ into $D$, by the Jordan curve theorem. But since this area is larger than that of $D$, we arrive at a contradiction.
  \begin{figure}
  \begin{minipage}[]{.45\textwidth}
    \centering
  	\includegraphics[width=0.5\linewidth]{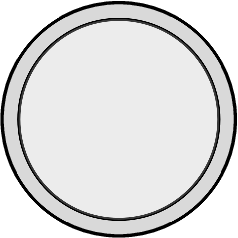}

  (a)
  \end{minipage}
  \begin{minipage}[]{.45\textwidth}
    \centering
    \includegraphics[width=0.6\linewidth]{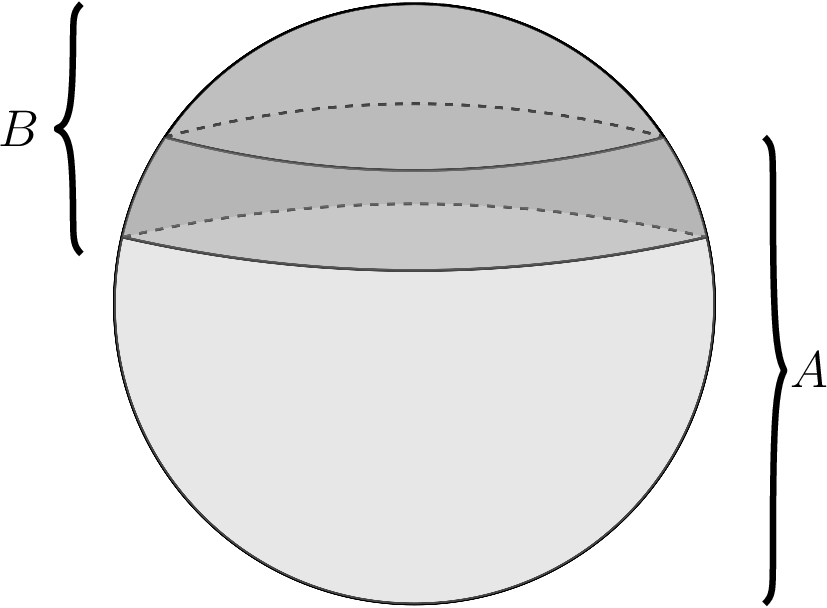}

    (b)
  \end{minipage}
  \caption{(a) disk and annulus; (b) two disks covering the sphere}
  \label{fig:disk_annulus_sphere}
  \end{figure}
\end{exam}

\begin{exam}\label{ex:sphere}
  Let $S^2 = \{z\in\R^3\,|\,\|z\|=1\}$ be the unit sphere with the standard area form, and for $-1 < b \leq a < 1$ define $A,B\subset S^2$ by $A = \{z\in S^2\,|\,z_3\leq a\}$, $B = \{z\in S^2\,|\,z_3\geq b\}$ (Figure \ref{fig:disk_annulus_sphere}(b)). If $\area(B) = 2\pi(1-b) \leq 2\pi(a+1)=\area(A)$, then there is a weak symplectophagy of $B$ into $A$ inside $S^2$, given, for example, by a suitable rotation. There is no genuine symplectophagy of $B$ into $A$ for any $a,b$, the obstruction again being purely topological, see Example \ref{ex:top_obstructions}.
\end{exam}

\begin{exam}
  Moving on to higher dimensions, for $0<b\leq a$ consider $A = [0,a]\times[0,1]^{2n-1}$, $B=[-b,0]\times[0,1]^{2n-1}$ (Figure \ref{fig:boxes}(a)). A suitably parametrized linear translation along the first coordinate of $\R^{2n}=\R\times\R^{2n-1}$ is a genuine symplectophagy of $B$ into $A$.
  \begin{figure}
  \begin{minipage}[]{.25\textwidth}
    \centering
  	\includegraphics[width=.9\linewidth]{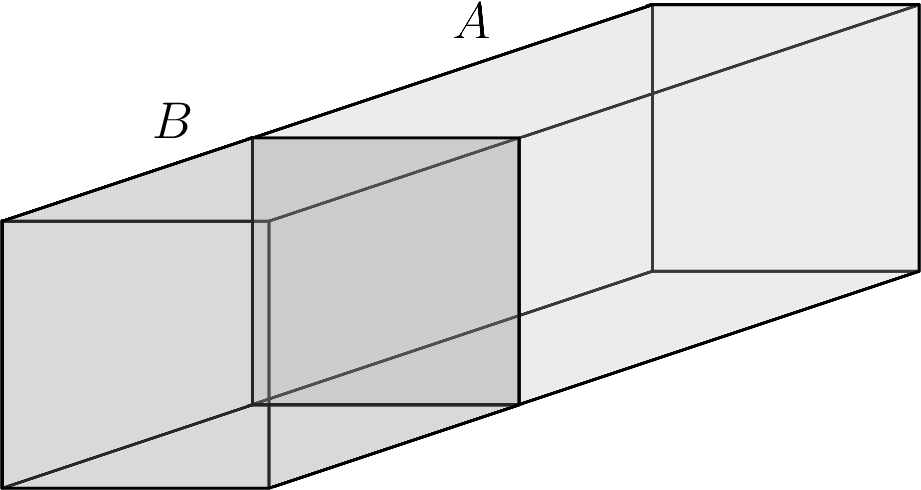}

  (a)
  \end{minipage}
  \begin{minipage}[]{.35\textwidth}
    \centering
    \includegraphics[width=.9\linewidth]{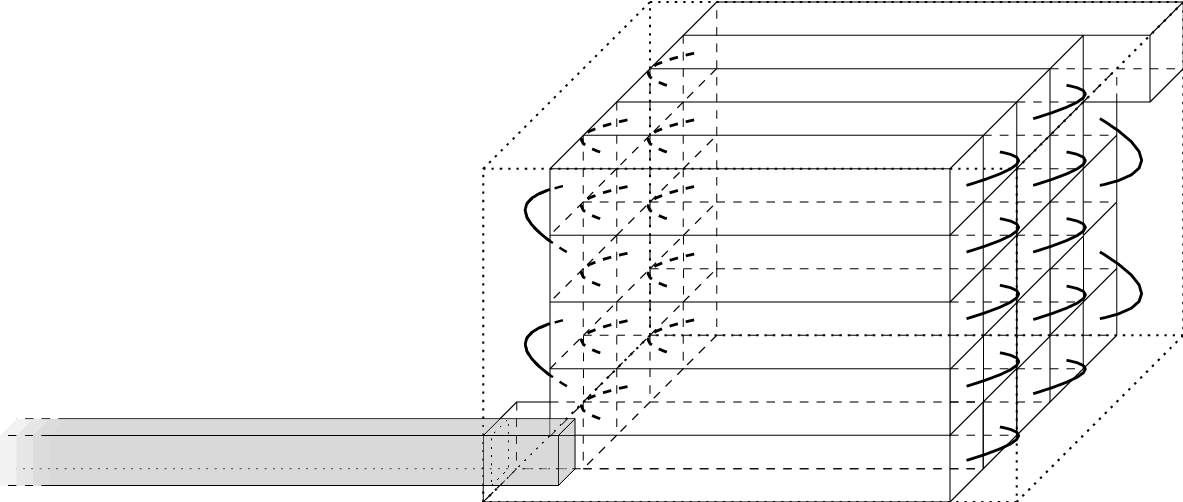}

    (b)
  \end{minipage}
  \begin{minipage}[]{.35\textwidth}
    \centering
    \includegraphics[width=.9\linewidth]{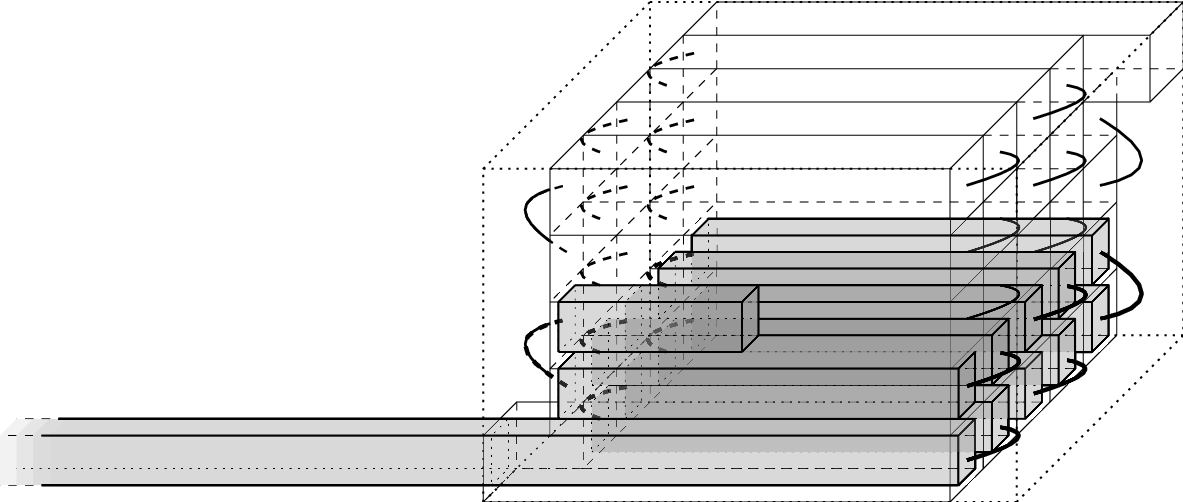}

    (c)
  \end{minipage}
  \caption{(a) boxes; (b) before symplectophagy; (c) during symplectophagy}
  \label{fig:boxes}
  \end{figure}
\end{exam}

\begin{exam}\label{ex:sympl_embed_polydisk_cube}
  Fix an odd integer $k\geq 5$. In this example for each $\epsilon\in(0,\frac 12)$, we will construct a genuine symplectophagy of a thin box $B = [1-\alpha+\epsilon,1-\epsilon]\times[\epsilon,1-\epsilon]^{2n-1}$ into the cube $A = (0,k)^{2n}$, where $\alpha$ is a certain positive number which is at least $k^{2n-1}(k-4)$ (Figure \ref{fig:boxes}(b-c)). In \cite{DZ_2026_constraints_on_symp_qs} there is a construction of a symplectic embedding $\iota$ of a polydisk $(0,\alpha)\times(0,1)^{2n-1}$ into the cube $(0,k)^{2n}$, for some $\alpha \geq k^{2n-1}(k-4)$. This embedding is the identity map when restricted to $(0,1)^{2n}$. Extend $\iota$ to a symplectic embedding $\wh\iota \fc (-\infty,\alpha)\times(0,1)^{2n-1} \to \R^{2n}$ by $\wh\iota|_{(0,\alpha)\times(0,1)^{2n-1}}=\iota$, $\wh\iota|_{(-\infty,1)\times(0,1)^{2n-1}}=\id$.

  For $t\in[0,\alpha-1]$ define $\iota_t \fc (1-\alpha,1)\times(0,1)^{2n-1}\to \R^{2n}$ to be the composition $\wh\iota\circ T_{(t,0,\dots,0)}|_{(1-\alpha,1)\times(0,1)^{2n-1}}$, where $T_x$ is the translation operator $y\mapsto y+x$ for $x \in \R^{2n}$. Clearly $\iota_0 = \id$ while $\iota_{\alpha-1}=\iota\circ T_{(\alpha-1,0,\dots,0)}$. The ``extension after restriction principle'' (\cite[p. 7]{Schlenk_book}) applies parametrically for parameters lying in a compact contractible space, here $[0,\alpha-1]$, and yields a Hamiltonian isotopy $(\phi_t)_{t\in[0,\alpha-1]}$ with $\phi_t|_B = \iota_t|_B$. A direct verification shows that $(\phi_{(\alpha-1)t})_{t\in[0,1]}$ is a genuine symplectophagy of $B$ into $A$.
  \begin{figure}
  \begin{minipage}[]{.45\textwidth}
    \centering
  	\includegraphics[width=.5\linewidth]{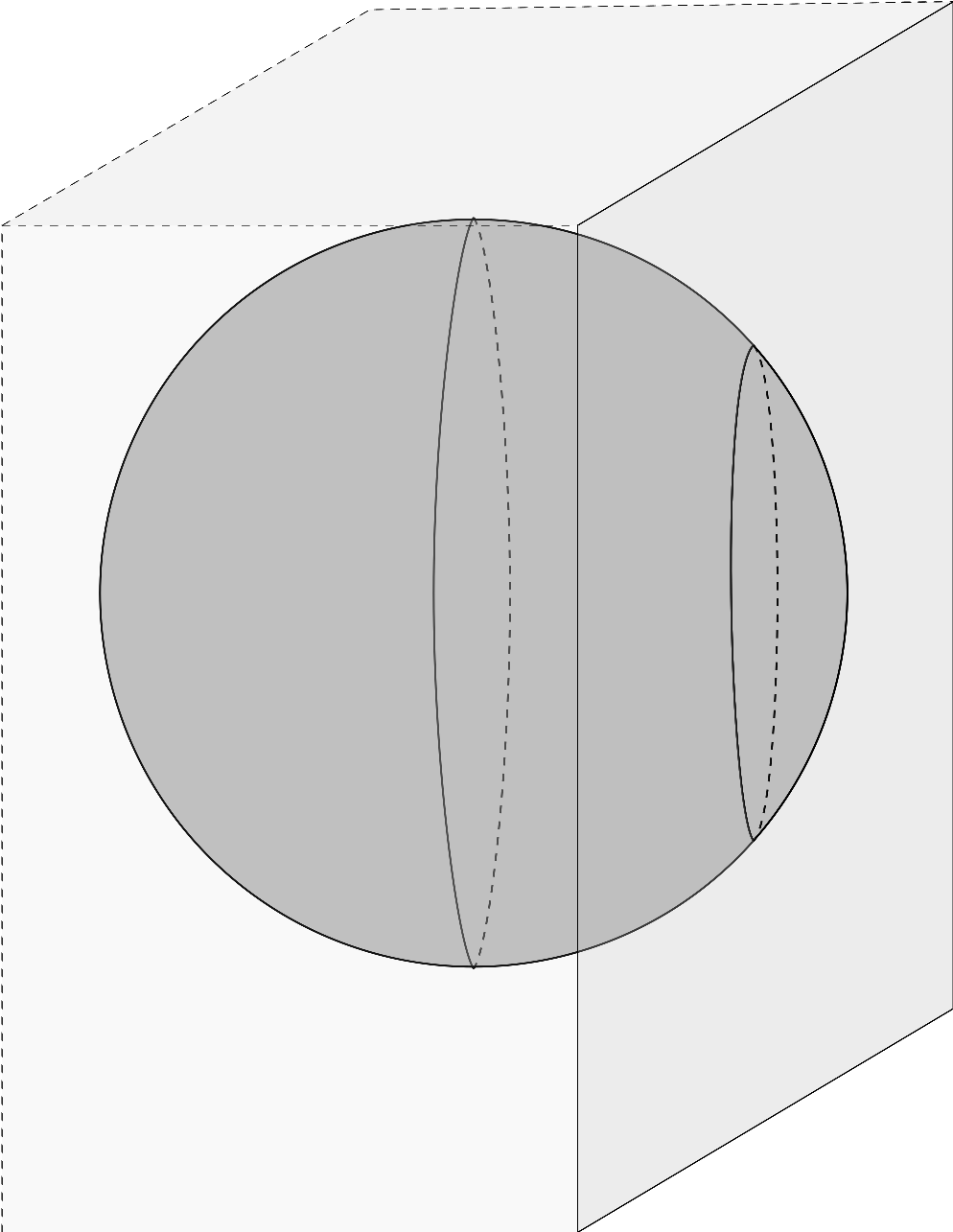}

    (a)
  \end{minipage}
  \begin{minipage}[]{.45\textwidth}
    \centering
  	\includegraphics[width=.6\linewidth]{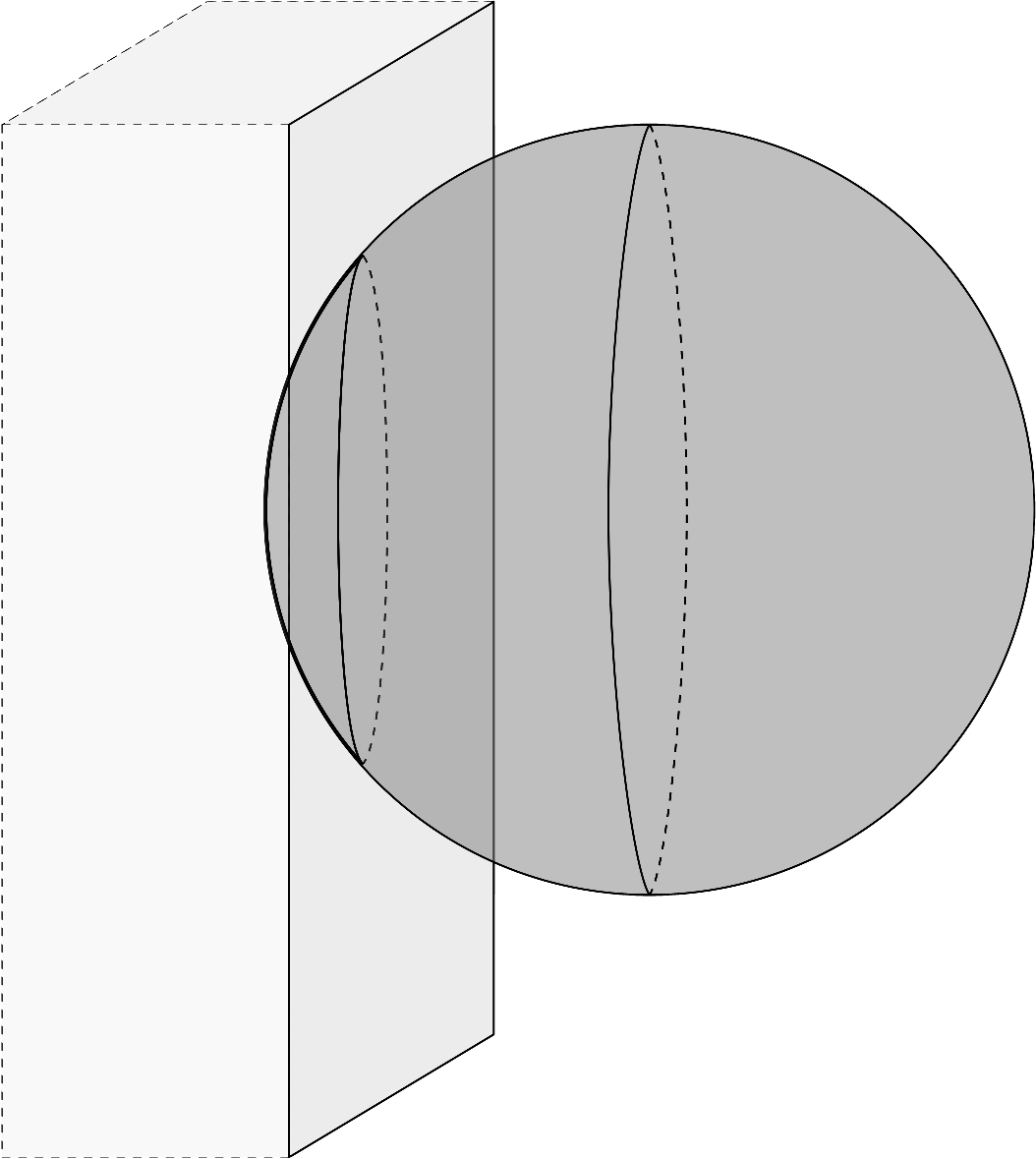}

    (b)
  \end{minipage}
  \caption{(a) $a\leq 0$; (b) $a\in (0,1)$}
  \label{fig:camel}
  \end{figure}
\end{exam}

\begin{exam}\label{ex:camel}
  Set $A = \{z\in\R^{2n}\,|\,z_1< 0\}$ and for $a\in \R$ let $B = \{z\in\R^{2n}\,|\,\|z-(a,0,\dots,0)\|\leq 1\}$ be the ball of radius $1$ around $(a,0,\dots,0)$ (Figure \ref{fig:camel}). For $a \leq 0$ it is obvious that a suitably parametrized translation along $z_1$ is a genuine symplectophagy of $B$ into $A$. For $a\geq 1$ there is no weak symplectophagy of $B$ into $A$  for the same reason as in item (i). For $a\in(0,1)$ we claim that there is no weak symplectophagy of $B$ into $A$, thanks to the symplectic camel theorem. We will use the following formulation of this theorem, implicitly contained in \cite{Viterbo_1992_GenFun}: the space of symplectic embeddings of $B_0:=\Int B$ into $X_b=\R^{2n}\setminus\{z\in\R^{2n}\,|\,z_1=0\,,\|z\|\geq b\}$ is disconnected for $b \in (0,1)$, where a standard embedding into $\{z_1 < 0\}$ and one into $\{z_1 > 0\}$ lie in different connected components. Here by a standard embedding we mean the restriction of a translation. Pick $b \in (\sqrt{1-a^2},1)$ and assume that there is a weak symplectophagy $(\phi_t)_{t\in[0,1]}$ of $B$ into $A$. Let us take a standard embedding of $B_0$ into $\{z_1>0\}$; it is obviously connected to the inclusion $B_0\hookrightarrow X_b$ by a path of embeddings. Concatenating it with the path $\phi_t|_{B_0}$ yields a path of symplectic embeddings starting at a standard one into $\{z_1>0\}$ and ending at $\phi_1|_{B_0}\fc {B_0}\to A = \{z_1<0\}$. This latter embedding is connected to a standard one thanks to Lemma \ref{lem:Alexander_trick} below. Thus we have related two standard embeddings of $B_0$ into $\{\pm z_1>0\}$ by a path, which is a contradiction to the symplectic camel theorem.
\end{exam}

\begin{lemma}\label{lem:Alexander_trick}
  Let $B\subset \R^{2n}$ be the closed unit ball and $\phi \fc B \to H = \{z\in\R^{2n}\,|\,z_1<0\}$ a symplectic embedding.\footnote{In the sense that $\phi$ extends to a symplectic embedding defined on some open set containing $B$.} Then $\phi$ is isotopic to a standard one, via embeddings $B\hookrightarrow H$.
\end{lemma}
\begin{proof}
  This is based on the Alexander trick: for $r\in[0,1]$, define $\psi_r \fc B \to \R^{2n}$ as follows: $\psi_r(x) = r\phi(0)+\frac 1 r(\phi(rx)-\phi(0))$ if $r>0$, and $\psi_0 = d_0\phi|_B$. It is easy to see that $(\psi_r)_{r\in[0,1]}$ is a smooth path of symplectic embeddings $B\hookrightarrow \R^{2n}$ connecting $\psi_0=d_0\phi|_B$ to $\psi_1=\phi$. Concatenating this with a path in the linear symplectic group connecting the differential $d_0\phi$ to the identity, and suitably reparametrizing, we obtain a path $(\iota_t)_{t\in[0,1]}$ of symplectic embeddings $B\hookrightarrow\R^{2n}$ starting at $\phi$ and terminating at the inclusion. Thanks to the compactness of $B$ and the continuity of our path of embeddings, there is $C>0$ large enough so that the path of symplectic embeddings $[0,1]\ni t\mapsto T_{(-Ct,0,\dots,0)}\circ\iota_t$ consists of maps with images in $H$. This path starts at $\phi$ and ends at a standard embedding, as required.
\end{proof}

\begin{exam}\label{ex:weak_no_genuine}
  Here we will construct subsets $A,B\subset\mathbb R^4$ so that there is a weak symplectophagy of $B$ into $A$, but not a genuine one. Let $e_4=(0,0,0,1)$, and let
  \[\overline{B}_x(a)=\{z\in\mathbb R^4\,|\,\pi\|z-x\|^2\leq a\}\]
  be the closed ball of capacity $a$ centered at $x$. Set
  \[R=\frac1{10},\qquad \varepsilon=\frac1{100},\]
  and consider the sets
  $$A = \bigcup_{s\in[0,\frac12]}\ol B_{se_4}(\pi R^2) \cup \bigcup_{s\in[\frac12,1]}\ol B_{se_4}(\pi\epsilon^2)\cup\ol B_{e_4}(\tfrac12\pi R^2)\,,\quad B=\bigcup_{s\in[\frac12,1]}\ol B_{se_4}(\tfrac34\pi R^2)\,,$$
  see Figure \ref{fig:weak_no_genuine}. Recall the translation operators $T_x$ from Example \ref{ex:sympl_embed_polydisk_cube}. Then the isotopy $(T_{-\frac t2e_4})_{t\in[0,1]}$ is a weak symplectophagy of $B$ into $A$, as is easily verified.
 \begin{figure}
    \centering
  	\includegraphics[width=.7\linewidth]{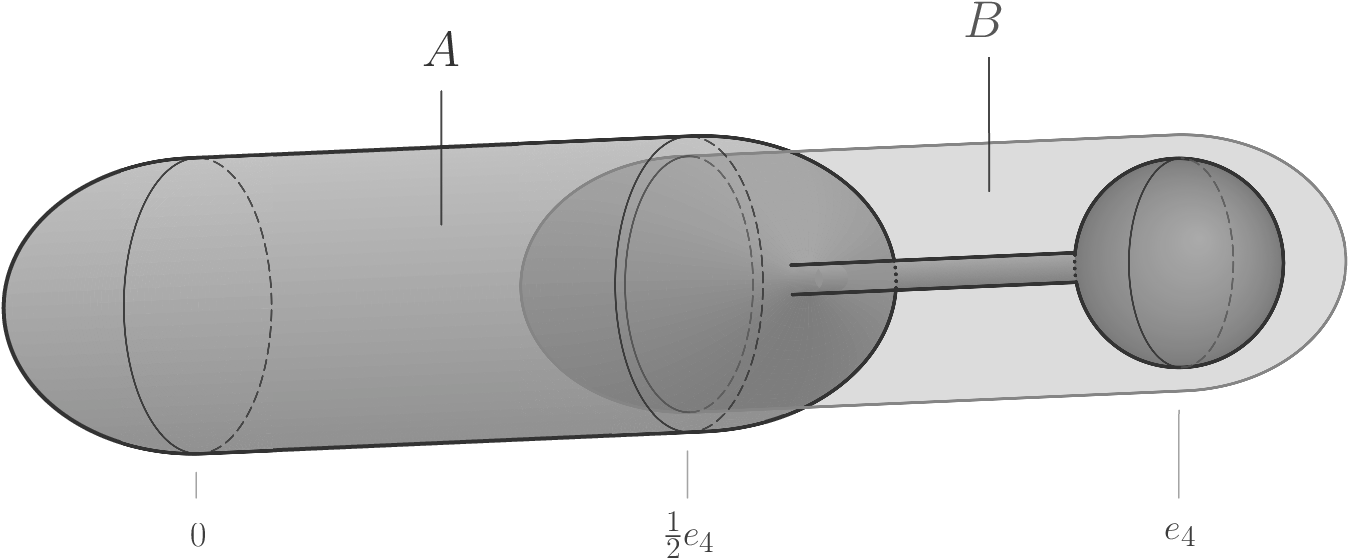}
    \caption{Example \ref{ex:weak_no_genuine}}
    \label{fig:weak_no_genuine}
  \end{figure}

  Let us show that there does not exist a genuine symplectophagy of $B$ into $A$. Arguing by contradiction, we assume that $(\phi_t)_{t\in[0,1]}$ is such a symplectophagy. The contradiction, like in the last example, will arise from the symplectic camel theorem. Let
  \[D=\operatorname{Int}\overline{B}_0(\tfrac12\pi R^2)\]
  and consider the symplectic embeddings
  \[\iota_t \fc D \hookrightarrow B\,,\qquad \iota_t(z)=z+\frac {t+1}2\,e_4 \quad \text{for }t\in[0,1]\,.\]
  Note that $\im\iota_t\subset B$ for all $t \in [0,1]$, and that the endpoints $\iota_0$ and $\iota_1$ both have image in $A\cap B$. Now, consider the following three paths of symplectic embeddings:
  $$\kappa_1=(\phi_t\circ\iota_1)_{t\in[0,1]}\,,\qquad \kappa_2=(\phi_1\circ\iota_{1-t})_{t\in[0,1]}\,,\qquad \kappa_3=(\phi_{1-t}\circ\iota_0)_{t\in[0,1]}\,.$$
  Note that all these paths have images in $A$. For the first and third paths this follows because $\iota_0,\iota_1$ have images in $A\cap B$ and $\phi_t(A\cap B)\subset A$ for all $t$. For the second path this is because $\im\iota_t\subset B$ and $\phi_1(B)\subset A$. Note that the endpoint of $\kappa_1$ equals $\phi_1\circ\iota_1$, which is the beginning point of $\kappa_2$, while the endpoint of $\kappa_2$ is $\phi_1\circ\iota_0$, which is the beginning of $\kappa_3$. It follows that the concatenation $\kappa_1\sharp\kappa_2\sharp\kappa_3$ is a continuous path of symplectic embeddings $D\hookrightarrow A$ beginning at $\iota_1$ and ending at $\iota_0$. Consider
  $$X = \R^4\setminus\{z\in\R^4\,|\,z_4=\tfrac34\,,\|z-\tfrac34e_4\|\geq2\epsilon\}\,,$$
  and note that $A\subset X$. It follows that the above concatenation is a path of symplectic embeddings into $X$ starting at a standard embedding into $\{z_4>\tfrac34\}$ and terminating at a standard embedding into $\{z_4<\tfrac34\}$. Since the hole in $X$ connecting these two half-spaces has radius $2\epsilon = \frac1{50} < \frac3{40}=\frac34\,R$, which is the radius of $D$, we have indeed obtained a contradiction to the symplectic camel theorem, as desired. This contradiction shows that there is no genuine symplectophagy of $B$ into $A$.
\end{exam}

\paragraph*{Organization of the paper} The rest of the paper is organized as follows. In Section \ref{ss:main_results} we formulate our main results, while in Section \ref{ss:obst_rel_SH} we briefly summarize relative symplectic cohomology and its role in our proofs. Section \ref{ss: obst for topophagy} deals with topological obstructions to the existence of certain homotopies, which in particular explain Examples \ref{ex:two_disks}, \ref{ex:sphere}. In Section \ref{s:relSH} we recall the notion of relative symplectic cohomology and some of its properties, and formulate and prove general obstructions to symplectophagies based on these. Section \ref{s:proofs_main_results} contains the proofs of the main results. Finally, Section \ref{s:excision_invariance} is dedicated to the proofs of the technical results---the excision and invariance properties for the relative symplectic cohomology of pairs---of which the latter is a highlight of this paper.

\subsection{Main results}\label{ss:main_results}

Let us now present our results. These consist of Theorems \ref{thm: main thm}, \ref{thm: main thm2}, \ref{thm: main thm3}, which furnish obstructions to weak and genuine symplectophagies in specific settings, and Theorem \ref{thm: general obst for symplectophagy in Cn}, which is a general result relating the existence of symplectophagies in $\C^n$ to the Gromov width and displacement energy of the sets in question. Theorem \ref{thm: main thm} is proved in this section immediately after the formulation of the general Theorem \ref{thm: general obst for symplectophagy in Cn}. The rest of the theorems are proved in Section \ref{s:proofs_main_results}, after we have recalled in Section \ref{s:relSH} the concept of relative symplectic cohomology and its properties.

We will use the standard toric moment map $\Phi_0 \fc \C^n \to \R^n$ given by
$$\Phi_0(z_1,\dots,z_n) = \pi(|z_1|^2,\dots,|z_n|^2)\,.$$
A \tb{toric set} in $\C^n$ for us is any subset of the form $\Phi_0^{-1}(Z)$ for $Z \subset \R^n$. For $a\geq 0$ we have the standard open and closed symplectic balls of capacity $a$:
$$B^{2n}(a) = \{z\in\C^n\,|\,\pi\|z\|^2 < a\}\,,\qquad \ol B{}^{2n}(a) = \{z\in\C^n\,|\,\pi\|z\|^2 \leq a\}\,.$$
These are, of course, toric sets. We also have the closed cylinder of capacity $a$:
$$\ol C{}^{2n}(a) = \ol B{}^2(a) \times \C^{n-1} \subset \C^n\,,\quad\text{which is also toric.}$$

\begin{figure}
  \begin{minipage}[]{.45\textwidth}
    \centering
  	\includegraphics[width=0.7\linewidth]{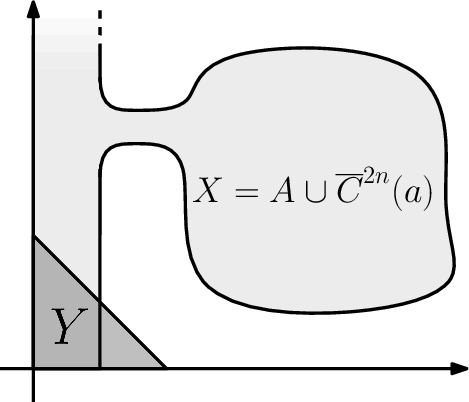}

    (a)
  \end{minipage}
  \begin{minipage}[]{.45\textwidth}
    \centering
  	\includegraphics[width=0.6\linewidth]{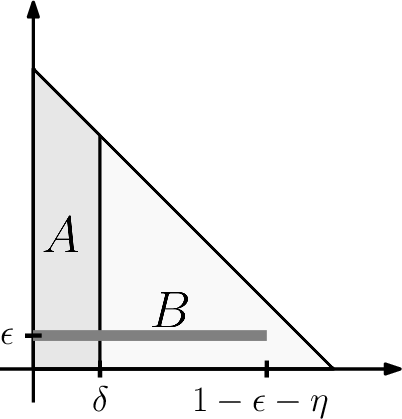}

    (b)
  \end{minipage}
  \caption{(a) Theorem \ref{thm: main thm}; (b) Theorem \ref{thm: main thm2}}
  \label{fig:thms12}
  \end{figure}

\begin{thm}\label{thm: main thm}
    Let $0<a<b$, and let $n\geq 2$. For any compact toric set $A\subset \C^n$ disjoint from $\ol B{}^{2n}(b)$, there is no weak symplectophagy of $Y=\ol B{}^{2n}(b)$ into $X=\ol C{}^{2n}(a)\cup A$.
\end{thm}
\noindent See Figure \ref{fig:thms12}(a) for an example in the moment picture for $n=2$.

\begin{rem}
  There are no topological obstructions here. Indeed, consider the path of radial contractions $R_t:=t\mapsto (1-t)\id$, $t\in [0,1)$. We have $R_t(Y)\subset X\cup Y$ for all $t$, and for $t$ close enough to $1$, $R_t(Y)\subset X$. Also, for $A$ large enough, there is a Hamiltonian isotopy moving $Y$ into $X$, so there is no symplectic embeddability obstruction.
\end{rem}

The next result deals with subsets of $\C P^2$. We will use the standard toric moment map $\Phi \fc \C P^n \to \R^n$ defined by
$$\Phi([z_0:\dots:z_n]) = \|z\|^{-2}(|z_1|^2,\dots,|z_n|^2)\quad\text{ for }z = (z_0,\dots,z_n)\in\C^{n+1}\setminus\{0\}\,.$$
Note that its image is the unit simplex $\big\{x\in\R^n\,|\,\forall i:x_i\geq 0\text{ and }\sum_ix_i\leq 1\big\}$.
\begin{thm}\label{thm: main thm2}
Equip $\C P^2$ with the Fubini--Study symplectic form. For $\epsilon\in(0,1/3)$, $\eta\in(0,\epsilon)$, and $\delta\in(0,1-2\epsilon)$, consider the subsets $A=\Phi^{-1}(A')$, $B=\Phi^{-1}(B')$ of $\C P^2$, where
$$A'=\{x\in \im\Phi\,|\,x_1\leq \delta\},\qquad B'=\{x\in \im\Phi\,|\,x_1\leq 1-\eta-\epsilon\text{ and }x_2=\epsilon\}\,.$$
Then there is no weak symplectophagy of $B$ into $A$.
\end{thm}
\noindent See Figure \ref{fig:thms12}(b) for an illustration.

\begin{rem}\phantom{M}
  \begin{itemize}
    \item There is a smooth isotopy keeping $B$ in $A\cup B$ and mapping $B$ into $A$ at time $1$, meaning there are no topological obstructions. Also, in case $\epsilon\leq\delta$, there is an isotopy composed of unitary maps which at time $1$ maps $B$ into $A$, that is there are no symplectic embeddability obstructions.
    \item It is plausible that a similar result holds in higher dimensions, but our proof is specific to dimension $4$. The reason is that in higher dimensions certain displacement energies are more difficult to bound from above.
  \end{itemize}
\end{rem}

The next result is an obstruction to genuine symplectophagy, and for this we go back to the standard complex linear space. A \tb{toric polyannulus} is a toric subset of $\C^n$ of the form $\Phi_0^{-1}(I_1\times\dots\times I_n)$, where $I_1,\dots,I_n\subset [0,\infty)$ are intervals. It is \tb{nondegenerate} if all the $I_j$ have positive length. A \tb{toric polydisk} is a toric polyannulus $\Phi_0^{-1}(I_1\times\dots\times I_n)$ where $0 \in I_j$ for all $j$.

\begin{figure}
  \begin{minipage}[]{.3\textwidth}
    \centering
  	\includegraphics[width=0.9\linewidth]{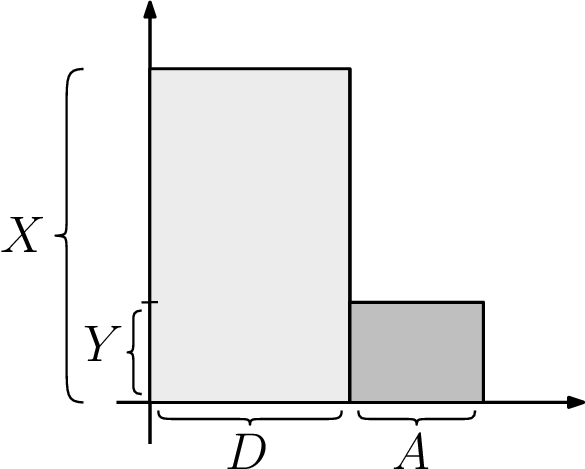}

  (a)
  \end{minipage}
  \begin{minipage}[]{.3\textwidth}
    \centering
    \includegraphics[width=0.9\linewidth]{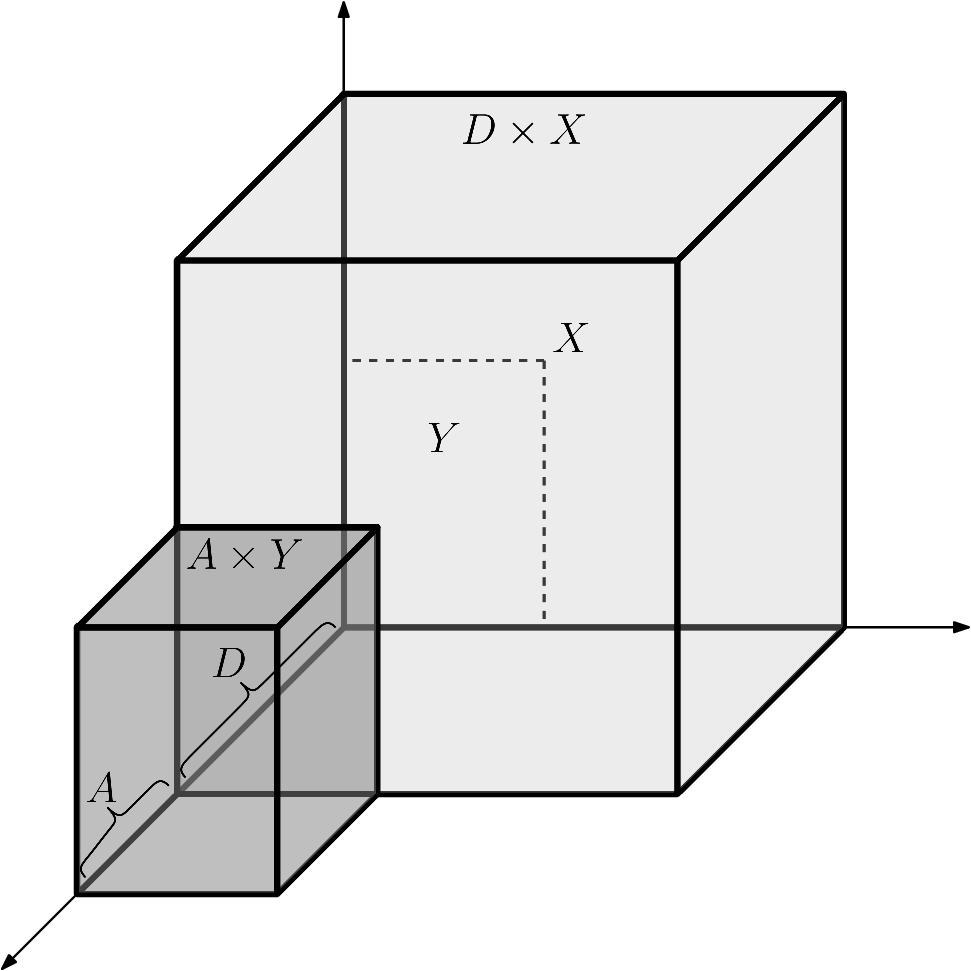}

    (b)
  \end{minipage}
  \begin{minipage}[]{.3\textwidth}
    \centering
    \includegraphics[width=0.9\linewidth]{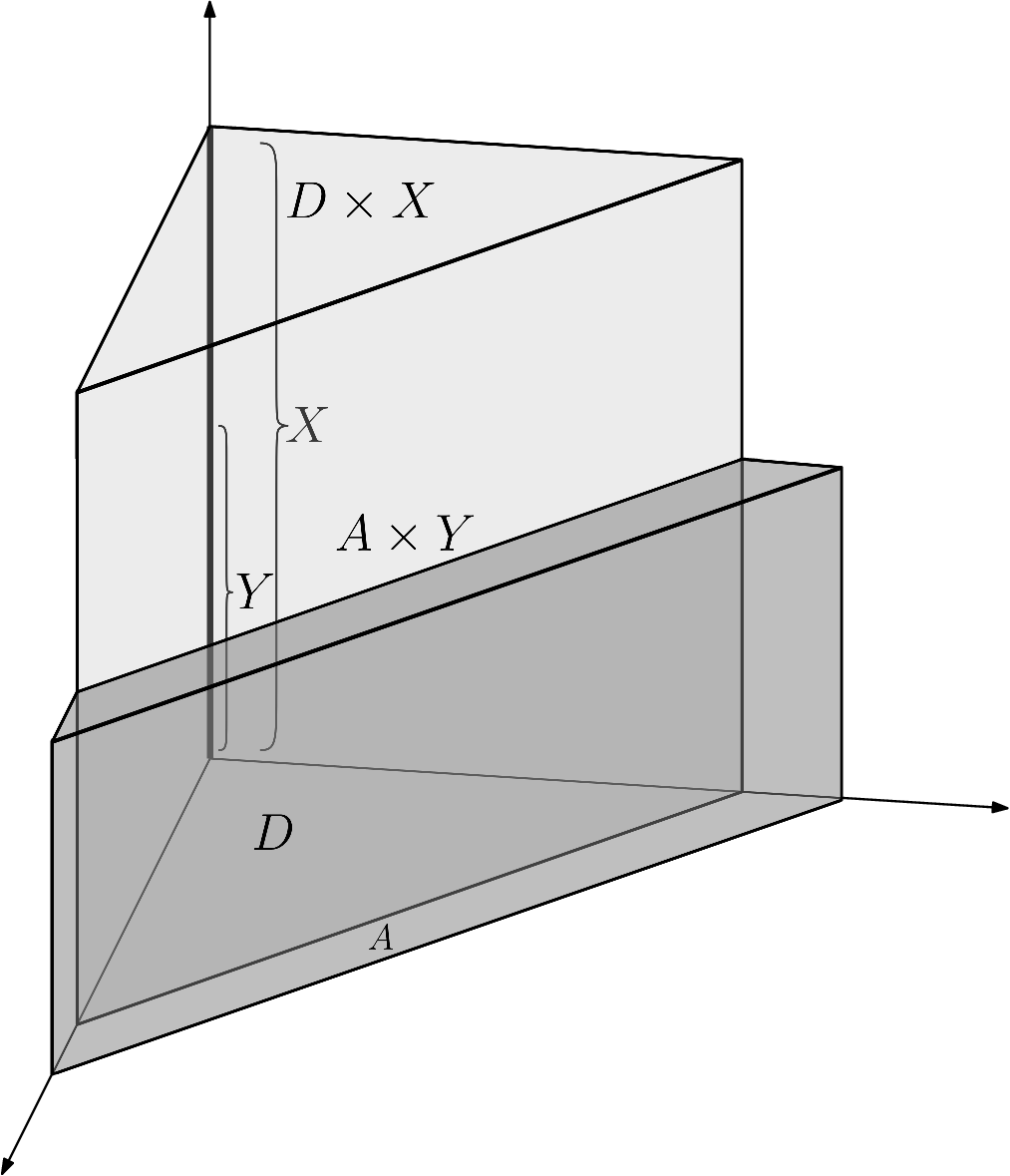}

    (c)
  \end{minipage}
  \caption{Theorem \ref{thm: main thm3}, (a) $k=n=1$; (b) $k=1$, $n=2$; (c) $k=2$, $n=1$}
  \label{fig:thm3}
  \end{figure}

\begin{thm}\label{thm: main thm3}
  Fix positive integers $k,n$, numbers $0 < a < b$, and consider
  $$D = \ol B{}^{2k}(a) = \{z\in \C^k\,|\,\pi\|z\|^2\leq a\} \quad\text{and}\quad A = \{z\in \C^k\,|\,\pi\|z\|^2 \in [a,b]\}\,.$$
  Then for any two nondegenerate closed toric polydisks $X,Y\subset \C^n$, there is no genuine symplectophagy of $A\times Y$ into $D\times X$ in $\C^{n+k}$.
\end{thm}
\noindent See Figure \ref{fig:thm3} for examples in the moment picture.

\begin{rem}\label{rem:rk_main_thm_3}
    Similarly to the first two theorems, there are no topological obstructions. Also, in case $k\leq n$, if $Y$ is small enough and $X$ is large enough, then there is a unitary transformation mapping $A\times Y$ into $D\times X$. See Remark \ref{rem:toric_ex_embeddability_criterion} for a precise statement.
\end{rem}
As mentioned above, these results follow from properties of the relative symplectic cohomology modules and their torsion. However, at this point we can formulate an obstruction to symplectophagy in more familiar terms. In order to do so, we recall here the notion of commutation of two subsets of a symplectic manifold \cite{DGPZ_2024_Symp_top_and_IVMs}, \cite{Varolgunes_2021_MV_and_relSH}.
\begin{defin}\label{def:commuting_sets}
  Let $(M,\omega)$ be a symplectic manifold. A smooth map $F \fc M \to Z$, where $Z$ is a manifold, is \tb{involutive} if for all $f,g \in C^\infty(Z)$, $F^*f,F^*g$ Poisson commute. Two closed subsets $A,B\subset M$ \tb{commute} if there is an involutive map $F\fc M \to Z$ and closed $A',B'\subset Z$ with $A=F^{-1}(A')$, $B=F^{-1}(B')$.
\end{defin}
\begin{exam}\label{ex:commuting_sets}
  The standard toric moment maps $\Phi,\Phi_0$ are involutive. It follows that any two preimages of closed sets by either one of these maps commute. Thus in particular all the sets in Theorems \ref{thm: main thm}, \ref{thm: main thm2}, \ref{thm: main thm3} commute.
\end{exam}
\noindent We will also need two classical quantities associated to a subset of a symplectic manifold, namely its Gromov width and displacement energy. For a subset $A$ of a symplectic manifold $(M,\omega)$ of dimension $2n$ its \tb{Gromov width} is the quantity
$$c_G(A) = \sup\{a\geq 0\,|\,\exists\text{ symplectic embedding }B^{2n}(a)\hookrightarrow \Int(A)\}\,.$$
The \tb{displacement energy} $e(A)$ of $A$ is defined as follows. First, for a Hamiltonian $H \fc [0,1]\times M \to \R$ with compact support consider its Hofer norm
$$\|H\|_{\mathrm{Hofer}} = \int_0^1\Big(\max_MH_t - \min_MH_t\Big)\,dt\,.$$
Then
$$e(A) = \inf\{\|H\|_{\mathrm{Hofer}}\,|\,H\text{ compactly supported}\,,\phi_H(A)\cap\ol A = \varnothing\}\,.$$

\begin{thm}\label{thm: general obst for symplectophagy in Cn}
    Let $X,Y\subset \C^n$ be commuting compact subsets.
    \begin{enumerate}
        \item If there is a weak symplectophagy of $Y$ into $X$ then $c_G(Y)\leq e(X\cap Y)$.
        \item If there is a genuine symplectophagy of $Y$ into $X$ then $c_G(X\cup Y)\leq e(X)$.
    \end{enumerate}
\end{thm}
\noindent Let us see how this easily implies Theorem \ref{thm: main thm}:
\begin{proof}[Proof of Theorem \ref{thm: main thm}]
  Let $(\phi_t)_t$ be a weak symplectophagy of $Y= \ol B{}^{2n}(b)$ into $X = \ol C{}^{2n}(a) \cup A$. Consider the trace of $Y$ under this isotopy:
  $$Y' = \bigcup_{t\in[0,1]}\phi_t(Y)\subset \C^n\,.$$
  Note that it follows from definitions that $Y'\subset X\cup Y$. Since $Y$ is compact, so is $Y'$, and therefore there exists $c>0$ such that $Y'\subset X'\cup Y$, where $X'= \big(\ol B{}^2(a)\times \ol B{}^{2n-2}(c)\big)\cup A$. It follows that $(\phi_t)_t$ is a weak symplectophagy of $Y$ into $X'$. We have $c_G(Y) = b$. On the other hand, $X'\cap Y$ is a compact subset of $\ol C{}^{2n}(a)$, and thus its displacement energy is at most $a$. Thus we have
  $$c_G(Y) = b > a \geq e(X'\cap Y)\,.$$
  Note that the sets $Y,X'$ are compact and that they commute, being both toric. Theorem \ref{thm: general obst for symplectophagy in Cn} then implies that there is no weak symplectophagy of $Y$ into $X'$.
\end{proof}

\subsection{Obstructions from relative symplectic cohomology}\label{ss:obst_rel_SH}

Here we will briefly summarize the tools and methods used in the proofs of our main results. A more detailed account appears in Section \ref{s:relSH}. Our obstructions are derived from the theory of relative symplectic cohomology due to Umut Varolg\"une\c s \cite{Varolgunes_2021_MV_and_relSH}. At a basic level, the relative symplectic cohomology assigns a graded module $SH^*_M(K)$ to each compact subset $K$ of a closed symplectic manifold $(M,\omega)$. This module is over the \tb{(universal) Novikov ring}
$$\Lambda_{\geq0} = \left\{\textstyle \sum_{i=1}^{\infty}a_iT^{\lambda_i}\,\Big|\, \forall i: a_i \in \Q\,,\lambda_i \in [0,\infty)\text{ and }\lambda_i\nearrow \infty\right\}\,.$$
The fundamental piece of information regarding $K$ which we extract from $SH^*_M(K)$ is its \tb{capacity} (see \cite[Remark 4.2.8]{Varolgunes_2018_PhD} and Definition \ref{def:capacity_sets})
$$c(K) = \inf\{\lambda\geq 0\,|\,T^\lambda SH^*_M(K)=0\}\,,$$
which quantifies the torsion of $SH^*_M(K)$.

One of the central features of relative symplectic cohomology is that it is a \emph{functor}, which means that to an inclusion of compact sets $K\subset K'$ there is associated the corresponding \tb{restriction map}
$$\res^{K'}_K \fc SH^*_M(K') \to SH^*_M(K)\,,$$
similarly to singular cohomology. The following theorem expresses the obstructions to symplectophagy in terms of relative symplectic cohomology and the capacity derived from it:
\begin{thm}\label{thm:summary_relSH_obstructions_intro}
  Let $(M,\omega)$ be a closed symplectic manifold and $A,B\subset M$ be compact commuting sets.
  \begin{itemize}
    \item (Proposition \ref{prop: obst for weak symplectophagy, V2}) If there is a weak symplectophagy of $B$ into $A$, then the restriction
    $$SH^*_M(B) \to SH^*(A\cap B)\quad\text{is injective,}$$
    and in particular $c(B)\leq c(A\cap B)$.
    \item (Proposition \ref{prop: SH obstruction for symplectophagy}) If there is a genuine symplectophagy of $B$ into $A$, then the restrictions
    $$SH^*_M(A\cup B) \to SH^*_M(A)\quad\text{and}\quad SH^*_M(B)\to SH^*(A\cap B)\quad\text{are isomorphisms,}$$
    and in particular $c(A\cup B)=c(A)$ and $c(B)=c(A\cap B)$.
  \end{itemize}
\end{thm}

The geometric obstructions of Theorem \ref{thm: general obst for symplectophagy in Cn} are a consequence of this theorem, together with Varolg\"une\c s's energy-capacity inequality (Theorem \ref{thm:energy_capacity_ineq}) relating the capacity of a set to its displacement energy, and the calculations of the relative symplectic cohomology of toric balls in complex projective spaces by Dickstein--Ganor \cite{DG_2026_relSH_in_CPn}. Another crucial ingredient in the proof of Theorem \ref{thm: general obst for symplectophagy in Cn} is the extension of relative symplectic cohomology to \emph{pairs of compact subsets}, defined by the authors together with L.\ Polterovich in \cite{DGPZ_2024_Symp_top_and_IVMs}. The main technical novelty of the present paper is Theorem \ref{thm: invariance}, which expresses an invariance property of the relative symplectic cohomology of pairs under Hamiltonian isotopies.

\subsection{Topological obstructions}\label{ss: obst for topophagy}

We mentioned above that there are purely topological obstructions to symplectophagy. We would like to explain here how these arise from algebraic topology, more specifically from cohomology. Curiously, similar techniques, later applied to relative symplectic cohomology, yield nontrivial obstructions to symplectophagy which are purely symplectic, such as Theorem \ref{thm:summary_relSH_obstructions_intro}.

In the symplectic case, we consider isotopies, which in particular are made of homeomorphisms. However, in order to obtain topological obstructions, we need only consider homotopies. We start with a general obstruction to the existence of certain homotopies. For a homotopy $\phi \fc [0,1]\times X \to Y$, where $X,Y$ are topological spaces, we denote $\phi_t=\phi(t,\cdot) \fc X \to Y$. Also, in this section for such homotopies $X$ will already be a subspace of $Y$, and \emph{we always assume that $\phi_0$ is the inclusion $X\hookrightarrow Y$}. In this section, $H^*(\cdot)$ means singular cohomology with fixed unspecified coefficients.
\begin{prop}\label{prop: obst for weak topophagy, V1}
      Let $X$ be a topological space and $A,B\subset X$. Let $\phi \fc [0,1]\times B \to A\cup B$ be a homotopy with $\phi_1(B)\subset A$. Then the restriction
      $$H^*(A\cup B,A)\to H^*(B,\varnothing)= H^*(B)\quad \text{vanishes.}$$
\end{prop}
\begin{proof}
  From the functoriality and homotopy invariance of cohomology, the diagram
  $$\xymatrix{H^*(A\cup B,A)\ar[r] & H^*(A\cup B)\ar[r]\ar[d] & H^*(A)\ar[d] \\ & H^*(B) & \ar[l]_{(\phi_1)^*} H^*\big(\phi_1(B)\big)}$$
  commutes. Here the unmarked arrows are restrictions. From the exact triangle of the cohomology of the pair $(A\cup B,A)$, the composition of the restriction maps
      $$H^*(A\cup B,A)\to H^*(A\cup B)\to H^*(A)$$
    is zero. Therefore, from the commutativity of the diagram we deduce that the composition
      $$H^*(A\cup B,A)\to H^*(A\cup B)\to H^*(B)$$
  vanishes. Now note that it equals the restriction $H^*(A\cup B,A) \to H^*(B)$.
\end{proof}

As we have mentioned in Example \ref{ex:commuting_sets}, in all of our results concerning obstructions to symplectophagy, the sets in question commute. Commutation implies an algebraic condition called \emph{descent}, see Section \ref{ss:proof_excision} and \cite{Varolgunes_2021_MV_and_relSH}. The algebraic-topological analog of sets in descent turns out to be \emph{excisive couples}, see the discussion in \cite{DGPZ_2024_Symp_top_and_IVMs}. Recall  \cite[p. 188]{Spanier_1966_Alg_Top} that two subsets $A,B\subset X$ form an excisive couple if the inclusion $C_*(A)+C_*(B)\hookrightarrow C_*(A\cup B)$ is a quasi-isomorphism, where $C_*(\cdot)$ is the singular functor. For us the main consequence is that the excision map $H^*(A\cap B,A)\to H^*(B,A\cap B)$ is an isomorphism, see Theorem 4 \emph{ibid.}
\begin{prop}\label{prop: obst for weak topophagy, V2}
    Let $X$ be a topological space, and let $A,B\subset X$ form an excisive couple. If $\phi \fc [0,1]\times B \to A\cup B$ is a homotopy with $\phi_1(B)\subset A$, then the restriction
    \[
    H^*(B)\to H^*(A\cap B)\quad\text{ is injective.}
    \]
\end{prop}
\noindent Let us see an application of this before passing to the proof.
\begin{exam}
  Let $n\geq 2$ and consider
  $$A = \{x\in\R^n\,|\,x_1\geq 0\}\quad\text{and} \quad B=\{x\in \R^n\,|\,\|x\|\in[1,2]\}\,.$$
  Let us show that there is no homotopy of the sort appearing in the proposition. In particular, there is no weak symplectophagy of $B$ into $A$ in case $n$ is even. Indeed, by Proposition \ref{prop: obst for weak topophagy, V2}, were there such a homotopy, the restriction morphism $H^{n-1}(B;\Z) \to H^{n-1}(A\cap B;\Z)$ would be injective. However, $B$ is homotopy equivalent to $S^{n-1}$, thus $H^{n-1}(B;\Z)\cong \Z$, while $A\cap B$ is contractible, meaning $H^{n-1}(A\cap B;\Z) = 0$, and this is a contradiction.
\end{exam}
\begin{proof}[Proof of Proposition \ref{prop: obst for weak topophagy, V2}]
  From the long exact sequences of the pairs $(A\cup B,A)$ and $(B,A\cap B)$ we obtain the diagram
    \[
    \xymatrix{
      H^*(A\cup B,A) \ar[r] \ar[d] & H^*(A\cup B) \ar[r] \ar[d] & H^*(A) \ar[d] \\
      H^*(B,A\cap B) \ar[r]        & H^*(B) \ar[r]              & H^*(A\cap B)
    }
    \]
    which commutes and has exact rows. Proposition~\ref{prop: obst for weak topophagy, V1} implies that the composition
    \[
    H^*(A\cup B,A)\to H^*(A\cup B)\to H^*(B)
    \]
    is zero. Since $A,B$ form an excisive couple, the left vertical arrow
    \[
    H^*(A\cup B,A)\to H^*(B,A\cap B)
    \]
    is an isomorphism. Consequently, the map $H^*(B,A\cap B)\to H^*(B)$ is also zero. Finally, from the exactness of the bottom row, we obtain:
    \[
    \ker\left(H^*(B)\to H^*(A\cap B)\right) = \im\left(H^*(B,A\cap B)\to H^*(B)\right) = 0.
    \]
    Therefore, the map $H^*(B)\to H^*(A\cap B)$ is injective.
\end{proof}

Lastly, we have
\begin{prop}\label{prop: obst_for_genuin_topophagy}
     Let $X$ be a topological space, and let $A,B\subset X$ form an excisive couple. If $\phi \fc [0,1] \times B \to A\cup B$ is a homotopy such that $\phi_1(B)\subset A$ and such that for all $t$ we have $\phi_t(A\cap B) \subset A$, then the restrictions
     $$H^*(A\cup B)\to H^*(A)\quad\text{and}\quad H^*(B)\to H^*(A\cap B)\quad\text{are isomorphisms}.$$
\end{prop}
\noindent This is proved after the following examples.
\begin{exam}\label{ex:top_obstructions}\phantom{M}
  \begin{enumerate}
    \item In the situation of Example \ref{ex:two_disks} in case $A\cap B$ has more than one connected component, we have $H^0(B;\Z)\cong \Z \not \cong H^0(A\cap B;\Z)$, meaning there is no homotopy of the kind appearing in the proposition, and in particular no genuine symplectophagy of $B$ into $A$.
    \item In Example \ref{ex:sphere} we have $H^2(A\cup B;\Z) = H^2(S^2;\Z)\cong \Z\not\cong 0=H^2(A;\Z)$, which again implies that there is no homotopy as above, and thus also no genuine symplectophagy of $B$ into $A$.
  \end{enumerate}
\end{exam}

\begin{proof}[Proof of Proposition~\ref{prop: obst_for_genuin_topophagy}]
By the homotopy invariance of cohomology, we have the following commutative diagram:
\begin{equation}\label{eq:diag_genuine_topophagy}
  \xymatrix{&H^*(A\cup B,A)\ar[dl]\ar[dr] & \\ H^*(B,A\cap B) && \ar[ll]_{(\phi_1)^*} H^*\big(\phi_1(B),\phi_1(A\cap B)\big)}
\end{equation}
Since $\phi_1(A\cap B)\subset \phi_1(B)\subset A$, the restriction
$$H^*(A\cup B,A)\to H^*(\phi_1(B),\phi_1(A\cap B))$$
factorizes through $H^*(A,A)=0$ and therefore it vanishes. Since the diagram~\eqref{eq:diag_genuine_topophagy} commutes, we deduce that the restriction map $$H^*(A\cup B,A)\to H^*(B,A\cap B)$$
is zero as well. Since $A,B$ are an excisive couple, $H^*(A\cup B,A)\to H^*(B,A\cap B)$ is an isomorphism. Since this restriction map is both zero and an isomorphism, it implies that its domain and target are zero:
$$H^*(A\cup B,A)=H^*(B,A\cap B)=0.$$
From the suitable long exact sequences for pairs, we deduce that the restriction maps
    $$H^*(A\cup B)\to H^*(A),\qquad H^*(B)\to H^*(A\cap B)$$
are isomorphisms, as required.
\end{proof}

\subsection{Discussion}\label{ss:discussion}

The notion of a symplectophagy introduced here is new, and thus not directly related to existing results. However, it continues the line of research devoted to various problems of the existence of symplectic or Hamiltonian isotopies satisfying certain conditions. Closest to our results is perhaps the question, motivated by the symplectic camel problem, whether the space of symplectic embeddings of one symplectic manifold into another is connected, in particular whether two such embeddings are related by a symplectic or Hamiltonian isotopy, see for example \cite{Hind_2016_Sympl_isotopies_dim_greater_four} and references therein. Also see \cite{Hind_Kerman_2026_Constraints_Lag_embeddings} for a related problem of mapping Lagrangians into a given domain. The symplectic camel problem is also related to the present paper, as can be seen from Examples \ref{ex:camel}, \ref{ex:weak_no_genuine}.

There is a natural generalization of symplectophagies where the second condition in Definition \ref{def:symplectophagy} is replaced by $\phi_t(B)\subset C$, where $C$ is a set containing $A\cup B$. Although the methods of the present paper can be used to obtain obstructions to the existence of such generalized isotopies, we have not found an example where these obstructions yield something interesting. Additional research is required, and such results will likely be proved by other methods. Another generalization is to consider general symplectic isotopies, not just Hamiltonian ones. It would be interesting to find nontrivial examples of obstructions to such isotopies, since the methods of the present paper rely on the isotopies being Hamiltonian.

The fundamental obstructions to the existence of symplectophagies appearing in Section \ref{ss:obstructions_symplectophagy} are proved using the logic of Section \ref{ss: obst for topophagy}, which deals with topological obstructions, extending the principle whereby methods of algebraic topology are translated into symplectic topology. Here the translation is achieved by replacing the singular cohomology functor by the relative symplectic cohomology, and excision by commutation of sets. This principle was also successfully utilized in \cite{DGPZ_2024_Symp_top_and_IVMs} to obtain interesting new symplectic rigidity results. One notable difference between that work and the present paper is that here we use the \emph{torsion} part of relative symplectic cohomology.

It is also interesting whether the use of additional structure on relative symplectic cohomology, such as products (\cite{Tonkonog_Varolgunes_2023_Super_rigidity_of_certain_skeleta_using_relSH}, \cite{DGPZ_2024_Symp_top_and_IVMs}) can furnish sharper results. Finally, in the present paper we only establish obstructions to symplectophagies when the sets in question commute, which allows us to use the power of the Mayer--Vietoris sequence for such sets, due to Varolg\"une\c s \cite{Varolgunes_2021_MV_and_relSH}. It is intriguing whether obstructions to symplectophagy for non-commuting sets can be obtained which do not rely either on existing results regarding spaces of symplectic embeddings, such as those in Examples \ref{ex:camel}, \ref{ex:weak_no_genuine}.

As a curious historical note, we would like to mention that the notion of a symplectophagy grew out of the attempt to understand the sequence
$$SH^*(A\cup B) \xrightarrow{(\res^{A\cup B}_A,\res^{A\cup B}_B)} SH^*(A)\oplus SH^*(B) \xrightarrow{\res^A_{A\cap B}-\res^B_{A\cap B}} SH^*(A\cap B)$$
in case $A,B$ are compact subsets of a symplectic manifold. It follows from the functoriality of the restriction maps (Sections \ref{ss:obst_rel_SH}, \ref{ss:single_sets}) that the above sequence is a chain complex, in the sense that the composition of the two maps vanishes.\footnote{We owe this observation to L.\ Polterovich.} The Mayer--Vietoris property of relative symplectic cohomology, proved by Varolg\"une\c s \cite{Varolgunes_2021_MV_and_relSH}, states that this sequence is exact if $A,B$ commute. In general, the cohomology of the above chain complex at the middle term can serve as a measure of failure of commutation, and thus it is interesting to be able to compute it. It turns out that in case there is a genuine symplectophagy of $B$ into $A$, this cohomology is isomorphic to $\coker \res^{A\cup B}_A\oplus \ker \res^B_{A\cap B}$.

\begin{acknow}We wish to express our gratitude to Dylan Cant, Octav Cornea, Ely Kerman, Emmanuel Opshtein, Felix Schlenk, and Egor Shelukhin for useful discussions, as well as to Leonid Polterovich for enlightening conversations during the early stages of the project, and for his interest and encouragement. We are also grateful to Fabian Ziltener for a question regarding Theorem \ref{thm: main thm3}, which led us to formulate and prove a more general statement.

The first author is supported by the Jenny Panitch Beckow Memorial Scholarship Canada, 2026. He was previously partially supported by the Milner Foundation and by the Israel Science Foundation (grant no.\ 1102/20) during his PhD studies at Tel Aviv University. He is currently partially supported by Fondation Courtois during his postdoctoral fellowship at the Centre de Recherches Math\'ematiques, Universit\'e de Montr\'eal. He would also like to thank the Centre de Recherches Math\'ematiques for its warm hospitality and excellent working environment.
\end{acknow}

\medskip
\noindent\textbf{Use of AI tools.}
Throughout this project, ChatGPT was used to assist in the
investigation of several mathematical questions, as detailed below.

\begin{itemize}

\item \textbf{Example~\ref{ex:weak_no_genuine}.}
Discussions with ChatGPT (GPT-5.6 Sol), guided by our geometric intuition, helped formulate the shape of the set $A$ in this example. Subsequently, a more advanced model (GPT-6 Astra) suggested the reduction to the symplectic camel problem used in the proof, building on our understanding that the camel obstruction was the key to ruling out a genuine symplectophagy.

\item \textbf{Theorem~\ref{thm: main thm}.}
This theorem is a far-reaching generalization, developed by us, of a specific four-dimensional example. The example emerged following our examination of
a TikZ figure generated by ChatGPT (GPT-5.6 Sol), during a search for subsets of $\C^2$ defined via the standard moment map that violate the equality $c(B) = c(A\cap B)$.

\item \textbf{First assertion of
Theorem~\ref{thm:summary_relSH_obstructions_intro}.}
The observation that the existence of a weak symplectophagy implies injectivity of the relevant restriction map originated from an AI-assisted search (using GPT-5.6 Sol) for applications of Proposition~\ref{prop: obst for weak symplectophagy, V1}. There, the AI ``cheated'' in achieving its assigned goal by utilizing consequences of commutation of the sets, which is not an assumption of Proposition~\ref{prop: obst for weak symplectophagy, V1}. We subsequently refined the observation into its present form.

\item \textbf{Theorem~\ref{thm:exists_parametrized_ray}.}
The development of the original proof was accelerated by AI assistance in completing technical details. The proof presented in the paper is a revised and reworked version by us.
\end{itemize}

We have carefully reviewed and fully understand the mathematical arguments as they appear in the paper, and take complete responsibility for their correctness. ChatGPT (GPT-6 Astra) was also used for proofreading and for checking the mathematical correctness of the final version of the manuscript.

\section{Relative symplectic cohomology}\label{s:relSH}

Here we briefly review Varolg\"une\c s's relative symplectic cohomology for compact sets \cite{Varolgunes_2018_PhD,Varolgunes_2021_MV_and_relSH} and its extension to pairs of compact sets, defined in \cite{DGPZ_2024_Symp_top_and_IVMs}. More details of these constructions appear in Section \ref{ss:basic_notions_defs}. In Section \ref{ss:obstructions_symplectophagy} we formulate and prove fairly general obstructions to the existence of symplectophagies in terms of relative symplectic cohomology. These obstructions are particularly potent when the sets in question commute.

The \tb{(universal) Novikov field} is
$$\Lambda = \left\{\textstyle \sum_{i=1}^{\infty}a_iT^{\lambda_i}\,\Big|\, \forall i: a_i \in \Q\,,\lambda_i \in \R\text{ and }\lambda_i\nearrow \infty\right\}\,.$$
It carries the \tb{valuation}
$$\val \fc \Lambda \to \R\cup\{\infty\}\,,$$
where $\val(0) = \infty$ and $\val\left(\sum_{i=1}^{\infty}a_iT^{\lambda_i}\right) = \lambda_1$ if $a_1\neq 0$ and $\lambda_1<\lambda_2<\dots$. The \tb{(universal) Novikov ring} is the subring $\Lambda_{\geq 0} = \val^{-1}\big([0,\infty]\big)$. More generally, we have the \tb{interval modules}
$$\Lambda_{>r} = \val^{-1}\big((r,\infty]\big)\,,\quad \Lambda_{\geq r} = \val^{-1}\big([r,\infty]\big)\,,\quad\text{and }\Lambda_{(r,s]} =\Lambda_{>r}/\Lambda_{>s}\text{ for }r\leq s\,.$$
Note that $\Lambda_{\geq 0}$ is a valuation ring and that $\Lambda_{>0}$ is its unique maximal ideal.

We refer to modules over the Novikov ring $\Lambda_{\geq 0}$ as \tb{Novikov modules}. Typical examples are the Novikov field, the Novikov ring itself, and the interval modules. A Novikov module $A$ is \tb{$k$-graded} for $k\in \N$ if $A = \bigoplus_{i\in\Z/(k)}A^i$, where $A^i$ are the graded components. When $k$ is clear from the context, we will simply write ``graded modules.'' Note that when $k=0$, we have the usual notion of ($\Z$-)graded modules. If $A,B$ are $k$-graded modules, a module map $f \fc A \to B$ is \tb{graded of degree $d \in \Z/(k)$} if $f(A^i)\subset B^{i+d}$ for all $i$. Graded modules form a category where the morphisms are graded maps of any degree; it has a natural subcategory whose morphisms are graded of degree zero.

\subsection{Single sets}\label{ss:single_sets}

Fix a closed symplectic manifold $(M,\omega)$, and let $N_M$ be the minimal Chern number of $M$, defined as the unique nonnegative generator\footnote{In case $c_1(M)$ vanishes on spheres, this means that $N_M = 0$, which differs from the convention, sometimes appearing in the literature, which sets $N_M=\infty$ in this case.} of the subgroup $\langle c_1(M),\pi_2(M)\rangle\subset \Z$. Let $\cC(M)$ be the category whose objects are the compact subsets of $M$ and whose morphisms are the inclusions of compact sets. The relative symplectic cohomology (in $M$) is a contravariant functor $SH^*_M$ from $\cC(M)$ to the category of $2N_M$-graded Novikov modules and degree zero graded module maps. If $K\subset M$ is compact, its relative symplectic cohomology is denoted by $SH^*_M(K)$. Since $SH^*_M$ is a functor, if $K'\subset K$ is another compact subset, we have the corresponding restriction map
$$\res^K_{K'} \fc SH^*_M(K) \to SH^*_M(K')\,,$$
which is a degree zero graded module map. Moreover, $\res^K_K = \id$ and $\res^{K'}_{K''}\circ\res^K_{K'} = \res^K_{K''}$ if $K''\subset K'$. As the context will usually determine $M$ unambiguously, we will mostly drop the subscript $_M$ from $SH^*_M$.

The relative symplectic cohomology satisfies various useful properties, some of which we now list. There are other properties such as a product (\cite{Tonkonog_Varolgunes_2023_Super_rigidity_of_certain_skeleta_using_relSH}, \cite{DGPZ_2024_Symp_top_and_IVMs}), and a chain-level E2-operad structure \cite{AGV_2024_framed_E2_structures}, but we won't use these.
\begin{itemize}
  \item $SH^*(M)\cong H^*(M;\Z)\otimes_\Z\Lambda_{>0}$ and $SH^*(\varnothing) = 0$;
  \item For commuting compact $A,B\subset M$ we have the long exact \tb{Mayer--Vietoris sequence}
  $$\dots \to SH^{*-1}(A\cap B)\xrightarrow{\delta^*} SH^*(A\cup B) \to SH^*(A)\oplus SH^*(B) \to SH^*(A\cap B) \xrightarrow{\delta^*}\dots\,,$$
  where $\delta^*$ is a certain degree $1$ graded module morphism while the rest of the arrows come from restrictions.
  \item The relative symplectic cohomology is natural in the sense of it being a symplectic invariant. Concretely it means that given a symplectomorphism $\phi$ of $M$, there is an isomorphism $\phi_* \fc SH^*(K) \to SH^*(\phi(K))$ for each $K\in\cC(M)$, which commutes with restrictions. We refer to $\phi_*$ as \tb{the naturality isomorphism} induced by $\phi$; see Section \ref{sss:naturality} for its construction.
  \item We have the following invariance property, crucial for applications: if $K\subset K'\subset M$ are compact and $(\phi_t)_{t\in[0,1]}$ is a Hamiltonian isotopy of $M$ with the property that for all $t$ we have $\phi_t(K)\subset K'$, then the following diagram commutes:
      \begin{equation}\label{eq:diagram_invariance_single_sets}
        \xymatrix{& SH^*(K')\ar[ld]_{\res} \ar[rd]^{\res} & \\ SH^*(K) \ar[rr]^{(\phi_1)_*} && SH^*\big(\phi_1(K)\big)}
      \end{equation}
      A proof of this property was sketched in Varolg\"une\c s's thesis \cite{Varolgunes_2018_PhD}. It can also be obtained as a particular case of our invariance result, Theorem \ref{thm: invariance}.
  \item If $K$ is displaceable, then $SH^*(K)$ is a torsion Novikov module, in the sense that each element of it is killed by some power $T^\lambda$, $\lambda \geq 0$. Equivalently
      $$SH^*(K)\otimes_{\Lambda_{\geq0}}\Lambda = 0\,.$$
\end{itemize}
\noindent The last property has a quantitative version, taking the form of an energy-capacity inequality (Theorem \ref{thm:energy_capacity_ineq} below). To formulate it, recall that Varolg\"une\c s \cite[Remark 4.2.8]{Varolgunes_2018_PhD} defined the \tb{capacity} of a Novikov module $A$ to be
$$c(A) = \inf\{\lambda \geq 0\,|\,T^\lambda A=0\}\,.$$
\begin{exam}
  For $s\geq r$, $c(\Lambda_{(r,s]})=s-r$.
\end{exam}
This induces the following capacity-like quantity for subsets of a symplectic manifold.
\begin{defin}\label{def:capacity_sets}
  Given a closed symplectic manifold $(M,\omega)$, the \tb{relative symplectic capacity} is the function $c \fc \cC(M) \to [0,\infty]$ defined by
  $$c(K):=c\big(SH^*(K)\big)\,.$$
\end{defin}
\begin{rem}
  Since the relative symplectic capacity is the only capacity-like quantity used in the paper, besides Gromov's width, we will simply refer to it as the \tb{capacity}. Also, by abuse of notation, we denote the capacity of a Novikov module and the capacity of a compact subset by the same letter $c$, but we trust that the context will resolve this ambiguity in each instance.
\end{rem}

\noindent We then have the following energy-capacity inequality.
\begin{thm}[{\cite[Remark 4.2.8]{Varolgunes_2018_PhD}}]\label{thm:energy_capacity_ineq}
  Let $K\subset M$ be compact. Then $c(K) \leq e(K)$. \qed
\end{thm}
\begin{rem}
  Dickstein--Ganor \cite{DG_2026_relSH_in_CPn} improve this to $c(K) \leq e_{\mathrm{st}}(K)$, where $e_{\mathrm{st}}$ is the \emph{stable} displacement energy.
\end{rem}

We have the following easy consequence of the Mayer--Vietoris sequence.
\begin{prop}\label{prop:capacity_estimate_commuting_sets}
  If $A,B\subset M$ are commuting compact subsets, then
  $$c(A\cup B) \leq \max\{c(A),c(B)\} + c(A\cap B)\,.$$
\end{prop}
\noindent For the proof we need the following lemma stating some elementary properties of the Novikov module capacity, which are straightforward consequences of the definitions and whose proofs are left to the reader.
\begin{lemma}\label{lem:elementary_properties_capacity}\phantom{M}
  \begin{itemize}
    \item We have $c\big(\bigoplus_{i\in I}A_i\big) = \sup_{i\in I}c(A_i)$ for a family $(A_i)_{i\in I}$ of Novikov modules; in particular the capacity of a graded module is the supremum of the capacities of its graded components;
    \item If $A \to B \to C$ is an exact sequence, then $c(B) \leq c(A) + c(C)$.\qed
  \end{itemize}
\end{lemma}

\begin{proof}[Proof of Proposition \ref{prop:capacity_estimate_commuting_sets}]
  For $j \in \Z/(2N_M)$, consider the following piece of the Mayer--Vietoris sequence for $A,B$:
  $$SH^{j-1}(A\cap B) \to SH^j(A\cup B) \to SH^j(A)\oplus SH^j(B)\,.$$
  Lemma \ref{lem:elementary_properties_capacity} implies that
  \begin{multline*}
    c(SH^j(A\cup B)) \leq c\big( SH^j(A)\oplus SH^j(B)\big) + c(SH^{j-1}(A\cap B))\\
       = \max\{c(SH^j(A)),c(SH^j(B))\} + c(SH^{j-1}(A\cap B))\,.
  \end{multline*}
  Again by Lemma \ref{lem:elementary_properties_capacity} we have $c(K) = \sup_jc\big(SH^j(K)\big)$. Thus the right-hand side in this inequality is at most $\max\{c(A),c(B)\} + c(A\cap B)$. Now taking the supremum over $j$ of the left-hand side yields the desired inequality.
\end{proof}

The following proposition, proved in Section \ref{ss: proof of monotonicity}, expresses a monotonicity property for the capacity of subsets of complex projective spaces, which we will need for the proof of Theorem \ref{thm: general obst for symplectophagy in Cn}. It is based on the calculation of relative symplectic cohomology for balls in complex projective spaces and relevant restriction maps in \cite{DG_2026_relSH_in_CPn}. Let us call a \tb{closed toric ball} any subset of $\C P^n$ of the form $\Phi^{-1}(\Delta)$, where $\Delta = \{x\in [0,\infty)^n\,|\,\sum_{j=1}^nx_j\leq a\}$ for $a \in (0,1)$.
\begin{prop}\label{prop: weak monotonicity of c_relSH}
  Let $B\subset K \subset \C P^n$, where $B$ is a closed toric ball while $K$ is compact. Then $c(B) \leq c(K)$.
\end{prop}
\begin{rem}
  One of the goals of \cite{DGZ_Var_conj} is to show that the relative symplectic capacity is monotone with respect to inclusions. This, of course, implies Proposition \ref{prop: weak monotonicity of c_relSH}. However, the proof requires the use of the product on relative symplectic cohomology together with an additional structure. We will not use this in the present paper.
\end{rem}

\subsection{Pairs}\label{ss:pairs}

In \cite{DGPZ_2024_Symp_top_and_IVMs}, the authors, together with L.\ Polterovich, extended the definition of relative symplectic cohomology to \emph{pairs of compact subsets}, which we will now briefly recall. Let $\cC\cP(M)$ denote the category whose objects are pairs $(K,K')$, where $K,K'\subset M$ are compact and $K'\subset K$, and where the only morphisms are inclusions of pairs $(K,K')\subset (L,L')$, that is $K\subset L$, $K'\subset L'$. In \cite{DGPZ_2024_Symp_top_and_IVMs} we defined a relative symplectic cohomology of pairs as any contravariant functor from $\cC\cP(M)$ to the category of $2N_M$-graded Novikov modules and degree zero morphisms satisfying a certain natural list of properties, and proved that such a functor indeed exists. In the present paper, by \tb{the relative symplectic cohomology of pairs} we mean the specific functor constructed in \cite{DGPZ_2024_Symp_top_and_IVMs}, see also Section \ref{sss:relSH_pairs}. The reason is that that functor satisfies an additional property not listed among the axioms in \cite{DGPZ_2024_Symp_top_and_IVMs}, namely invariance with respect to Hamiltonian isotopies, similar to diagram \eqref{eq:diagram_invariance_single_sets}, see Theorem \ref{thm: invariance}.

The relative symplectic cohomology of pairs, being a functor, associates to each compact pair $(K,K')\in\cC\cP(M)$ a $2N_M$-graded Novikov module $SH^*_M(K,K')$, and to each inclusion of pairs $(K,K') \subset (L,L')$ the corresponding restriction map
$$\res^{(L,L')}_{(K,K')} \fc SH^*_M(L,L') \to SH^*_M(K,K')\,,$$
which is a degree zero graded module morphism. These are subject to the usual functoriality identities. Moreover, there is a natural isomorphism $SH^*_M(K,\varnothing)\cong SH^*_M(K)$ compatible with restrictions. We use the notation $\res^{(L,L')}_K \fc SH^*_M(L,L') \to SH^*_M(K)$ to mean the composition of this natural isomorphism with the restriction $\res^{(L,L')}_{(K,\varnothing)}$.

From now on, we will mostly drop the subscript $_M$ also from the relative symplectic cohomology of pairs. For each compact pair $(K,K')$ we have the corresponding long exact sequence
$$\dots \to SH^{*-1}(K') \xrightarrow{\star} SH^*(K,K') \to SH^*(K) \to SH^*(K')\xrightarrow{\star} \dots\,,$$
where $\star$ denotes the connectant map, the rest of the arrows being restrictions.
\begin{rem}\label{rem:SH_K_K_vanishes}
  Note that from this exact sequence one obtains that $SH^*(K,K)=0$.
\end{rem}

The following two properties of the relative symplectic cohomology of pairs---\emph{excision} and \emph{invariance}---are paramount for our applications.
\begin{prop}[Excision for relative symplectic cohomology]\label{prop:excision}
  Let $A,B\subset M$ be commuting compact subsets. Then the restriction map
  $$SH^*(A\cup B,A) \to SH^*(B,A\cap B) \quad \text{is an isomorphism.}$$
\end{prop}

To formulate the invariance property, we mention here that, just like for single sets, a symplectomorphism $\phi$ of $M$ induces \tb{naturality isomorphisms}
$$\phi_* \fc SH^*(K,K') \to SH^*(\phi(K),\phi(K')) \quad \text{for compact pairs }(K,K')\,.$$
These are defined in Section \ref{sss:naturality}.
\begin{thm}[Invariance for relative SH of pairs under isotopies]  \label{thm: invariance}
  Let $(M,\omega)$ be a closed symplectic manifold, let $(K,K')\subset (L,L')$ be compact pairs in $M$, and assume that $(\phi_t)_{t\in[0,1]}$ is a Hamiltonian isotopy of $M$ such that for all $t$ we have $\big(\phi_t(K),\phi_t(K')\big)\subset (L,L')$. Then with $\phi=\phi_1$, the following diagram commutes:
  $$\xymatrix{& SH^*(L,L') \ar[ld]_{\res} \ar[rd]^{\res} & \\ SH^*(K,K')\ar[rr]^{\phi_*} && SH^*\big(\phi(K),\phi(K')\big)}$$
\end{thm}
\noindent Proposition \ref{prop:excision} is proved in Section \ref{ss:proof_excision} and Theorem \ref{thm: invariance} in Section \ref{ss:proof_invariance}.
\begin{rem}
  Note that if we put $L'=K'=\varnothing$ in the theorem, we recover the aforementioned invariance property for single sets.
\end{rem}

\subsection{Obstructions to symplectophagy}\label{ss:obstructions_symplectophagy}

Here we formulate obstructions to symplectophagy, based on the properties of the relative symplectic capacity and the relative symplectic cohomology of pairs. In this section, $(M,\omega)$ is a fixed closed symplectic manifold. We wish to draw the reader's attention to the fact that the proofs of the three propositions in this section parallel, logically and ideologically, those appearing in Section \ref{ss: obst for topophagy}.

\begin{prop}\label{prop: obst for weak symplectophagy, V1}
  Let $A,B\subset M$ be compact. If there is a weak symplectophagy of $B$ into $A$, then the restriction map
  $$SH^*(A\cup B,A) \to SH^*(B)\quad\text{vanishes.}$$
\end{prop}

\begin{proof}
  Let $(\phi_t)_{t\in[0,1]}$ be a weak symplectophagy of $B$ into $A$. In the following diagram the upper triangle commutes by functoriality, while the lower one by the invariance of relative symplectic cohomology (Theorem~\ref{thm: invariance}, with $K=B$, $L=A\cup B$, $K'=L'=\varnothing$):
  \[
  \xymatrix{
    SH^*(A\cup B,A)\ar[r] & SH^*(A\cup B)\ar[r]\ar[d] \ar[dr]& SH^*(A)\ar[d] \\
    & SH^*(B)\ar[r]^-{\phi_*} & SH^*\big(\phi(B)\big)
  }
  \]
  Here the unmarked arrows are restrictions, and $\phi_*$ is the naturality isomorphism induced by $\phi=\phi_1$.
  From the exact triangle of the relative symplectic cohomology of the
  pair $(A\cup B,A)$, the composition of the restriction maps
  \[
    SH^*(A\cup B,A)\to SH^*(A\cup B)\to SH^*(A)
  \]
  is zero. Therefore, from the commutativity of the diagram we deduce
  that the composition
  \[
    SH^*(A\cup B,A)\to SH^*(A\cup B)\to SH^*(B)
  \]
  vanishes. Now note that it equals the restriction
  $SH^*(A\cup B,A)\to SH^*(B)$.
\end{proof}

\begin{prop}\label{prop: obst for weak symplectophagy, V2}
  If $A,B\subset M$ are compact commuting subsets and there is a weak symplectophagy of $B$ into $A$, then the restriction map
  $$SH^*(B) \to SH^*(A\cap B) \quad \text{is injective,}$$
  and in particular $c(B) \leq c(A\cap B)$.
\end{prop}

\begin{proof}
  From the long exact sequences of the pairs $(A\cup B,A)$, $(B,A\cap B)$, and functoriality we obtain the diagram
  \[
  \xymatrix{
    SH^*(A\cup B,A)\ar[r]\ar[d] & SH^*(A\cup B)\ar[r]\ar[d] & SH^*(A)\ar[d] \\
    SH^*(B,A\cap B)\ar[r] & SH^*(B)\ar[r] & SH^*(A\cap B)
  }
  \]
  which commutes and has exact rows.
  Proposition~\ref{prop: obst for weak symplectophagy, V1} implies
  that the composition
  \[
    SH^*(A\cup B,A)\to SH^*(A\cup B)\to SH^*(B)
  \]
  is zero. Since $A,B$ commute, the left vertical arrow
  \[
    SH^*(A\cup B,A)\to SH^*(B,A\cap B)
  \]
  is an isomorphism by excision (Proposition~\ref{prop:excision}).
  Consequently, the restriction map $SH^*(B,A\cap B)\to SH^*(B)$ also vanishes.
  Finally, from the exactness of the bottom row, we obtain:
  \[
    \ker\left(SH^*(B)\to SH^*(A\cap B)\right)
    =\im\left(SH^*(B,A\cap B)\to SH^*(B)\right)=0.
  \]
  Therefore, the map $SH^*(B)\to SH^*(A\cap B)$ is injective.
  The inequality $c(B)\leq c(A\cap B)$ follows from
  Lemma~\ref{lem:elementary_properties_capacity}, applied to the
  exact sequence $0\to SH^*(B)\to SH^*(A\cap B)$.
\end{proof}
%
%
%

\begin{prop}\label{prop: SH obstruction for symplectophagy}
  If $A,B\subset M$ are commuting compact subsets and there is a genuine symplectophagy of $B$ into $A$, then the restriction maps
  $$SH^*(A\cup B) \to SH^*(A)\quad\text{and}\quad SH^*(B)\to SH^*(A\cap B)$$
  are isomorphisms. In particular, $c(A\cup B) = c(A)$ and $c(B) = c(A\cap B)$.
\end{prop}

\begin{proof}
Assume that there exists a genuine symplectophagy $(\phi_t)_{t\in [0,1]}$ of $B$ into $A$, and let $\phi=\phi_1$. By Theorem~\ref{thm: invariance} the diagram
 \begin{equation}\label{diag: triangle of symplectophagy}
     \xymatrix{
    & SH^*(A\cup B,A) \ar[dl]_{\res} \ar[dr]^{\res} \\
    SH^*(B,A\cap B)\ar[rr]^{\phi_*} & & SH^*(\phi(B),\phi(A\cap B))
}
 \end{equation}
commutes.
Since $\phi(A\cap B)\subset \phi(B)\subset A$, the restriction map $$\res\fc SH^*(A\cup B,A)\to SH^*(\phi(B),\phi(A\cap B))$$
factorizes through $SH^*(A,A)=0$ (Remark \ref{rem:SH_K_K_vanishes}), hence this map vanishes. Since the diagram~\eqref{diag: triangle of symplectophagy} commutes, we deduce that the restriction map
$$\res\fc SH^*(A\cup B,A)\to SH^*(B,A\cap B)$$
is zero as well. Since $A,B$ commute, we deduce that by the excision isomorphism, Proposition~\ref{prop:excision}, this restriction map is actually an isomorphism, which implies that its domain and target are both zero:
$$SH^*(A\cup B,A)=SH^*(B,A\cap B)=0.$$
By the exact triangles for relative symplectic cohomology of pairs, we deduce that the restriction maps
    $$SH^*(A\cup B)\to SH^*(A),\qquad SH^*(B)\to SH^*(A\cap B)$$
are isomorphisms, as required.
\end{proof}

\section{Proofs of main results}\label{s:proofs_main_results}

In this section, we present proofs of all the results unproved thus far, with the exception of Proposition \ref{prop:excision} (excision) and Theorem \ref{thm: invariance} (invariance), for which we need to go into the construction of the relative symplectic cohomology of pairs. The proofs of these appear in Section \ref{s:excision_invariance}.

\subsection{Proof of Proposition~\ref{prop: weak monotonicity of c_relSH}}\label{ss: proof of monotonicity}

The proof of Proposition~\ref{prop: weak monotonicity of c_relSH} is based on the recent computation of the relative symplectic cohomology of toric balls in $\C P^n$ together with appropriate restriction maps, due to Dickstein--Ganor \cite{DG_2026_relSH_in_CPn}. For $a \in(0,1)$ let $\Delta(a) = \{x\in[0,\infty)^n\,|\,\sum_{j=1}^{n}x_j\leq a\}$ and $B(a)=\Phi^{-1}(\Delta(a))\subset \C P^n$. The following result is a particular case of Theorem 1.1 and Proposition 1.8 from \cite{DG_2026_relSH_in_CPn}:
    \begin{thm}[Theorem 1.1 and Proposition 1.8, \cite{DG_2026_relSH_in_CPn}]\label{thm: SH of ball in CP^n + res maps}
         Normalize the Fubini--Study form on $\C P^n$ so that complex lines have area $1$. Then for $a \in (0,\frac n{n+1})$ we have
            $$SH^*_{\CP^n}(B(a))\cong \bigoplus\limits_{i=0}^\infty \Lambda_{(0,a]}\quad\text{for even }*$$
            and $SH^*_{\C P^n}(B(a))=0$ for odd $*$. Moreover, the restriction map
            $$\res^{\CP^n}_{B(a)} \fc SH_{\CP^n}^0(\CP^n) = \Lambda_{>0} \to SH_{\CP^n}^0(B(a))$$
            is given by
            $$\res^{\CP^n}_{B(a)}(\lambda) = \big((-T^{n(1-a)-a})^{i}\ol\lambda\big)_{i=0}^\infty$$
            for $\lambda \in \Lambda_{>0}$, where $\ol\lambda \in \Lambda_{(0,a]}$ is the image of $\lambda$ by the quotient map $\Lambda_{>0} \to \Lambda_{(0,a]}$. \qed
    \end{thm}
\begin{coroll}\label{cor:capacity_nonheavy_ball_CPn}
  The capacity $c(B(a))$ equals $a$. \qed
\end{coroll}

\begin{proof}[Proof of Proposition~\ref{prop: weak monotonicity of c_relSH}]Normalize the symplectic form on $\C P^n$ so that lines have area $1$. Let $a\in(0,1)$ be such that $B(a)\subset K$. The proof separates into two cases.

\textbf{Case I:} $a\geq\frac{n}{n+1}$. In this case $B(a)$ and $K$ both contain the Clifford torus, and thus are both heavy. From \cite{MSV_2024_heavy_sets_and_SH} we deduce that
$$c(B(a))=c(K)=\infty.$$

\textbf{Case II:} $a<\frac{n}{n+1}$. In this case we have $c(B(a))=a$ by Corollary \ref{cor:capacity_nonheavy_ball_CPn}. Let us show that for every $\lambda\in[0,a)$ we have $T^\lambda SH^*_{\CP^n}(K)\neq0$, which will then imply that $a \leq c(K)$, as required.

Set $x=T^{a-\lambda}\in \Lambda_{>0}= SH_{\CP^n}^0(\CP^n)$, and consider
$$y=\res_K^{\CP^n}(x)\in SH_{\CP^n}^0(K).$$
It suffices to show that $T^\lambda y\neq0$. By functoriality we have
$$\res^{K}_{B(a)}(T^\lambda y)=T^\lambda\res^{K}_{B(a)}(y)=T^\lambda\res^{K}_{B(a)}\big(\res^{\CP^n}_{K}(x)\big)=T^\lambda \res^{\CP^n}_{B(a)}(x),$$
and by Theorem~\ref{thm: SH of ball in CP^n + res maps} we deduce that
$$\res^{K}_{B(a)}(T^\lambda y)=T^\lambda \res^{\CP^n}_{B(a)}(T^{a-\lambda})=\res^{\CP^n}_{B(a)}(T^a)=(\ol{T^a},0,0,\ldots)\neq0\in \bigoplus_{i=0}^\infty \Lambda_{(0,a]}.$$
Therefore $T^\lambda y\neq0$ and this completes the proof.
\end{proof}

\subsection{Proof of Theorem~\ref{thm: general obst for symplectophagy in Cn}}\label{ss:proof_thm_gen_obst_symplectophagy}

For the proof we will need a bit of preparation. The main point here is that commutation of compact subsets is a local property, meaning it is preserved under symplectic embeddings. Also, we will mention a couple of elementary facts regarding displacement energy. Recall Definition \ref{def:commuting_sets} for commutation of closed subsets in a symplectic manifold $(M,\omega)$.
\begin{rem}
  It is an exercise, which we leave for the reader, to show that two closed subsets $A,B\subset M$ commute if and only if there are Poisson-commuting $f,g\in C^\infty(M)$ with $A = f^{-1}(0)$, $B = g^{-1}(0)$.
\end{rem}
\begin{prop}\label{prop: invariance of commuting sets}
    Let $\iota \fc M \hookrightarrow N$ be a symplectic embedding, and let $A,B\subset M$ be compact sets which commute in $M$. Then $\iota(A),\iota(B)$ commute in $N$.
\end{prop}
\begin{proof}
  Without loss of generality we may assume that $M$ is an open symplectic submanifold of $N$ and that $\iota$ is the inclusion map. Moreover, we may assume that $N$ is connected, because otherwise the following argument can be applied to each connected component of $N$ separately. Also, if $M=N$, there is nothing to prove, thus we assume that $M\neq N$, and in particular $A,B\neq N$.  Finally, since the empty set trivially commutes with any other set, we may assume that $A,B$ are nonempty.

  Let $f,g\in C^\infty(M)$ be Poisson-commuting functions such that $A = f^{-1}(0)$, $B=g^{-1}(0)$. Squaring $f,g$, if necessary, we may assume $f,g\geq 0$. Note that squaring preserves Poisson commutation. We will construct Poisson-commuting $F,G \in C^\infty(N)$ with $A = F^{-1}(0)$, $B=G^{-1}(0)$.

  Let $U\subset M$ be an open neighborhood of $A$ with compact closure. Note that $U$, being a proper nonempty open set of the connected topological space $N$, cannot be closed, meaning $\partial U\neq \varnothing$. Since $\partial U = \ol U \setminus U$ is disjoint from $A$, we see that $f|_{\partial U}>0$. Let $a=\inf_{\partial U} f$. Since $\partial U$ is compact, we have $a > 0$. Define $X = \ol U \cap f^{-1}\big([0,\frac a2]\big)$. It clearly is a compact set, and thus closed in $N$. Since $f|_{\partial U} \geq a$, we have $X\cap \partial U=\varnothing$ and thus $X\subset U$. It follows that $N = (N\setminus X)\cup U$ is an open cover of $N$. Let $\chi \fc [0,\infty)\to \R$ be a smooth function such that $\chi(t)=t$ for $t\in[0,\frac a4]$, $\chi|_{[\frac a2,\infty)}\equiv \frac{a}2$, and $\chi(t) > 0$ for $t>0$. Define $F \fc N \to \R$ by $F|_U:=\chi\circ f$, $F|_{N\setminus X}\equiv \frac a2$. Let us check that $F$ is smooth. Since it is given by its restrictions to sets in an open cover, it suffices to check smoothness on each open subset. On $N\setminus X$ it is constant, and thus trivially smooth. On $U$, $F$ is the composition of two smooth functions, and therefore smooth. It remains to check that $F$ is well-defined, meaning that $\chi\circ f|_{(N\setminus X)\cap U}\equiv \frac a2$. Let $x \in (N\setminus X)\cap U = U\setminus X$. Since $X$ is precisely the set of points in $U$ where $f \leq \frac a2$, it follows that $f(x) > \frac a2$, therefore indeed $(\chi\circ f)(x) = \frac a2$, as claimed.

  In a similar fashion we can construct a function $G \in C^\infty(N)$: we pick a precompact neighborhood $V\subset M$ of $B$, let $b=\min_{\partial V}g > 0$, set $Y=\ol V\cap g^{-1}\big([0,\frac b2]\big) = V\cap g^{-1}\big([0,\frac b2]\big)$, thus obtaining an open cover $N = (N\setminus Y)\cup V$. Picking a smooth function $\eta \fc [0,\infty) \to \R$ with $\eta(t)=t$ for $t\in[0,\frac b4]$, $\eta|_{[\frac b2,\infty)}\equiv \frac b2$, and $\eta(t) > 0$ for $t>0$, we can set $G|_V:=\eta\circ g$, $G|_{N\setminus Y} = \frac b2$. Arguing as above, we see that $G$ is smooth.

  Now let us show that $F,G$ Poisson commute. It suffices to do so on $U\cap V$, since elsewhere at least one of $F,G$ is constant. On $U\cap V$ we have $F=\chi\circ f$ and $G=\eta\circ g$, therefore $\{F,G\}=(\chi'\circ f)(\eta'\circ g)\{f,g\} \equiv 0$, since $\{f,g\}\equiv 0$. Finally, it is clear from the construction that $A = F^{-1}(0)$, $B = G^{-1}(0)$.
\end{proof}

\begin{rem}\label{rem:elementary_properties_displ_E} We record here a couple of elementary facts about displacement energy. For a symplectic manifold $(M,\omega)$ we let $e_M$ stand for displacement energy in $M$.
\begin{enumerate}
  \item The displacement energy is monotone nonincreasing under symplectic embeddings: if $\iota \fc M \hookrightarrow N$ is a symplectic embedding, then for compact $A\subset M$ we have $e_M(A) \geq e_N(\iota(A))$. Indeed, if $H$ is a compactly supported Hamiltonian on $M$ displacing $A$, then the extension by zero of $H\circ\iota^{-1}$ is a compactly supported Hamiltonian on $N$ displacing $\iota(A)$.
  \item If $A\subset \C^n$ is compact, then for any $\epsilon > 0$ there is a standard open ball $B\subset \C^n$ such that $e_B(A) \leq e_{\C^n}(A) + \epsilon$. Indeed, let $H$ be a compactly supported Hamiltonian on $\C^n$ with $\|H\|_{\mathrm{Hofer}} \leq e_{\C^n}(A) + \epsilon$. Now take $B$ to be a ball containing both $A$ and the support of $H$, and deduce that $e_B(A) \leq \|H\|_{\mathrm{Hofer}} \leq e_{\C^n}(A) + \epsilon$. In particular if $\C P^n$ is obtained from $\ol B$ by the symplectic cut construction \cite{Lerman_1995_symplectic_cut}, that is we collapse $\partial B$ by the Hopf map onto the corresponding divisor $\C P^{n-1} \subset \C P^n$, then $e_{\C P^n}(A) \leq e_{\C^n}(A)+\epsilon$, where we identify $B=\C P^n\setminus\C P^{n-1}$.
\end{enumerate}
\end{rem}

\begin{prop}\label{prop:capacities_under_symplectophagy}
  Let $X,Y \subset \C^n$ be compact commuting sets, $(\phi_t)_t$ a compactly supported Hamiltonian isotopy, and $a > 0$.
  \begin{enumerate}
    \item If $\ol B{}^{2n}(a)\subset Y$ and $(\phi_t)_t$ is a weak symplectophagy of $Y$ into $X$, then $a \leq e(X\cap Y)$;
    \item If $\ol B{}^{2n}(a)\subset X\cup Y$ and $(\phi_t)_t$ is a genuine symplectophagy of $Y$ into $X$, then $a \leq e(X)$.
  \end{enumerate}
\end{prop}
\begin{proof}
  Let $G$ be a compactly supported Hamiltonian on $\C^n$ generating $(\phi_t)_t$ that is $\phi_t=\phi^t_G$. Fix $\epsilon > 0$ and let $R>\frac{n+1}{n}\,a$ be large enough so that $B^{2n}(R) \supset X\cup Y \cup \supp G$, and so that for $\star \in \{X\cap Y,X\}$ we have $e_{B^{2n}(R)}(\star) \leq e_{\C^n}(\star)+\epsilon$. Let $\C P^n(R)$ be obtained from $\ol B{}^{2n}(R)$ by the symplectic cut construction. We can view $\ol B{}^{2n}(a),X,Y$ as subsets of this $\C P^n(R)$. Since $R>\frac{n+1}{n}\,a$, the ball $\ol B{}^{2n}(R)$ does not contain the Clifford torus, Corollary $\ref{cor:capacity_nonheavy_ball_CPn}$ together with a suitable renormalization of the symplectic form implies that $c(\ol B{}^{2n}(a)) = a$. Note that $X,Y$ still commute in $\C P^n(R)$, thanks to Proposition \ref{prop: invariance of commuting sets}, that $(\phi_t)_t$ is a well-defined Hamiltonian isotopy of $\C P^n(R)$. For the rest of the proof, $c$ is the capacity of compact subsets of $\C P^n(R)$.
  \begin{enumerate}
    \item Note that $(\phi_t)_t$ is still a weak symplectophagy of $Y$ into $X$ in $\C P^n(R)$. We have
    \begin{align*}
      a &= c(\ol B{}^{2n}(a)) && \text{by the above} \\
        &\leq c(Y)            && \text{by Proposition \ref{prop: weak monotonicity of c_relSH}}\\
        &\leq c(X\cap Y)      && \text{by Proposition \ref{prop: obst for weak symplectophagy, V2}}\\
        &\leq e_{\C P^n(R)}(X\cap Y) && \text{by Theorem \ref{thm:energy_capacity_ineq}}\\
        &\leq e_{B^{2n}(R)}(X\cap Y) && \text{by Remark \ref{rem:elementary_properties_displ_E}}\\
        &\leq e_{\C^n}(X\cap Y) + \epsilon && \text{by the choice of }R\,.
    \end{align*}
    Taking $\epsilon\to 0$, we obtain the desired conclusion.
    \item Note that $(\phi_t)_t$ is a genuine symplectophagy of $Y$ into $X$ in $\C P^n(R)$. We have
    \begin{align*}
      a &= c(\ol B{}^{2n}(a))       && \text{by the above} \\
        &\leq c(X\cup Y)            && \text{by Proposition \ref{prop: weak monotonicity of c_relSH}}\\
        &= c(X)                     && \text{by Proposition \ref{prop: SH obstruction for symplectophagy}}\\
        &\leq e_{\C P^n(R)}(X)      && \text{by Theorem \ref{thm:energy_capacity_ineq}}\\
        &\leq e_{B^{2n}(R)}(X)      && \text{by Remark \ref{rem:elementary_properties_displ_E}}\\
        &\leq e_{\C^n}(X) + \epsilon && \text{by the choice of }R\,.
    \end{align*}
    Now take $\epsilon\to 0$.
  \end{enumerate}
\end{proof}

\begin{proof}[Proof of Theorem~\ref{thm: general obst for symplectophagy in Cn}]
Fix $\lambda >0$ and let $\psi \fc B^{2n}(\lambda) \hookrightarrow\C^n$ be a symplectic embedding. Assume $a \in (0,\lambda)$. Then the ``extension after restriction'' principle \cite[p. 7]{Schlenk_book} implies the existence of a Hamiltonian diffeomorphism $\Psi\fc \C^n \to \C^n$ such that $\Psi|_{\ol B{}^{2n}(a)} = \psi|_{\ol B{}^{2n}(a)}$. Assume that $(\phi_t)_t$ is a compactly supported Hamiltonian isotopy. We observe that $\big(\wt\phi_t:=\Psi^{-1}\circ\phi_t\circ\Psi\big)_t$ is again a compactly supported Hamiltonian isotopy, and that $\wt X:=\Psi^{-1}(X),\wt Y:=\Psi^{-1}(Y) \subset \C^n$ are compact and still commute, since commutation is obviously preserved by symplectomorphisms.
\begin{enumerate}
  \item To prove the first assertion of the theorem, assume $\im\psi\subset Y$. The assumption that there is a weak symplectophagy of $Y$ into $X$ implies that there is a compactly supported one, see Remark \ref{rem:symplectophagy_name_cpt_suppt}. Let $(\phi_t)_t$ be such a compactly supported weak symplectophagy of $Y$ into $X$. It follows that $\ol B{}^{2n}(a)\subset \wt Y$ and that $(\wt\phi_t)_t$ is a weak symplectophagy of $\wt Y$ into $\wt X$. Proposition \ref{prop:capacities_under_symplectophagy}, item (i) implies that $a \leq e(\wt X\cap \wt Y)$. Since displacement energy is preserved by symplectomorphisms, it follows that $a \leq e(X\cap Y)$. It is clear that $c_G(Y)$ is the supremum of such $a$ and thus we obtain the desired conclusion.
  \item To prove the second assertion of the theorem, assume $\im\psi\subset X\cup Y$. As in the previous item, we may assume that we are given a compactly supported genuine symplectophagy $(\phi_t)_t$ of $Y$ into $X$. It follows that $\ol B{}^{2n}(a)\subset \wt X \cup \wt Y$ and that $(\wt\phi_t)_t$ is a genuine symplectophagy of $\wt Y$ into $\wt X$. Proposition \ref{prop:capacities_under_symplectophagy}, item (ii) implies that $a \leq e(\wt X) = e(X)$, again by the symplectic invariance of displacement energy. The supremum of such $a$ is exactly $c_G(X\cup Y)$ and the proof is complete.
\end{enumerate}
\end{proof}

\subsection{Proof of Theorem~\ref{thm: main thm2}}\label{ss: Proof of main thm2}
Recall the statement of the theorem: for $\epsilon \in (0,\frac13)$, $\delta \in(0,1-2\epsilon)$, $\eta \in (0,\epsilon)$, and
$$A=\Phi^{-1}\left(\{x\in \im\Phi\,|\,x_1\leq \delta\}\right),\quad B=\Phi^{-1}\left(\{x\in \im\Phi\,|\,x_1\leq 1-\eta-\epsilon,\,x_2=\epsilon\}\right)\,,$$
there is no weak symplectophagy of $B$ into $A$. For the proof we will use the capacity of the boundary of a large ball in $\C P^n$, which can be immediately deduced from the following
\begin{thm}[{\cite[Theorem 1.15]{DG_2026_relSH_in_CPn}}]
  Let $n$ be a positive integer, and normalize the symplectic form on $\C P^n$ so that lines have area $1$. Then for $a \in (\frac n{n+1},1)$ we have
  $$SH^*_{\C P^n}(\partial B(a)) = \bigoplus_{k=0}^\infty\Lambda_{(0,n(1-a)]}\quad \text{if }*\text{ is even}\,,$$
  and $SH^*_{\C P^n}(\partial B(a))=0$ if $*$ is odd. \qed
\end{thm}
\begin{coroll}\label{cor:capacity_bdry_large_ball_CPn}
  For $a\in(\frac{n}{n+1},1)$, we have $c(\partial B(a))=n(1-a)$.\qed
\end{coroll}
\begin{proof}[Proof of Theorem~\ref{thm: main thm2}]
We will show that
$$c(A\cap B) \leq \epsilon < c(B)\,,$$
which by Proposition \ref{prop: obst for weak symplectophagy, V2} implies that there exists no weak symplectophagy of $B$ into $A$. For the first inequality, note that for any $\gamma \in \big(0,\frac12(1-\delta-2\epsilon)\big)$, the set
$$P = \Phi^{-1}([0,\delta+\gamma)\times[0,2\epsilon+\gamma))$$
is an embedded open toric polydisk, inside which the set $A\cap B$ has displacement energy $\leq \epsilon$. Thus
$$c(A\cap B) \leq e_{\C P^2}(A\cap B) \leq e_P(A\cap B)\leq \epsilon\,.$$

Now let us show that $c(B) > \epsilon$. Define
$$S = \Phi^{-1}([0,1-\epsilon]\times\{\epsilon\})\,,\quad K = \Phi^{-1}([1-\epsilon-\eta,1-\epsilon]\times\{\epsilon\})\,.$$
Note that $S$ is the boundary of the ball $\Phi^{-1}(\R\times[\epsilon,\infty))$, which is symplectomorphic to the toric ball $\Phi^{-1}(\Delta(1-\epsilon))$. Since its interior contains the Clifford torus, Corollary \ref{cor:capacity_bdry_large_ball_CPn} implies that $c(S) = 2\epsilon$. We have $S = B\cup K$, and since $B,K$ commute, Proposition \ref{prop:capacity_estimate_commuting_sets} implies that
$$2\epsilon = c(S) \leq \max(c(B),c(K)) + c(B\cap K)\,.$$
By an argument similar to the one appearing at the beginning of the proof, $K$ has displacement energy $e(K) \leq \eta<\epsilon$ and so $e(B\cap K) \leq e(K) < \epsilon$. We claim that
$$\max(c(B),c(K)) = c(B)\,.$$
Indeed, otherwise we would have $c(B) < c(K)$, implying that
$$2\epsilon \leq c(K) + c(B\cap K) \leq e(K) + e(B\cap K) <\epsilon+\epsilon = 2\epsilon\,,$$
which is a contradiction. We deduce that
$$2\epsilon \leq \max(c(B),c(K)) + c(B\cap K) = c(B) + c(B\cap K) \leq c(B) + e(B\cap K) < c(B) + \epsilon\,,$$
whence
$$c(B) > 2\epsilon - \epsilon = \epsilon\,,$$
as claimed.
\end{proof}

\subsection{Proof of Theorem \ref{thm: main thm3}}\label{ss:proof_main_thm_3}

Let us recall the statement of the theorem: if $k,n$ are positive integers, $0<a_0<b_0$ are real numbers, and
$$D = \ol B{}^{2k}(a_0) = \{z_0\in \C^k\,|\,\pi\|z_0\|^2\leq a_0\}\,,\quad A = \{z_0\in \C^k\,|\,a_0\leq \pi\|z_0\|^2\leq b_0\}\,,$$
then for any two closed nondegenerate toric polydisks $X,Y \subset \C^n$, there is no genuine symplectophagy of $A\times Y$ into $D\times X$ in $\C^{n+k} = \C^k \times \C^n$.

For $a_1,\dots,a_n>0$ let
$$\ol P(a_1,\dots,a_n) = \prod_{j=1}^{n}\ol B{}^2(a_j)\subset \C^n$$
denote the corresponding nondegenerate closed toric polydisk. We have
$$c_G\big(\ol P(a_1,\dots,a_n)\big)= e_{\C^n}\big(\ol P(a_1,\dots,a_n)\big) = \min_{1\leq j\leq n}a_j$$
due to Gromov's nonsqueezing \cite{Gromov_1985_J-hol_curves}.

The proof of the theorem uses the following elementary observation, which immediately follows from the definitions:
\begin{claim}
  Let $(M,\omega)$ be a symplectic manifold, $A,B\subset M$, and $(\phi_t)_t$ be a genuine symplectophagy of $B$ into $A$. Let $A'\subset M$ be such that $A\subset A'$ and $A\cap B = A'\cap B$. Then $(\phi_t)_t$ is also a genuine symplectophagy of $B$ into $A'$. \qed
\end{claim}

\begin{figure}
	\centering
	\includegraphics[width=0.25\linewidth]{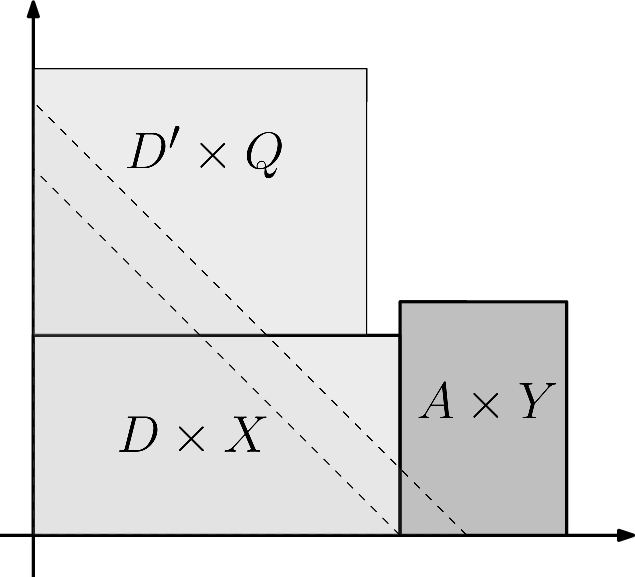}
	\caption{Proof of Theorem \ref{thm: main thm3}, the case $k=n=1$. The smaller dashed triangle symbolizes a toric ball with the same displacement energy as $D\times(X\cup Q)$. The larger triangle refers to a strictly larger toric ball embeddable in $D\times X \cup D'\times Q \cup A \times Y$}
	\label{fig:pfThm3}
\end{figure}

We argue by contradiction and assume that $(\phi_t)_t$ is a genuine symplectophagy of $A\times Y$ into $D\times X$. If $a_0' \in (0,a_0)$, then, thanks to the claim, for any compact $Q\subset \C^n$, $(\phi_t)_t$ is also a genuine symplectophagy of $A\times Y$ into the larger set $D\times X \cup D'\times Q$, where $D'=\ol B{}^{2k}(a_0')$, since
$$(D\times X \cup D'\times Q)\cap(A\times Y) = D\times X \cap A \times Y\,,\quad\text{as an immediate check shows.}$$
In case $Q$ is also toric, it then follows from Theorem \ref{thm: general obst for symplectophagy in Cn}, item (ii) that
$$c_G(D\times X \cup D'\times Q \cup A \times Y) \leq e(D\times X \cup D'\times Q)\,.$$
Since $D\times X \cup D'\times Q \subset D\times(X\cup Q)$, it follows that
$$e(D\times X \cup D'\times Q) \leq e\big(D\times(X\cup Q)\big) \leq e(D)=a_0\,.$$
Thus, if we can find $a_0'\in(0,a_0)$ and a compact toric $Q\subset \C^n$ such that
$$c_G(D\times X \cup D'\times Q \cup A \times Y) > a_0\,,$$
then we will arrive at a contradiction, thereby proving the theorem. We refer the reader to Figure \ref{fig:pfThm3} for an illustration of the geometry.

Let the numbers $a_1,\dots,a_n,b_1,\dots,b_n>0$ be given by
$$X = \ol P(a_1,\dots,a_n)\,,\quad Y = \ol P(b_1,\dots,b_n)\,.$$
Define
$$a=\min_{1\leq j\leq n}a_j\quad\text{and} \quad b = \min_{1\leq j\leq n}b_j\,,\quad\text{and note that }a,b>0\,.$$
Pick any $a_0' \in \big(\max\{0,a_0-a\},a_0\big)$ and set
$$r=\min\{b_0,a_0+b,a_0'+a\}\quad \text{and} \quad Q=\ol P(r,\dots,r) \subset \C^n\,.$$
We claim:
\begin{enumerate}
  \item $r>a_0$;
  \item for $D'=\ol B{}^{2k}(a_0')$ we have $\ol B{}^{2(n+k)}(r)\subset D\times X \cup D'\times Q\cup A \times Y$.
\end{enumerate}
First, note that these imply that
$$c_G(D\times X \cup D'\times Q\cup A \times Y) \geq r > a_0\,,$$
leading to the desired contradiction. It therefore remains to prove the assertions. For (i) note that $b_0>a_0$, $a_0+b>a_0$, since $b>0$, and that by our choice
$$a_0'+a > \max\{0,a_0-a\}+a \geq a_0-a+a=a_0\,,$$
and so these inequalities imply that $r>a_0$.

Let us verify assertion (ii). Pick $z=(z_0,z_1,\dots,z_n) \in \ol B{}^{2(k+n)}(r)$, where $z_0\in \C^k$ and $z_1,\dots,z_n\in \C$. Note that $\pi\|z_0\|^2\leq r$. The verification is split into three cases:
\begin{itemize}
  \item $\pi\|z_0\|^2 \leq a_0'$: in this case $z_0 \in D'$; also for each $j=1,\dots,n$, $\pi|z_j|^2 \leq r$, meaning $(z_1,\dots,z_n) \in Q$, that is $z \in D'\times Q$.
  \item $\pi\|z_0\|^2 \in [a_0',a_0]$: here we have $z_0 \in D$ and
  $$\pi\big(|z_1|^2+\dots+|z_n|^2\big) = \pi\|z\|^2-\pi\|z_0\|^2 \leq r-a_0'\leq a_0'+a-a_0'=a\,,$$
  which implies that for each $j=1,\dots,n$ we have $\pi|z_j|^2\leq a \leq a_j$, that is $(z_1,\dots,z_n) \in X$, and ultimately $z\in D\times X$.
  \item $\pi\|z_0\|^2 \in [a_0,r]$: here $z_0 \in A$, but also we have
  $$\pi\big(|z_1|^2+\dots+|z_n|^2\big) = \pi\|z\|^2-\pi\|z_0\|^2 \leq r-a_0\leq a_0+b-a_0=b\,,$$
  meaning that for each $j=1,\dots,n$ we have $\pi|z_j|^2 \leq b \leq b_j$, therefore $(z_1,\dots,z_n) \in Y$, that is $z \in A\times Y$.
\end{itemize}
We see that in all cases $z \in D\times X\cup D'\times Q\cup A\times Y$, as claimed. The proof is complete.

\begin{rem}\label{rem:toric_ex_embeddability_criterion}
As we mentioned in Remark \ref{rem:rk_main_thm_3}, in case $k\leq n$, for $Y$ small enough and $X$ large enough, there is a Hamiltonian isotopy of $A\times Y$ into $D\times X$, which means that there is no symplectic embeddability obstruction to the existence of symplectophagy. Precisely, if $b_1+\dots+b_k \leq a_0$, $b_{k+1}\leq a_{k+1},\dots,b_n\leq a_n$, and $b_0 \leq \min(a_1,\dots,a_k)$, then the unitary transformation $U\in U(k+n)$, defined by $U(e_j)=e_{j+k}$ and $U(e_{j+k})=e_j$ for $j=1,\dots,k$, and $U(e_j)=e_j$ for $j>2k$, maps $A\times Y$ into $D\times X$, as an elementary calculation shows. Since unitary transformations are Hamiltonian, we have thus obtained a Hamiltonian isotopy mapping $A\times Y$ into $D\times X$, and by cutting off the Hamiltonian, its support can be made compact.
\end{rem}

\section{Excision and invariance}\label{s:excision_invariance}

This section is dedicated to the proof of Proposition \ref{prop:excision} (excision) and Theorem \ref{thm: invariance} (invariance). For this we will briefly recall the definitions of the relative symplectic cohomology of pairs and the associated notions. Throughout, $(M,\omega)$ is a fixed closed symplectic manifold.

\subsection{Relative symplectic cohomology of pairs}\label{ss:basic_notions_defs}

In this section we will recall the relevant notions and constructions used in the definition of the relative symplectic cohomology of compact sets and pairs thereof as in \cite{Varolgunes_2021_MV_and_relSH,DGPZ_2024_Symp_top_and_IVMs}.

\subsubsection{Cubes}\label{sss:cubes}

We briefly review the necessary notions following Varolg\"une\c s \cite{Varolgunes_2021_MV_and_relSH} and \cite{DGPZ_2024_Symp_top_and_IVMs}. Throughout we use the standard Euclidean $n$-cube $[0,1]^n$, where $n$ is a nonnegative integer. A \tb{face} of the $n$-cube is any subset $F\subset[0,1]^n$ given by setting some of the coordinates to $0$ or to $1$. The rest of the coordinates are \tb{free}, and the number of the free coordinates is the dimension $\dim F$ of $F$. Each face $F$ determines and is determined by its \tb{initial} and \tb{terminal} vertices $\ini F,\ter F\in F$, which are the vertices closest to, respectively farthest from the origin $0\in[0,1]^n$, with respect to the usual Euclidean distance.

We will use throughout the Morse function $\rho_1\fc[0,1] \to \R$ given by $\rho_1(x)=\cos\pi x$ and the induced Morse function $\rho_n$ on $[0,1]^n$ given by $\rho_n(x)=\sum_{j=1}^{n}\rho_1(x_j)$. We will use the gradient vector field $\nabla \rho_n$ of $\rho_n$ with respect to the Euclidean metric. Note that it is tangent to every face of the cube. In particular, the critical points of $\rho_n$ are precisely the vertices of $[0,1]^n$. We will also need to use slits. The \tb{$2$-slit} is the subset $\Slit^2\subset[0,1]^2$ which is the union of the closures of the negative gradient trajectories of $\rho_2$ passing through the points $(t,1-t)$ for $t\in[\frac13,\frac23]$. For $n > 2$, the \tb{$n$-slit} is $\Slit^n=[0,1]^{n-2}\times \Slit^2$.

For a nonnegative integer $n$, a \tb{Hamiltonian $n$-cube} is a map $\cG \fc [0,1]^n \to C^\infty(S^1\times M)$ which is constant on a neighborhood of each vertex, and which is smooth in the sense that the resulting map $[0,1]^n\times S^1\times M \to \R$ is smooth. A Hamiltonian $n$-cube $\cG$ is \tb{monotone} if for each negative gradient trajectory $\gamma \fc \R \to [0,1]^n$ of $\rho_n$, $\cG\circ\gamma$ is a nondecreasing function relative to the usual order on $C^\infty(S^1\times M)$. Monotone Hamiltonian cubes are of particular importance in the theory. A \tb{Hamiltonian $n$-slit} is a smooth map $\cG \fc \Slit^n \to C^\infty(S^1\times M)$, which is constant in a neighborhood of each vertex. We say that $\cG$ is monotone if it is nondecreasing along negative gradient trajectories of $\rho_n$.

A \tb{Hamiltonian $n$-ray} is a collection $\cG_\bullet=(\cG_i)_{i\in \N}$ of Hamiltonian $n$-cubes such that for each $i\in \N$, $\cG_i(1,x_2,\dots,x_n)=\cG_{i+1}(0,x_2,\dots,x_n)$ for all $(x_2,\dots,x_n)\in[0,1]^{n-1}$, or in words, the $x_1=1$ wall of $\cG_i$ coincides with the $x_1=0$ wall of $\cG_{i+1}$. A ray is monotone if its constituent cubes are. For $n\geq 3$, a \tb{Hamiltonian slit-like $n$-ray} is a collection $\cG_\bullet=(\cG_i)_{i\in \N}$ of Hamiltonian $n$-slits with the same compatibility condition at the $x_1$ walls, but now for $(x_2,\dots,x_n)\in\Slit^{n-1}$. A Hamiltonian slit-like ray is monotone if its constituent Hamiltonian slits are.

Given $A\subset M$, we denote $C^\infty_{A\subset M} = \{H\in C^\infty(S^1\times M)\,|\;H|_{S^1\times A}<0\}$. If $K\subset M$ is compact, an \tb{acceleration datum\footnote{There are two distinct meanings of this term in the literature, one as defined here, and the other meaning what we here call an acceleration $1$-ray. We decided to separate the two notions by giving them distinct names.} for $K$} is any pointwise nondecreasing sequence of nondegenerate Hamiltonians $(H_i)_{i\in \N}$ such that $H_i \in C^\infty_{K\subset M}$ for all $i$ and such that we have the pointwise limits $H_i(t,x) \to 0$ for $x\in K$ and $H_i(t,x) \to \infty$ for $x\notin K$, as $i\to \infty$.

We will also find use for the notion of acceleration rays. For $n\geq 0$, a \tb{system of $2^n$ sets}, or a system of sets for short, is a map $\{0,1\}^n\to\cC(M)$, $\star\mapsto K^\star$ satisfying $\star\leq\star'\Rightarrow K^\star\supset K^{\star'}$, where the order on $\{0,1\}^n$ is the product order induced from the standard one on $\{0,1\}$. Given such a system of sets $(K^\star)_\star$, a Hamiltonian $(n+1)$-ray $\cG_\bullet=(\cG_i)_{i\in\N}$ is an \tb{acceleration ray for $(K^\star)_\star$} if it is monotone and such that for each $\star\in\{0,1\}^n$, the sequence $(\cG_i(0,\star))_{i\in\N}$ is an acceleration datum for $K^\star$ (note that $(0,\star)\in\{0\}\times\{0,1\}^n\subset[0,1]^{n+1}$). In the special case $n=0$ we have $\{0,1\}^0=$ a singleton, so that $(K^\star)_\star$ in this case is just given by one set $K$. Then a $1$-ray $\cG_\bullet=(\cG_i)_{i\in\N}$ is an acceleration $1$-ray for $K$ if and only if it is monotone and $(\cG_i(0))_{i\in\N}$ is an acceleration datum for $K$. Other special cases of particular interest are an \tb{acceleration $2$-ray for a compact pair $(K,K')$}, which is an acceleration $2$-ray for the system $\{0,1\}\to\cC(M)$, $0\mapsto K$, $1\mapsto K'$, and an \tb{acceleration $3$-ray for an inclusion of pairs $(K,K')\subset (L,L')$}, by which we mean an acceleration $3$-ray for the system $\{0,1\}^2 \to \cC(M)$, $(0,0)\mapsto L$, $(0,1)\mapsto L'$, $(1,0)\mapsto K$, $(1,1)\mapsto K'$.

Acceleration rays always exist:
\begin{prop}[{Varolg\"une\c s \cite[Section 3.4 and Proposition 3.2.13]{Varolgunes_2021_MV_and_relSH}}]
  For any system of $2^n$ sets as above there exists an acceleration $(n+1)$-ray. \qed
\end{prop}

We have similar notions---cubes, slits, rays, and slit-like rays---for families of almost complex structures. Let $\cJ(M,\omega)$ be the space of $\omega$-compatible almost complex structures and $C^\infty(S^1,\cJ(M,\omega))$ the space of smooth $S^1$-families of such structures. An $n$-cube of almost complex structures is a smooth map $\cI\fc [0,1]^n \to C^\infty(S^1,\cJ(M,\omega))$, which is assumed to be constant on a neighborhood of each vertex. The rest of the definitions are similarly formulated.

Next we recall the relevant algebraic notions. Fix a base ring, which is assumed unital and commutative. An \tb{algebraic $n$-cube} is a pair $\cC=\big((\cC^v)_{v\in\{0,1\}^n},(f^\cC_F)_{F\subset[0,1]^n}\big)$, where the $\cC^v$ are graded\footnote{The grading may be over $\Z/(2k)$. In any case it will be clear from the context.} modules over our base ring, called the \tb{vertices of $\cC$}, while for each face $F\subset[0,1]^n$, $f^\cC_F \fc \cC^{\ini F} \to \cC^{\ter F}$ is the corresponding \tb{face map}, which is a graded module map of degree $1-\dim F$. These maps are subject to the condition
$$\text{for each face }F\subset[0,1]^n\qquad \sum_{F=F'\cdot F''}(-1)^{\dim F'}\sgn(F',F'')f^\cC_{F''}\circ f^\cC_{F'} = 0\,.$$
Here $F=F'\cdot F''$ denotes the situation in which $F',F''\subset [0,1]^n$ are faces satisfying
$$\ini F=\ini F'\,,\quad \ter F'=\ini F''\,,\quad \ter F'' = \ter F\,.$$
The sign $\sgn(F',F'')$ is defined as follows: equip the faces $F',F''$ with orientations coming from the natural ordering of coordinates; then $\sgn(F',F'')$ is the oriented intersection number of $F',F''$ in $F$. See \cite[Section 2.1.1, Appendix A]{Varolgunes_2021_MV_and_relSH} for another interpretation.
\begin{rem}
  The face map identities mean in particular that for each vertex $v$, $f^\cC_{\{v\}}$ is a differential on the module $\cC^v$, making it into a chain complex. For a one-dimensional face $F$, $f^\cC_F$ is a chain map between the chain complexes at the ends of $F$, while the face map corresponding to a two-dimensional face is a homotopy between the two compositions of the chain maps associated to its edges, and so on.
\end{rem}

Given algebraic $n$-cubes $\cC = \big((\cC^v)_v,(f^\cC_F)_F\big)$, $\cC' = \big((\cC'^v)_v,(f^{\cC'}_F)_F\big)$, an \tb{isomorphism} $\cC \to \cC'$ is a collection $\Phi = (\Phi^v)_v$, where for each $v$, $\Phi^v \fc \cC^v \to\cC'^v$ is a degree zero graded module isomorphism, such that $\Phi^{\ter F}\circ f^\cC_F = f^{\cC'}_F\circ\Phi^{\ini F}$ for each face $F\subset [0,1]^n$. The $\Phi^v$ are the \tb{vertex maps} of the isomorphism. In the language of \cite{DGPZ_2024_Symp_top_and_IVMs}, an isomorphism of cubes is a straight cube whose structure maps are all isomorphisms.

There is a special kind of cubes we'll need which are called slits in \cite[p. 559]{Varolgunes_2021_MV_and_relSH}. Here we will instead use the term \tb{slit-like cube}. This is a cube of the form
$$\xymatrix{\cA \ar[r]^{\id} \ar[d] \ar[rd] & \cA \ar[d] \\ \cB \ar[r]^{\id} & \cB}$$
where $\cA,\cB$ are algebraic $(n-2)$-cubes, whose coordinates are the first $n-2$ ones of the slit-like cube, while the two displayed coordinates are the last two. The top and bottom arrows denote $(n-1)$-cubes with $\cA$, respectively $\cB$ as the corresponding $(n-2)$-dimensional faces, and where $\id$ stands for face maps where the only nonzero maps are identities connecting identical vertices of $\cA,\cB$.
\begin{rem}\label{rem:slit_like_2_cube}
  A slit-like $2$-cube is equivalently given by two chain complexes, two chain maps between them, and a homotopy between those chain maps.
\end{rem}

We will also use an operation on cubes called the cone. It transforms an $n$-cube into an $(n-1)$-cube, see \cite[Section 2.1.4]{Varolgunes_2021_MV_and_relSH}, \cite[Section 4.4]{DGPZ_2024_Symp_top_and_IVMs}, and it can be performed in any one of the $n$ coordinate directions. The cocone, which we denote by $\co$, is defined as the cone but shifted by $-1$ in degree, see \cite[Section 4.4]{DGPZ_2024_Symp_top_and_IVMs}.

An \tb{algebraic $n$-ray} is a collection of algebraic $n$-cubes $\cC_\bullet=(\cC_i)_{i\in \N}$ such that for each $i\in\N$, the face of $\cC_i$ corresponding to $x_1=1$ coincides with the face of $\cC_{i+1}$ corresponding to $x_1=0$. Likewise we have \tb{slit-like algebraic $n$-rays}, which are just algebraic $n$-rays composed of slit-like $n$-cubes. In \cite{Varolgunes_2021_MV_and_relSH,DGPZ_2024_Symp_top_and_IVMs} a procedure is described, called the \tb{telescope}, which produces an algebraic $n$-cube from an algebraic $(n+1)$-ray. We will not repeat the definition here. The only thing to note is that the telescope of a slit-like ray is a slit-like cube.

\subsubsection{The Floer functor}\label{sss:Floer_fctr}

The most important point of this section is the existence of an operation, which we will here refer to as the \tb{Floer functor},\footnote{It is not an actual functor, strictly speaking, but can be interpreted as an $\infty$-functor between suitable $\infty$-categories.} and which is defined in \cite{Varolgunes_2021_MV_and_relSH}. This operation transforms pairs $(\cG,\cI)$, where $\cG$ is a monotone Hamiltonian $n$-cube such that all the Hamiltonians at the vertices are nondegenerate, and where $\cI$ is an $n$-cube of almost complex structures, into an algebraic $n$-cube over the Novikov ring $\Lambda_{\geq 0}$ denoted by $\mathfrak{F}(\cG,\cI)$. Since later we will need a slight modification of this construction in case $\cG$ is replaced with a cube of Hamiltonian $1$-forms, we will briefly review it.

\begin{rem}\label{rem:Floer_thry_conventions}\phantom{M}
  \begin{itemize}
    \item We use the sign conventions of \cite{DGPZ_2024_Symp_top_and_IVMs} for Hamiltonian vector fields, action functionals, almost complex structures, and so on. As a result, our Floer equation takes the form
        $$\partial_su - J(\partial_tu-X_H) = 0$$
        for the positive gradient flow of the action functional of $H$. It also follows that when bubbling occurs, we obtain $J$-\emph{antiholomorphic} spheres rather than $J$-holomorphic ones. Nevertheless, in order to make the exposition consistent with the accepted terminology, and because we only use antiholomorphic spheres, we will still refer to them as \emph{holomorphic} spheres. Finally, when discussing energy contributions, we will always use the positive quantity $|\int_{S^2}u^*\omega|$ for such an antiholomorphic sphere $u$. For the reader who is uncomfortable with these conventions, we suggest using the conventions of \cite{Varolgunes_2021_MV_and_relSH}, which are related to ours by time-reversal, see \cite{DGPZ_2024_Symp_top_and_IVMs} for details.
    \item We use the virtual fundamental class approach of Pardon \cite{Pardon_2016_alg_vir_counting,Pardon_2019_contact} in the definition of the Floer functor and related structures. In what follows we shall state several claims regarding certain virtual counts arising from this approach. A rigorous verification of those claims would entail going into the details of Pardon's construction at a considerable depth, and even stating the requisite definitions would likely substantially increase the length of the paper. On the other hand, the claims we make are more or less trivial, in the sense that they are all of the form ``certain symmetries induce natural identifications between certain objects.'' We will indicate how this happens in each specific instance, in hopes that the interested reader will be able to fill in the details.
  \end{itemize}
\end{rem}

For a nondegenerate Hamiltonian $H\in C^\infty(S^1\times M)$ let $\cP^\circ(H)$ denote the set of its contractible $1$-periodic orbits. It is graded modulo $2N_M$ by the Conley--Zehnder index. The module underlying the Floer complex of $H$ is the $\Z/(2N_M)$-graded Novikov module
$$CF^*(H) = \bigoplus_{x\in\cP^\circ(H)}\Lambda_{\geq0}\cdot x\,.$$
Fix a $1$-periodic time-dependent $\omega$-compatible almost complex structure $J$. There exist additional choices (see \cite{Varolgunes_2021_MV_and_relSH}, \cite[Definition 7.5.3]{Pardon_2016_alg_vir_counting}) such that for any $x_\pm\in\cP^\circ(H)$, a number $n(x_-,x_+,H,J,A) \in \Q$ is defined, where $A$ is a homotopy class of maps $[-1,1]\times S^1\to M$ whose restrictions to $S^1\cong\{\pm1\}\times S^1$ coincide with $x_\pm$.

The number $n(x_-,x_+,H,J,A)$ is the signed virtual count associated to the compactified moduli space of stable Floer trajectories from $x_-$ to $x_+$ having total homotopy class $A$, using fixed coherent orientations and auxiliary virtual data. The smooth unbroken trajectories are maps $u \fc \R \times S^1\to M$ satisfying
$$\partial_su - J_t(u)\big(\partial_tu-X_H(u)\big) = 0\,,\qquad \lim_{s\to\pm\infty} u(s,\cdot)=x_\pm\,,$$
considered modulo translation in the $s$-coordinate. The compactification also includes configurations with trees of sphere bubbles attached to the cylinders, as well as broken trajectories consisting of successive cylinders with matching asymptotic orbits. Each cylinder in a broken trajectory may carry sphere bubbles. Every sphere in a tree attached at $(s,t)$ is $J_t$-holomorphic\footnote{That is, antiholomorphic, according to our convention in Remark \ref{rem:Floer_thry_conventions}.}. These configurations are considered up to translation of the cylindrical components and reparametrization of the sphere bubbles, subject to stability. The total class of such a configuration is required to equal $A$ (see \cite[Section 3]{Varolgunes_2021_MV_and_relSH} regarding this). See \cite[Definitions 10.2.2, 10.2.3 and Sections 10.3, 10.4]{Pardon_2016_alg_vir_counting} for the precise definitions and construction.

The Floer differential is then defined by the formula
$$dx_- = \sum_{x_+\in\cP^\circ(H)}\sum_A n(x_-,x_+,H,J,A)T^{E_{\mathrm{top}}(x_\pm,H,A)}x_+\,,$$
where $A$ runs over the homotopy classes associated to $x_\pm$. Here
$$E_{\mathrm{top}}(x_\pm,H,A) = \int_{S^1}\big(H_t(x_+(t))-H_t(x_-(t))\big)\,dt - \int_A\omega$$
is the so-called topological energy. It is the ``increase in action'' of $H$ along the class $A$. In terms of the Floer cylinder $u$ and the attached holomorphic spheres, it is the sum of the areas of the spheres plus the actual increase in action along $u$.\footnote{There is a canonically defined $1$-form $\mathsf A_H(x)(\eta)=\int_{S^1}\omega\big(\dot x(t) - X^t_H(x(t)),\eta(t)\big)\,dt$ on the loop space of $M$, and, viewing $u$ as a map from $\R$ to the loop space, this increase is simply the integral $\int_\R u^*\mathsf A_H$.} The main point here is that the topological energy is always nonnegative, since the holomorphic spheres contribute a positive number, while the cylindrical component satisfies the Floer equation, which is the positive gradient flow of the action. It follows that the above differential is well-defined over the Novikov \emph{ring}.

Let us now describe what happens in the case of cubes. Let us fix a monotone Hamiltonian $n$-cube $\cG$ with nondegenerate Hamiltonians at the vertices, and a cube $\cI$ of compatible $S^1$-dependent almost complex structures. Again, some extra choices having been fixed, to these data there is associated an algebraic cube $\mathfrak F(\cG,\cI)$, whose vertices are the underlying graded Novikov modules of the Hamiltonians $\cG(v)$:
$$\mathfrak F(\cG,\cI)^v = CF^*(\cG(v))\,.$$
Let us describe the corresponding face maps. For a face $F\subset [0,1]^n$ with $v=\ini F$, $v'=\ter F$, any $x\in\cP^\circ(\cG(v))$, $x'\in\cP^\circ(\cG(v'))$, there are numbers $n_F(x,x',\cG,\cI,A) \in \Q$, which provide virtual counts of objects of a similar nature to what is described above. Namely, these objects consist of a pair $(\gamma,u)$, where $\gamma \fc \R \to F$ is a negative gradient line of $\rho_n$ asymptotic to $v,v'$, and $u$ is a smooth cylinder together with holomorphic spheres attached, where $u$ satisfies the equation
$$\partial_su - \cI(\gamma(s))(u)(\partial_tu-X_{\cG(\gamma(s))}(u))=0\,.$$
Again, a sphere attached to $u$ at $u(s,t)$ is $\cI(\gamma(s))_t$-holomorphic. Pairs $(\gamma,u)$ obtained from one another by simultaneous shift in the $s$-coordinate in $\gamma$ and in $u$ are considered the same for the purposes of the count. One then defines the corresponding face map
$$f_F^{\mathfrak F(\cG,\cI)}(x) = \sum_{x'\in\cP^\circ(\cG(v'))}\sum_An_F(x,x',\cG,\cI,A)T^{E_{\mathrm{top}}(x,x',\cG,F,A)}x'\,.$$
Here again we have the topological energy
\begin{equation}\label{eq:top_energy_cube}
  E_{\mathrm{top}}(x,x',\cG,F,A) = \int_{S^1}\big(\cG(v')_t(x'(t)) - \cG(v)_t(x(t))\big)\,dt - \int_A\omega\,,
\end{equation}
which is always nonnegative, thanks to the monotonicity of $\cG$.

If we now specialize to the case where $F = \{v\}$ is a single vertex, then there is a unique negative gradient trajectory of $\rho_n$ connecting $v$ to itself, namely the constant trajectory at $v$. It follows that the above equation reduces to the Floer equation for the pair $(\cG(v),\cI(v))$, that is the corresponding face map is the Floer differential. Varolg\"une\c s \cite{Varolgunes_2021_MV_and_relSH} proves that $\big((\mathfrak{F}(\cG,\cI)^v)_v,(f_F^{\mathfrak{F}(\cG,\cI)})_F\big)$ is indeed an algebraic $n$-cube.

If we have a pair $(\cG_\bullet,\cI_\bullet)$ consisting of a ray of Hamiltonians and a ray of almost complex structures, then applying the Floer functor cube-wise, we arrive at an algebraic ray $\mathfrak{F}(\cG_\bullet,\cI_\bullet) = (\mathfrak{F}(\cG_i,\cI_i))_{i\in\N}$. The Floer functor can also be defined for pairs $(\cG,\cI)$ where $\cG$ is a monotone Hamiltonian $n$-slit and $\cI$ is an $n$-slit of almost complex structures. It yields a slit-like algebraic cube $\mathfrak{F}(\cG,\cI)$ (\cite{Varolgunes_2021_MV_and_relSH}). Similarly, if we have a pair consisting of a slit-like ray of Hamiltonians and a slit-like ray of almost complex structures, the Floer functor will transform it into a slit-like algebraic ray.

\begin{notation}\label{not:Floer_1_ray}
  If $\cG_\bullet=(\cG_i)_{i\in\N}$ is a monotone $1$-ray with nondegenerate vertices and $\cI_\bullet=(\cI_i)_{i\in\N}$ is a $1$-ray of almost complex structures, we will use an alternative notation for the algebraic $1$-ray $\mathfrak{F}(\cG_\bullet,\cI_\bullet)$, namely $CF^*(\cG)$. Although slightly imprecise, it is suggestive, and serves to remind the reader that the vertices of the ray $\mathfrak{F}(\cG_\bullet,\cI_\bullet)$ are the Floer complexes of the Hamiltonians $\cG_i(0)$.
\end{notation}

\subsubsection{Relative symplectic cohomology}\label{sss:relSH_pairs}

Here we will briefly recall the definition of the relative symplectic cohomology of a compact set and of a pair of compact sets in $M$, following \cite{Varolgunes_2021_MV_and_relSH,DGPZ_2024_Symp_top_and_IVMs}, as well as the relevant restriction maps. First, if $A$ is a Novikov module, its \tb{completion} is the Novikov module
$$\wh A = \varprojlim_{r\to\infty} A/T^rA\,.$$
We refer the reader to \cite[Section 2.3]{Varolgunes_2021_MV_and_relSH}, \cite[Section 4.1]{DGPZ_2024_Symp_top_and_IVMs} for details, examples, and properties. Here we will briefly describe the features we need. Completion is a functor, meaning a module map $A\to B$ naturally induces a module map $\wh A \to \wh B$. For graded Novikov modules completion is defined degreewise. Since it is a functor, the completion of a chain complex (graded or not) is again a chain complex. More generally, we can complete an algebraic cube by completing at each vertex and plugging in the naturally induced face maps. Thanks to the functoriality of completion, the result is again a cube. Finally, completion commutes with (co)cones.

If $K\subset M$ is compact, fix an acceleration $1$-ray $\cG_\bullet$ for it and a $1$-ray $\cJ_\bullet$ of almost complex structures. Applying the Floer functor to $(\cG_\bullet,\cJ_\bullet)$ we obtain an algebraic $1$-ray $\mathfrak{F}(\cG_\bullet,\cJ_\bullet)$. Its telescope is a $0$-cube over the Novikov ring, or in other words, a chain complex. Completing it, we obtain a chain complex over $\Lambda_{\geq0}$ whose cohomology is $SH^*(K)$. It is shown in \cite{Varolgunes_2021_MV_and_relSH} that the cohomology is independent of the choices. The chain complex we have just described does depend on the various choices. However, we will just use the simplified notation $SC^*(K)$ for it, remembering that it in fact depends on choices other than $K$.

Given a compact pair $(K,K')$ in $M$, fix an acceleration $2$-ray $\cG_\bullet$ for $(K,K')$. Pick a $2$-ray of almost complex structures $\cJ_\bullet$ and apply the Floer functor to the pair $(\cG_\bullet,\cJ_\bullet)$. There results an algebraic $2$-ray $\mathfrak F(\cG_\bullet,\cJ_\bullet)$. We can now take its telescope, which yields an algebraic $1$-cube, that is a chain map between two chain complexes. Its Novikov completion is a chain map $SC^*(K) \to SC^*(K')$. The resulting map on cohomology is the restriction
$$\res^K_{K'} \fc SH^*(K) \to SH^*(K')\,.$$
In \cite{Varolgunes_2021_MV_and_relSH}, it is shown that this map is independent of all the choices.

Now, the complex underlying the relative symplectic cohomology $SH^*(K,K')$ of the pair $(K,K')$ is the cocone $\co(SC^*(K)\to SC^*(K'))$ of the above chain map. That is
$$SH^*(K,K')=H\big(\co(SC^*(K)\to SC^*(K'))\big)\,.$$
It is shown in \cite{DGPZ_2024_Symp_top_and_IVMs} that this cohomology module is independent of the choices in the sense that the cohomology for any other choice is canonically isomorphic to this one.

In order to define restriction maps, let $(K,K')\subset (L,L')$ be compact pairs in $M$, and let $\cG_\bullet$ be an acceleration $3$-ray for the inclusion $(K,K')\subset(L,L')$. Again, we pick a $3$-ray of almost complex structures $\cJ_\bullet$, apply the Floer functor to $(\cG_\bullet,\cJ_\bullet)$, apply telescope and then completion to the resulting algebraic $3$-ray to get a $2$-cube of the form
\begin{equation}\label{eq:diagram_rest_L_to_K}
  \xymatrix{SC^*(L) \ar[r] \ar[rd] \ar[d] & SC^*(L')\ar[d] \\ SC^*(K)\ar[r] & SC^*(K')}
\end{equation}
We can now take the cocone in the horizontal direction, which yields a chain map
$$\big(\co(SC^*(L)\to SC^*(L'))\big) \to \big(\co(SC^*(K)\to SC^*(K'))\big)\,,$$
that is
$$SC^*(L,L') \to SC^*(K,K')\,.$$
\begin{notation}\label{not:restriction_depending_on_data}
  For future use we denote this chain map by
  $$\res(\cG_\bullet,\cJ_\bullet)\fc SC^*(L,L') \to SC^*(K,K')$$
  to explicitly specify the data going into its definition.
\end{notation}
The induced map on cohomology is the desired restriction map
$$\res^{(L,L')}_{(K,K')} \fc SH^*(L,L') \to SH^*(K,K')\,.$$
Again, in \cite{DGPZ_2024_Symp_top_and_IVMs} it is shown that this map is independent of all the choices, and that it satisfies the properties announced in Section \ref{ss:pairs}.

\subsection{Proof of Proposition \ref{prop:excision}}\label{ss:proof_excision}

Recall that the proposition states that if $A,B\subset M$ are commuting compact sets, then the restriction map $SH^*(A\cup B,A) \to SH^*(B,A\cap B)$ is an isomorphism. It is proved in \cite{Varolgunes_2021_MV_and_relSH} that in case $A,B$ commute, they are \emph{in descent}, which means the $2$-cube
$$\xymatrix{SC^*(A\cup B) \ar[r] \ar[d] & SC^*(A) \ar[d] \\ SC^*(B) \ar[r] & SC^*(A\cap B)}$$
which we denote by $\cC$, is acyclic in the sense that its iterated cone is an acyclic chain complex. The cocone of this $2$-cube in the horizontal direction yields a chain map
$$r \fc \co(SC^*(A\cup B) \to SC^*(A)) \to \co (SC^*(B)\to SC^*(A\cap B))\,.$$
Since $\co$ is the cone shifted in degree by $-1$, the chain complex $\cone(r)$ is isomorphic, via a standard sign change, to the iterated cone of $\cC$, shifted by $-1$. It is therefore acyclic and thus $r$ is a quasi-isomorphism. By the definition of the relative symplectic cohomology of pairs, the source and target of $r$ are complexes computing $SH^*(A\cup B,A)$ and $SH^*(B,A\cap B)$, respectively. Moreover, the map induced by $r$ on cohomology is, up to sign, none other than the restriction map
$$SH^*(A\cup B,A) \to SH^*(B,A\cap B)\,,$$
which is therefore an isomorphism, as claimed.

\subsection{Proof of Theorem \ref{thm: invariance}}\label{ss:proof_invariance}

Recall the formulation of the theorem: Given compact pairs $(K,K')\subset(L,L')$ and a Hamiltonian isotopy $(\phi_t)_{t\in[0,1]}$ satisfying $(\phi_t(K),\phi_t(K'))\subset (L,L')$ for all $t\in[0,1]$, the following diagram commutes:
$$\xymatrix{& SH^*(L,L') \ar[ld]_{\res} \ar[rd]^{\res} \\ SH^*(K,K') \ar[rr]^{\phi_*} && SH^*(\phi(K),\phi(K'))}$$
that is
\begin{equation}\label{eq:identity_invariance}
  \phi_*\circ\res^{(L,L')}_{(K,K')} = \res^{(L,L')}_{(\phi(K),\phi(K'))}\,,
\end{equation}
where $\phi=\phi_1$, while $\phi_*$ is the so-called naturality isomorphism, defined below.

We need quite a bit of preparation for the proof, and so for the convenience of the reader, we will shortly present an overview thereof. Before we do so, let us introduce the central notion used in it, namely \emph{$\bs\phi$-parametrized cubes and rays}, where $\bs\phi$ is a Hamiltonian isotopy.

\begin{figure}[htbp]
	\centering
	\includegraphics[width=0.4\linewidth]{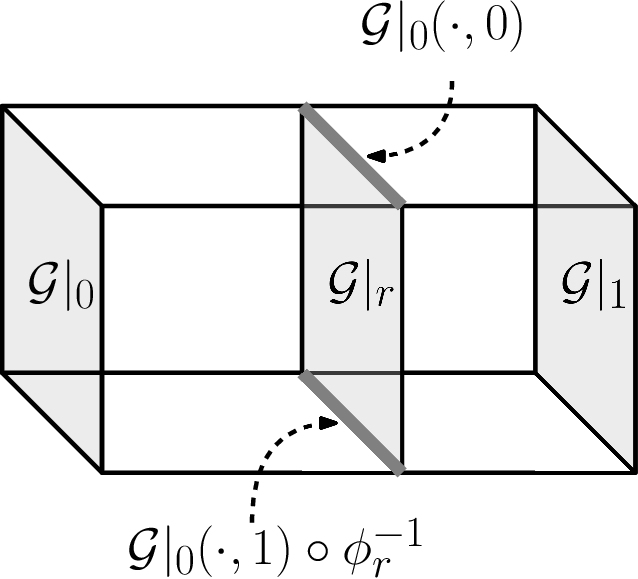}
	\caption{An illustration of a $\bs\phi$-parametrized Hamiltonian $3$-cube.}
	\label{fig: Hybrid Floer-alg cube 1}
\end{figure}

Let us denote the last coordinate of $[0,1]^n=[0,1]^{n-1}\times[0,1]$ by $r$ and view it as a parameter. An $n$-cube of Hamiltonians $\cG \fc [0,1]^n \to C^\infty(S^1\times M)$ induces for each $r$ the $(n-1)$-cube\footnote{Strictly speaking, $\cG|_r$ is not in general a Hamiltonian cube according to our definition, but rather just a smooth map, because it need not be constant near the vertices. However, we do not rely on the assumption of constancy for the $\cG|_r$ in what follows, and, on the other hand, for $\bs\phi$-parametrized cubes this assumption is implicitly a part of the definition.} $\cG|_r\fc[0,1]^{n-1}\to C^\infty(S^1\times M)$, called its \tb{$r$-slice}, which is given by $\cG|_r:=\cG(\cdot,r)$. This way, $\cG$ can be viewed as a smooth family $\cG=(\cG|_r)_{r\in[0,1]}$. In words, an $n$-cube is the same thing as a $[0,1]$-family of $(n-1)$-cubes. Similarly, if $\cG_\bullet=(\cG_i)_{i\in \N}$ is a Hamiltonian $n$-ray, we can consider its $r$-slice $\cG_\bullet|_r=(\cG_i|_r)_{i\in\N}$, which is an $(n-1)$-ray. We will need a special kind of such parametric cubes and rays.
\begin{notation}
  If $\cG \fc [0,1]^n\to C^\infty(S^1\times M)$ and $\phi\fc M \to M$ is a map, we define $\cG\circ\phi \fc [0,1]^n \to C^\infty(S^1\times M)$ by
  $$(\cG\circ\phi)(x)(t,m)=\cG(x)(t,\phi(m)) \quad\text{for }x\in[0,1]^n\,,t\in S^1\,,m\in M\,.$$
  If $\cG_\bullet=(\cG_i)_{i\in\N}$ is a sequence of maps as above, we let
  $$\cG_\bullet\circ\phi:=(\cG_i\circ\phi)_{i\in\N}\,.$$
\end{notation}
\begin{defin}\label{def:phi_parametrized_cubes}
  Let $n\geq 2$ and let us denote the last coordinate of $[0,1]^{n-1}=[0,1]^{n-2}\times[0,1]$ by $z$. Let $\bs\phi=(\phi_r)_{r\in[0,1]}$ be a Hamiltonian isotopy. A Hamiltonian $n$-cube $\cG$ is called \tb{$\bs\phi$-parametrized} if there is $\epsilon \in (0,\frac12)$ such that
  \begin{itemize}
    \item For $z \in [0,\epsilon]$ and $r\in[0,1]$, $\cG|_r(\cdot,z) = \cG|_0(\cdot,0)$;
    \item For $z \in [1-\epsilon,1]$ and $r\in[0,1]$, $\cG|_r(\cdot,z) = \cG|_0(\cdot,1)\circ\phi_r^{-1}$.
  \end{itemize}
  See Figure \ref{fig: Hybrid Floer-alg cube 1}. A $\bs\phi$-parametrized cube $\cG=(\cG|_r)_r$ is \tb{monotone} if each $\cG|_r$ is monotone. For $n\geq 3$, if $\cG_\bullet=(\cG_i)_{i\in \N}$ is a Hamiltonian $n$-ray, we call it \tb{$\bs\phi$-parametrized} (respectively, \tb{monotone}), if each constituent cube $\cG_i$ is (all for the same $\epsilon$).
\end{defin}

\noindent In words, for $z$ near $0$, $\cG|_r(\cdot,z)$ is independent of both $r,z$, and for $z$ near $1$, $\cG|_r(\cdot,z)$ is $z$-independent, and depends on $r$ in the specific manner indicated in the definition. The main result here is the following existence theorem for parametrized rays in the geometric context of Theorem \ref{thm: invariance}. It is proved in Section \ref{sss:proof_exists_parametrized_ray}.
\begin{thm}\label{thm:exists_parametrized_ray}
  Given compact pairs $(K,K')\subset(L,L')$ and a Hamiltonian isotopy $\bs\phi=(\phi_t)_{t\in[0,1]}$ such that $(\phi_t(K),\phi_t(K'))\subset(L,L')$ for all $t$, there exists a monotone $\bs\phi$-parametrized Hamiltonian $4$-ray $\cG_\bullet=(\cG_i)_{i\in\N}$ such that $\cG_\bullet|_0$ is an acceleration $3$-ray for the inclusion of pairs $(K,K')\subset (L,L')$.
\end{thm}
\begin{rem}\label{rem:phi_paramd_ray_slices_acc_rays}
  Note that since each $\cG_\bullet|_r$ is a monotone $3$-ray, it follows from the definition of a $\bs\phi$-parametrized ray that it is in fact an acceleration $3$-ray for the inclusion of pairs $(\phi_r(K),\phi_r(K'))\subset (L,L')$. This will be used in Section \ref{sss:proof_invariance} in the proof of Theorem \ref{thm: invariance}.
\end{rem}
\noindent Given this, we can now present our overview.
\paragraph*{Overview of the proof of Theorem \ref{thm: invariance}} We will construct a diagram of the form
\begin{equation}\label{eq:diag_overview_complexes}
  \xymatrix{SC^*(L,L') \ar@{}[rd]|-{\mathrm{I}} \ar[r]^{\id} \ar[d]_{\res_0} & SC^*(L,L') \ar[r]^{\id} \ar@{}[rd]|-{\mathrm{II}} \ar[d]_{\res_1} & SC^*(L,L') \ar[d]^{\res} \\ SC^*(K,K') \ar[r]_{\id} & SC^*(K,K')\ar[r]_-{\phi_*} & SC^*(\phi(K),\phi(K'))}
\end{equation}
consisting of two juxtaposed algebraic $2$-cubes, of which the cube II commutes while the cube I is a slit-like algebraic cube. The maps $\res_0,\res$ are usual restriction maps, while $\res_1$ is a certain map described below. Before we indicate how to obtain such a diagram, let us show that it does yield the desired identity \eqref{eq:identity_invariance}. Indeed, on passing to cohomology,  we obtain the following \emph{commutative} diagram (see Remark \ref{rem:slit_like_2_cube}):
$$\xymatrix{SH^*(L,L') \ar[r]^{\id} \ar[d]_{\res^{(L,L')}_{(K,K')}} & SH^*(L,L') \ar[r]^{\id} \ar[d] & SH^*(L,L') \ar[d]^{\res^{(L,L')}_{(\phi(K),\phi(K'))}} \\ SH^*(K,K') \ar[r]^{\id} & SH^*(K,K')\ar[r]^-{\phi_*} & SH^*(\phi(K),\phi(K'))}$$
and chasing around it yields
$$\phi_*\circ\res^{(L,L')}_{(K,K')} = \res^{(L,L')}_{(\phi(K),\phi(K'))}\,,\quad \text{as claimed.}$$

Briefly, the diagram \eqref{eq:diag_overview_complexes} is obtained in two steps, one for each square. For the first square we \emph{use the containment assumption, $(\phi_t(K),\phi_t(K'))\subset (L,L')$ for all $t$}, in order to be able to apply Theorem \ref{thm:exists_parametrized_ray}, which yields a suitable $\bs\phi$-parametrized ray. This $\bs\phi$-parametrized ray is then transformed into a $4$-ray of \emph{Hamiltonian $1$-forms} (see Section \ref{sss:Ham_1_forms}), from which the first square arises upon applying a generalization of the Floer functor, and then telescope, completion, and cocone. The second square is obtained as a particular case of Proposition \ref{prop:conjugated_cube_isomorphic}. This step is considerably easier, and \emph{does not rely on the containment assumption.}

Let us describe this in more detail. For any Hamiltonian $n$-cube $\cG$, $n\geq 1$, and any Hamiltonian isotopy $\boldsymbol{\psi}=(\psi_t)_{t\in[0,1]}$ which is stationary near $t=0,1$, in Section \ref{sss:conj_isotopy} we will describe a procedure, here called \emph{conjugation by $\boldsymbol\psi$}, which produces an $n$-cube of Hamiltonian $1$-forms $\cG^{\boldsymbol\psi}$. Likewise, given a cube of almost complex structures $\cI$, there is a naturally associated conjugated cube of almost complex structures $\cI^{\boldsymbol\psi}$. Again, in Section \ref{sss:conj_isotopy} we will describe how to generalize the Floer functor so that it is applicable to $(\cG^{\bs\psi},\cI^{\bs\psi})$, in case $\cG$ is monotone and nondegenerate at the vertices. It then yields the corresponding algebraic cube $\mathfrak F(\cG^{\bs\psi},\cI^{\bs\psi})$, which by Proposition \ref{prop:conjugated_cube_isomorphic} is isomorphic to $\mathfrak F(\cG,\cI)$ by an isomorphism whose vertex maps are identities and naturality isomorphisms $\psi_*$, where $\psi = \psi_1$.

We now take the isotopy $\boldsymbol\phi$ from the conditions of the theorem, and produce a family of isotopies $(\boldsymbol\phi_r)_{r\in[0,1]}$, all parametrized by $[0,1]$, with the property that all of them are stationary near $0,1$, and such that the time-$1$ map of $\boldsymbol\phi_r$ is $\phi_r$. Let $\cG_\bullet = (\cG_i)_{i\in\N}$ be the $\bs\phi$-parametrized $4$-ray from Theorem \ref{thm:exists_parametrized_ray}. We conjugate each cube $\cG_i|_r$ by the isotopy $\boldsymbol\phi_r$ to obtain a $3$-cube of Hamiltonian $1$-forms $\cK_i|_r = \big(\cG_i|_r\big)^{\boldsymbol\phi_r}$, which then assemble into a $4$-cube $\cK_i = (\cK_i|_r)_{r\in[0,1]}$, which is smooth and \emph{$z$- and $r$-independent for $z$ near $0,1$}, thanks to the fact that $\cG$ is $\bs\phi$-parametrized. These cubes then form the corresponding $4$-ray $\cK_\bullet=(\cK_i)_i$ of $1$-forms.

There exists a smooth map $\Xi :=\Xi_4 \fc [0,1]^4\to\Slit^4$ (Corollary \ref{cor:existence_properties_Xi_n}) collapsing the top and bottom faces and suitably compatible with the gradient flows of our Morse functions $\rho_n$. Thanks to the $z$- and $r$-independence of $\cK_\bullet$ for $z$ near $0,1$, there exists a smooth slit-like $4$-ray of Hamiltonian $1$-forms $\wt\cK_\bullet = (\wt\cK_i)_{i\in\N}$ with $\cK_i=\wt\cK_i\circ\Xi$ for all $i$. Now we pick a slit-like $4$-ray of almost complex structures $\wt\cJ_\bullet$ and apply the Floer functor to the pair $(\wt\cK_\bullet,\wt\cJ_\bullet)$. It yields a slit-like algebraic $4$-ray $\mathfrak{F}(\wt\cK_\bullet,\wt\cJ_\bullet)$. Taking its telescope and completing, we obtain an algebraic slit-like $3$-cube, whose cocone in the suitable direction is a slit-like algebraic $2$-cube which is the required square I in the diagram \eqref{eq:diag_overview_complexes}.

To identify the maps $\res_{0,1}$ in \eqref{eq:diag_overview_complexes}, set $\cJ_\bullet=(\wt\cJ_i\circ\Xi)_{i\in\N}$. Unraveling the definitions, it turns out (Lemma \ref{lem:computing_res_01}) that $\res_0 = \res(\cK_\bullet|_0,\cJ_\bullet|_0)$ (see Notation \ref{not:restriction_Ham_1_forms}), which is a usual restriction map, since $\cK_\bullet|_0=\cG_\bullet|_0$. On the other hand, $\res_1$ turns out to equal $\res(\cK_\bullet|_1,\cJ_\bullet|_1)$ (Lemma \ref{lem:computing_res_01} again). Since $\cK_\bullet|_1 = (\cG_\bullet|_1)^{\bs\phi_1}$ and the time-$1$ map of $\bs\phi_1$ is $\phi_1=\phi$, it follows from Proposition \ref{prop:conjugated_cube_isomorphic} that the map $\res(\cK_\bullet|_1,\cJ_\bullet|_1)$ is isomorphic to the restriction $\res\fc SC^*(L,L') \to SC^*(\phi(K),\phi(K'))$ by an isomorphism which involves the naturality isomorphism $\phi_*$. This yields the square II in \eqref{eq:diag_overview_complexes}. The desired identity $\phi_*\circ\res^{(L,L')}_{(K,K')} = \res^{(L,L')}_{(\phi(K),\phi(K'))}$ is obtained as a result. Let us now present the full proof.

\subsubsection{Naturality isomorphisms}\label{sss:naturality}

Any symplectomorphism $\phi$ of $M$ induces \tb{naturality isomorphisms}
\begin{equation}\label{eq:naturality_isos}
  \phi_* \fc SH^*(K) \to SH^*(\phi(K))\,,\quad \phi_* \fc SH^*(K,K') \to SH^*(\phi(K),\phi(K'))
\end{equation}
for $K\in\cC(M)$ and $(K,K') \in \cC\cP(M)$, compatible with restrictions. They are defined as follows. If $H$ is a nondegenerate Hamiltonian and $J$ an $S^1$-family of almost complex structures, and we have fixed extra data as above, then we have the module isomorphism
$$\phi_* \fc CF^*(H) \to CF^*(H\circ\phi^{-1})\quad\text{given by }x\mapsto \phi\circ x\,.$$
We claim that there are equalities of virtual counts $n(x_-,x_+,H,J,A)=n(\phi\circ x_-,\phi\circ x_+,H\circ\phi^{-1},\phi_*J,\phi_*A)$, which then yield an isomorphism of chain complexes
$$(CF^*(H),d) \cong (CF^*(H\circ\phi^{-1}),d)\,.$$
To see why this is so, let us very briefly recall how these virtual counts come about. One has the moduli space of solutions of the Floer equation $\ol\cM(x_-,x_+,H,J,A)$, suitably compactified (see the discussion in Section \ref{sss:Floer_fctr}), on which there is constructed an implicit atlas $\cA$ (\cite{Pardon_2016_alg_vir_counting}). This implicit atlas is canonically associated to the data and is independent of all choices, with the exception of a certain lifting map \cite[Definition 7.5.3]{Pardon_2016_alg_vir_counting}. The composition with $\phi$, $u\mapsto \phi\circ u$, applied to Floer cylinders and holomorphic spheres, induces a homeomorphism
$$\ol\cM(x_-,x_+,H,J,A)\cong \ol\cM(\phi\circ x_-,\phi\circ x_+,H\circ\phi^{-1},\phi_*J,\phi_*A)\,,$$
but also induces a natural bijection between $\cA$ and the implicit atlas for the latter moduli space, thanks to its functorial nature; the choice of a lifting map for $(H,J)$ yields, through $\phi$, a choice of a lifting map for $(H\circ\phi^{-1},\phi_*J)$. There then results an equality of virtual fundamental classes and therefore of the virtual counts. The topological energy of the configurations whose virtual counts enter the definition of $d$ is obviously preserved, leading to the above chain isomorphism.

The same procedure can be applied to cubes of Hamiltonians and almost complex structures: given a Hamiltonian cube $\cG$ and a cube $\cI$ of almost complex structures, we obtain the pushed-forward cubes $\phi_*\cG=\cG\circ\phi^{-1}$, $\phi_*\cI$. If $\cG$ has nondegenerate vertices and is monotone, then the same holds for $\phi_*\cG$. For a face $F$ with initial and terminal vertices $v,v'$, respectively, periodic orbits $x\in\cP^\circ(\cG(v))$, $x' \in \cP^\circ(\cG(v'))$, the number $n_F(x,x',\cG,\cI,A)$ is, as indicated above, a virtual count of pairs $(\gamma,u)$ consisting of a negative gradient line $\gamma \fc \R\to F$ asymptotic to $v,v'$ and a smooth cylinder $u$ satisfying
$$\partial_su-\cI(\gamma(s))(u)(\partial_tu-X_{\cG(\gamma(s))}(u))=0\,,$$
with some $\cI$-holomorphic spheres attached. The composition with $\phi$ maps such objects to objects $(\gamma,\phi\circ u)$ which now satisfy
$$\partial_s(\phi\circ u) - \phi_*\cI(\gamma(s))(\phi\circ u)(\partial_t(\phi\circ u) - X_{\phi_*\cG(\gamma(s))}(\phi\circ u))=0\,,$$
with the holomorphic spheres being transferred by $\phi$ as well, and this induces a homeomorphism of the corresponding moduli spaces. Again, the implicit atlases on both sides are in a natural bijection induced by $\phi$, which yields an equality of virtual counts
$$n_F(x,x',\cG,\cI,A)=n_F(\phi\circ x,\phi\circ x',\phi_*\cG,\phi_*\cI,\phi_*A)\,.$$
It follows that the face maps of the cubes $\mathfrak{F}(\cG,\cI), \mathfrak{F}(\phi_*\cG,\phi_*\cI)$ are intertwined by the isomorphisms $\phi_*$ at each vertex, that is there is an isomorphism of algebraic cubes $\mathfrak{F}(\cG,\cI)\to \mathfrak{F}(\phi_*\cG,\phi_*\cI)$ given by $\phi_*$ at each vertex.

Let now $(K^\star)_\star\fc \{0,1\}^n \to\cC(M)$ be a system of sets, and fix an acceleration ray $\cG_\bullet$ for it. Let $\cI_\bullet$ be an $(n+1)$-ray of almost complex structures. We can then form the pushforward rays $\phi_*\cG_\bullet$, $\phi_*\cI_\bullet$ in an obvious manner, where $\phi_*\cG_\bullet$ is now an acceleration ray for the system $(\phi(K^\star))_\star$. There result algebraic $(n+1)$-rays $\mathfrak{F}(\cG_\bullet,\cI_\bullet)$, $\mathfrak F(\phi_*\cG_\bullet,\phi_*\cI_\bullet)$. Thanks to the above discussion, these rays are cube-wise isomorphic via isomorphisms being given by $\phi_*$ at each vertex. Since the telescope construction is functorial, it follows that the telescopes of $\mathfrak{F}(\cG_\bullet,\cI_\bullet)$, $\mathfrak F(\phi_*\cG_\bullet,\phi_*\cI_\bullet)$ are likewise isomorphic, the isomorphism being induced from the $\phi_*$. Applying completion to the resulting algebraic $n$-cubes, we obtain cubes with the $SC^*(K^\star)$ and $SC^*(\phi(K^\star))$ at the vertices, and moreover, thanks to the functoriality of completion, these cubes are isomorphic where the vertex isomorphisms are induced from $\phi_*$, and which we will still denote by $\phi_*$ by abuse of notation.

The desired naturality isomorphisms \eqref{eq:naturality_isos} are now obtained as particular cases of this construction. Namely, if we only have one compact set, the construction, applied to the system of sets consisting only of $K$ yields a chain isomorphism $SC^*(K) \to SC^*(\phi(K))$, which induces the naturality isomorphism
$$\phi_* \fc SH^*(K) \to SH^*(\phi(K))\quad\text{on passing to cohomology.}$$
One shows that this is independent of the choices. If we have a compact pair $(K,K')$, the above construction yields two chain maps $SC^*(K) \to SC^*(K')$, $SC^*(\phi(K))\to SC^*(\phi(K'))$, which appear in the following commutative diagram:
$$\xymatrix{SC^*(K) \ar[r]^{\phi_*} \ar[d]_{\res} & SC^*(\phi(K)) \ar[d]^{\res} \\ SC^*(K') \ar[r]^{\phi_*} & SC^*(\phi(K'))}$$
Passing to cohomology, we see that the naturality isomorphisms respect restriction. If in the last diagram we instead take the cocone in the vertical direction, we will obtain a chain isomorphism
$$\phi_* \fc SC^*(K,K') \to SC^*(\phi(K),\phi(K'))\,,$$
which on cohomology yields the other naturality isomorphism in \eqref{eq:naturality_isos}. Similarly, using an inclusion of pairs $(K,K')\subset (L,L')$, we see that naturality isomorphisms commute with restrictions for pairs.

\subsubsection{Hamiltonian $1$-forms}\label{sss:Ham_1_forms}

We need a slight generalization of the Floer functor, which is applicable to cubes of \emph{Hamiltonian $1$-forms} on the cylinder $\R\times S^1$. Such a form $K$ is simply a $1$-form on $\R\times S^1$ with values in $C^\infty(M)$, and can be written as
$$K(s,t) = F(s,t)\otimes ds + G(s,t)\otimes dt\quad\text{for }(s,t)\in\R\times S^1\,,$$
where $F(s,t),G(s,t)$ are smooth functions on $M$, smoothly dependent on $(s,t)$. Equivalently, $K$ can be represented by a $1$-form on the manifold $\R\times S^1\times M$ which vanishes on vectors tangent to the fibers of the projection $\R\times S^1\times M \to \R\times S^1$. To such a Hamiltonian $1$-form there is associated the corresponding $1$-form with values in Hamiltonian vector fields
$$X_K(s,t) = X_{F(s,t)}\otimes ds + X_{G(s,t)}\otimes dt\,.$$
Given a smooth family of almost complex structures $J(s,t)$, we have the corresponding generalized Floer equation for $u \in C^\infty(\R\times S^1,M)$:
$$\big(d_{(s,t)}u - X_K(s,t)(u(s,t))\big)^{1,0}_{J(s,t)}=0\,.$$
Here for a linear map $\eta\fc T_{(s,t)}(\R\times S^1) \to T_{u(s,t)}M$ we use its $J(s,t)$-complex linear part
$$\eta^{1,0}_{J(s,t)} = \eta - J(s,t)(u(s,t))\eta\circ i\,.$$
Substituting $\partial_s$ into the above equation, we arrive at the following coordinate expression:
$$\partial_su - X_F(u) - J(u)(\partial_tu - X_G(u)) = 0\,.$$
This reverts to the usual Floer equation whenever $F\equiv 0$.

Even though a more general treatment certainly makes sense, for our purposes we will only consider Hamiltonian $1$-forms which are \emph{$s$-independent} in the sense that they have the form
$$K(s,t) = F(t)\otimes ds + G(t)\otimes dt\,.$$
Just like for Hamiltonian functions, we have the notion of an \tb{$n$-cube of Hamiltonian $1$-forms}, which is a smooth map $\cK \fc [0,1]^n\to\Omega^1(\R\times S^1,C^\infty(M))$, constant in a neighborhood of each vertex. Since we will only consider $s$-independent $1$-forms, such a cube can be uniquely written as
$$\cK = \cF \otimes ds + \cG\otimes dt\quad\text{for Hamiltonian cubes }\cF\,,\cG\,.$$
The main point here is that, similarly to Hamiltonian cubes, cubes of Hamiltonian $1$-forms, satisfying suitable nondegeneracy and compatibility conditions, paired with cubes of almost complex structures, give rise to algebraic cubes by means of the Floer functor.

To define the Floer functor for such a $1$-form $\cK=\cF\otimes ds + \cG\otimes dt$, we need to require that $\cF\equiv 0$ in a neighborhood of each vertex, and that $\cG$ be nondegenerate at each vertex. The description of the functor is nearly identical to the case of Hamiltonian cubes. Namely, in addition to $\cK$, let us also fix a cube of almost complex structures $\cI$. For a face $F\subset [0,1]^n$ with initial and terminal vertices $v,v'$, respectively, and $x \in \cP^\circ(\cG(v))$, $x'\in\cP^\circ(\cG(v'))$, we have numbers $n_F(x,x',\cK,\cI,A)\in\Q$ which are interpreted as a virtual count of objects consisting of pairs $(\gamma,u)$, where $\gamma \fc \R \to F$ is a negative gradient line of $\rho_n$ from $v$ to $v'$, $u$ is a smooth cylinder, to which there are attached holomorphic spheres, where now $u$ solves the suitably generalized equation
\begin{equation}\label{eq:Floer_eq_Ham_forms}
  (du - X_{\cK\circ\gamma}(u))^{1,0}_{\cI\circ\gamma}=0\,,
\end{equation}
or in coordinates
$$\partial_su - X_{\cF(\gamma(s))}(u) - \cI(\gamma(s))(\partial_tu - X_{\cG(\gamma(s))}(u)) = 0\,.$$
Here again pairs $(\gamma,u)$ obtained from one another by a simultaneous $s$-shift are considered the same object.
\begin{rem}\label{rem:map_at_vertex_Ham_form_Floer_diff}
  Note that if $F=\{v\}$ consists solely of a single vertex, due to the condition that $\cF\equiv 0$ near $v$, the corresponding face map is simply the Floer differential on $CF^*(\cG(v))$.
\end{rem}
The rest of the definition is a \emph{verbatim} repetition of the construction in \ref{sss:Floer_fctr}, except we need to require that \emph{the topological energy be nonnegative} for each such solution. The topological energy in question is expressed by the same formula \eqref{eq:top_energy_cube}. This imposes a restriction on the Hamiltonian form $\cK$, which in cases relevant to us is satisfied.
\begin{rem}
  Although Hamiltonian $1$-forms are not directly considered in \cite{Pardon_2016_alg_vir_counting,Pardon_2019_contact}, the virtual cycle formalism there carries over \emph{mutatis mutandis}. Indeed, it is based on flow diagrams, which are identical in both cases, and on topological gluing arguments, which again apply to both cases.
\end{rem}

\begin{rem}\label{rem:Ham_forms_acc_data_induce_restrictions}
  We will close this section by considering a special case of the application of the Floer functor to cubes of Hamiltonian forms, which will play a crucial role in the proof of Theorem \ref{thm: invariance}. Let $(K,K')\subset (L,L')$ be an inclusion of compact pairs. Let $\cG_\bullet^{(K,K')}$ and $\cG_\bullet^{(L,L')}$ be acceleration $2$-rays for the pairs $(K,K')$, $(L,L')$, respectively. Number the coordinates of $[0,1]^3$ as follows: the first coordinate is used for the semi-infinite direction of rays, the second coordinate goes from $K$ to $K'$ and from $L$ to $L'$; name the third coordinate $z$. Assume that $\cK_\bullet$ is a $3$-ray of Hamiltonian $1$-forms such that for some $\epsilon >0$, its restriction to $\{z\in[0,\epsilon]\cup[1-\epsilon,1]\}$ is independent of $z$, and for $z=0$ it coincides with $\cG_\bullet^{(L,L')}$, as a function of the first two coordinates of $[0,1]^3$, while for $z=1$ it coincides with $\cG_\bullet^{(K,K')}$. Let us now pick a $3$-ray $\cJ_\bullet$ of almost complex structures. Assuming as well that all the relevant solutions have nonnegative topological energy as in the discussion above, we can apply the Floer functor to $(\cK_\bullet,\cJ_\bullet)$ cube-wise to obtain an algebraic $3$-ray $\mathfrak{F}(\cK_\bullet,\cJ_\bullet)$. Thanks to Remark \ref{rem:map_at_vertex_Ham_form_Floer_diff}, the vertices of this algebraic ray are the Floer complexes of the acceleration data coming from $\cG_\bullet^{(L,L')}$, $\cG_\bullet^{(K,K')}$. Moreover, since the restrictions of $\cK_\bullet$ to the corner rays yield actual acceleration rays for the sets $L,L',K,K'$, it follows that the algebraic $1$-rays corresponding to these corners yield, upon taking the telescope and completion, complexes $SC^*(A)$, $A=L,L',K,K'$. Thus, the completion of the telescope of $\mathfrak{F}(\cK_\bullet,\cJ_\bullet)$ is a $2$-cube of the form
  $$\xymatrix{SC^*(L) \ar[r] \ar[rd] \ar[d] & SC^*(L')\ar[d] \\ SC^*(K)\ar[r] & SC^*(K')}$$
  Since $\cK_\bullet$ is required to coincide with actual acceleration $2$-rays near $z=0,1$, the horizontal arrows in the last diagram are in fact genuine restriction maps, which means that when we take the cocone of the diagram in the horizontal direction, we obtain a map $SC^*(L,L') \to SC^*(K,K')$.
\end{rem}
\noindent Since $\cK_\bullet$ is not necessarily a monotone Hamiltonian ray, this map will not in general induce the restriction map $\res^{(L,L')}_{(K,K')}$ on cohomology. Nevertheless, for convenience's sake, we will employ
\begin{notation}\label{not:restriction_Ham_1_forms}
  The map $SC^*(L,L') \to SC^*(K,K')$ described in the preceding paragraph will be denoted by $\res(\cK_\bullet,\cJ_\bullet)$, as will be the induced map on cohomology.
\end{notation}
\noindent This is in analogy to Notation \ref{not:restriction_depending_on_data}.

\subsubsection{Conjugation by an isotopy}\label{sss:conj_isotopy}

The appearance of Hamiltonian $1$-forms is explained as follows. We will shortly describe a procedure, which we refer to as \emph{conjugation}, which takes a Hamiltonian cube $\cG$ and a Hamiltonian isotopy $\bs\phi$ and produces a cube of $s$-independent Hamiltonian $1$-forms $\cG^{\bs\phi}$. A related construction produces from a cube of almost complex structures $\cI$ another such cube $\cI^{\bs\phi}$. The result here is Proposition \ref{prop:conjugated_cube_isomorphic} which says that the algebraic cubes $\mathfrak{F}(\cG,\cI)$, $\mathfrak{F}(\cG^{\bs\phi},\cI^{\bs\phi})$ are isomorphic by an isomorphism whose vertex maps are identity maps and the naturality isomorphisms $\phi_*$, where $\phi$ is the time-$1$ map of $\bs\phi$. Let us define these precisely.

\begin{defin}
  Let $\cG \fc [0,1]^n \to C^\infty(S^1\times M)$ be a Hamiltonian cube and $\cI \fc [0,1]^n\to C^\infty(S^1,\cJ(M,\omega))$ be a cube of almost complex structures. Let $\bs\phi=(\phi_\sigma)_{\sigma\in[0,1]}$ be a Hamiltonian isotopy, which is stationary near $\sigma=0,1$, and let $H\fc [0,1]\times M \to \R$ be its normalized generating Hamiltonian, that is $\phi_\sigma=\phi^\sigma_H$. The \tb{conjugation of $\cG$ by $\bs\phi$} is the $n$-cube of Hamiltonian $1$-forms $\cG^{\bs\phi}$ defined by
  $$\cG^{\bs\phi}(x) = \nabla_n\rho_n(x)H_{x_n}\circ\phi_{x_n}\otimes ds + \cG(x)\circ\phi_{x_n}\otimes dt\,,$$
  where $\nabla_n\rho_n$ is the $n$-component of $\nabla\rho_n$, and where $H_{x_n},\phi_{x_n}$ mean $H,\bs\phi$ \emph{at time $x_n$}. Similarly, the \tb{conjugation of $\cI$ by $\bs\phi$} is the $n$-cube of almost complex structures $\cI^{\bs\phi}$ defined by
  $$\cI^{\bs\phi}(x) = \phi_{x_n}^*\cI(x) = d\phi_{x_n}^{-1}\circ\cI(x)\circ d\phi_{x_n}\,.$$
\end{defin}
\noindent We similarly define the conjugation of a Hamiltonian $n$-ray $\cG_\bullet$ as the $n$-ray $\cG_\bullet^{\bs\phi} = (\cG_i^{\bs\phi})_{i\in\N}$, and similarly for an $n$-ray $\cI_\bullet$ of almost complex structures we have the conjugated $n$-ray $\cI_\bullet^{\bs\phi} = (\cI_i^{\bs\phi})_{i\in\N}$.

\begin{rem}\phantom{M}
\begin{itemize}
  \item Note that $\cG^{\bs\phi}$ is indeed a cube of $s$-independent Hamiltonian $1$-forms, since it only depends on $t$ through $\cG(x)$. The dependence on $s$ reappears in the associated Floer equation when a negative gradient line is inserted into it.
  \item Note that since $H$ is assumed normalized, it is uniquely determined by $\bs\phi$, and thus the above definition only depends on $\bs\phi$.
\end{itemize}
\end{rem}
The main point here is as follows.
\begin{prop}\label{prop:conjugated_cube_isomorphic}
  Let $(\cG,\cI)$ be a pair consisting of a monotone Hamiltonian cube with nondegenerate vertices, and a cube of almost complex structures. Then the algebraic cubes $\mathfrak F(\cG,\cI)$, $\mathfrak F(\cG^{\bs\phi},\cI^{\bs\phi})$ are isomorphic by an isomorphism $\mathfrak F(\cG^{\bs\phi},\cI^{\bs\phi}) \to \mathfrak F(\cG,\cI) $ whose vertex maps are identities in the top layer $x_n=0$ and the naturality isomorphisms $\phi_*$ in the bottom layer $x_n=1$, where $\phi=\phi_1$.
\end{prop}
\begin{proof}
  We will only sketch the proof, indicating the relevant points from Pardon's construction of the virtual fundamental classes \cite{Pardon_2016_alg_vir_counting}. The assertion boils down to the following: fix a face $F\subset[0,1]^n$ with initial and terminal vertices $v,v'$, respectively, periodic orbits $x\in\cP^\circ(\cG(v))$, $x'\in\cP^\circ(\cG(v'))$, and a class $A$ relative to $x,x'$; we then must show that
  $$n_F(x,x',\cG,\cI,A) = n_F(y,y',\cG^{\bs\phi},\cI^{\bs\phi},\bs\phi^*A)\,,$$
  where $y = x$ if $x_n(v)=0$ and $y=\phi^{-1}\circ x$ if $x_n(v)=1$, and similarly for $y'$. The class $\bs\phi^*A$ is defined as follows: if $x_n(v')=0$, $\bs\phi^*A=A$; if $x_n(v)=1$, $\bs\phi^*A = \phi^*A$; if $x_n(v)=0$, $x_n(v')=1$, then $\bs\phi^*A$ is the class of $[0,1]\times S^1 \ni (\sigma,\tau) \mapsto \phi_\sigma^{-1}(w(\sigma,\tau))$ for some representative $w \fc [0,1]\times S^1\to M$ of $A$.

  We have the following three cases.

  \tb{Case I: $F\subset\{x_n=0\}$.} Here the moduli spaces on both sides (conjugated and unconjugated) only involve the data
    $$(\cG|_{x_n=0},\cI|_{x_n=0}) \quad \text{and}\quad (\cG^{\bs\phi}|_{x_n=0},\cI^{\bs\phi}|_{x_n=0}) = (\cG|_{x_n=0},\cI|_{x_n=0})\,,$$
    meaning $y=x$, $y'=x'$, $\bs\phi^*A=A$. Thus, we have the same data on both sides, and so $n_F(y,y',\cG^{\bs\phi},\cI^{\bs\phi},\bs\phi^*A) = n_F(x,x',\cG,\cI,A)$.

  \tb{Case II: $F\subset \{x_n=1\}$.} Here we need only consider the data
    $$(\cG^{\bs\phi}|_{x_n=1},\cI^{\bs\phi}|_{x_n=1}) = (\phi^*\cG|_{x_n=1},\phi^*\cI|_{x_n=1})\,,$$
    that is $y=\phi^{-1}\circ x$, $y'=\phi^{-1}\circ x'$, $\bs\phi^*A=\phi^*A$. Here the equality $n_F(y,y',\cG^{\bs\phi},\cI^{\bs\phi},\bs\phi^*A) = n_F(x,x',\cG,\cI,A)$ is a particular case of what goes on into the construction of the naturality isomorphism, see Section \ref{sss:naturality}.

  \tb{Case III: $x_n(v) = 0$ and $x_n(v')=1$.} Here we will describe how to transform an element in the moduli space for the data $(\cG,\cI)$ to an element in the moduli space for $(\cG^{\bs\phi}, \cI^{\bs\phi})$. Recall that $n_F(x,x',\cG,\cI,A)$ is a virtual count of objects consisting of a pair $(\gamma,u)$, where $\gamma \fc \R \to F$ is a negative gradient line of $\rho_n$ from $v$ to $v'$, $u$ is a smooth cylinder solving
  \begin{equation}\label{eq:Floer_eq_along_gamma}
    \partial_su - \cI(\gamma(s))(u)(\partial_tu - X_{\cG(\gamma(s))}(u))=0\,,
  \end{equation}
  and where there are holomorphic spheres attached to $u$. Again, pairs $(\gamma,u)$ differing by simultaneous $s$-shifts are considered the same. Let us define $u^{\bs\phi}$ by
  $$u^{\bs\phi}(s,t) = \phi_{\gamma_n(s)}^{-1}(u(s,t))\,,\quad\text{where }\gamma(s)=(\gamma_1(s),\dots,\gamma_n(s))\,.$$
  We have the following key property, proved below.
  \begin{lemma}\label{lem:homeomorphism_moduli_spaces_conjugated}
    The map $(\gamma,u)\mapsto (\gamma,u^{\bs\phi})$ is a homeomorphism between the space $\{(\gamma,u)\}$ of solutions of \eqref{eq:Floer_eq_along_gamma} and the space $\{(\gamma,u^{\bs\phi})\}$ of solutions of
    \begin{equation}\label{eq:Floer_eq_conjugated}
      \big(du^{\bs\phi} - X_{\cG^{\bs\phi}\circ\gamma}(u^{\bs\phi})\big)^{1,0}_{\cI^{\bs\phi}\circ\gamma}=0\,.
    \end{equation}
  \end{lemma}
  \noindent Of course, if we have a sphere $w$ attached to $u$ at $u(s,t)$, then it is $\cI(\gamma(s))_t$-holomorphic, and it is then mapped to $\phi_{\gamma_n(s)}^{-1}\circ w$, which is now attached to $u^{\bs\phi}$ at $u^{\bs\phi}(s,t)$ and is $\cI^{\bs\phi}(\gamma(s))_t$-holomorphic. It follows that this operation induces a homeomorphism between the relevant moduli spaces of solutions. It also induces isomorphisms between the linearized operators at $u$ and $u^{\bs\phi}$, and therefore a bijection between the implicit atlases on both sides, the atlases being constructed out of stabilizations of these linearizations. Ultimately, this results in the sought-for equality of virtual counts
  $$n_F(y,y',\cG^{\bs\phi},\cI^{\bs\phi},\bs\phi^*A) = n_F(x,x',\cG,\cI,A)\,.$$
  It only remains to note the equality of topological energies
  $$E_{\mathrm{top}}(x,x',\cG,A) = E_{\mathrm{top}}(y,y',\cG^{\bs\phi},\bs\phi^*A)\,,$$
  which follows from the Stokes formula. This establishes the desired isomorphism of cubes.
\end{proof}

\begin{proof}[Proof of Lemma \ref{lem:homeomorphism_moduli_spaces_conjugated}]
  Assume $u$ solves
  $$\partial_su - \cI(\gamma(s))(u)(\partial_tu - X_{\cG(\gamma(s))}(u))=0$$
  and $u(s,t) = \phi_{\gamma_n(s)}(u^{\bs\phi}(s,t))$. It follows that
  $$\partial_su = \gamma_n'(s)X_{H_{\gamma_n(s)}}(\phi_{\gamma_n(s)}(u^{\bs\phi})) + d\phi_{\gamma_n(s)}(\partial_su^{\bs\phi})\,,\quad\partial_tu=d\phi_{\gamma_n(s)}(\partial_tu^{\bs\phi})\,.$$
  Therefore
  \begin{align*}
    0 &= \partial_su - \cI(\gamma(s))(u)(\partial_tu - X_{\cG(\gamma(s))}(u)) \\
      &= d\phi_{\gamma_n(s)}(\partial_su^{\bs\phi}) + \gamma_n'(s)X_{H_{\gamma_n(s)}}(\phi_{\gamma_n(s)}(u^{\bs\phi}))\\ &\qquad - \cI(\gamma(s))(\phi_{\gamma_n(s)}(u^{\bs\phi}))\big(d\phi_{\gamma_n(s)}(\partial_t u^{\bs\phi}) - X_{\cG(\gamma(s))}(\phi_{\gamma_n(s)}(u^{\bs\phi}))\big) \\
      &= d\phi_{\gamma_n(s)}\Big(\partial_su^{\bs\phi} + \gamma_n'(s)d\phi_{\gamma_n(s)}^{-1}X_{H_{\gamma_n(s)}}(\phi_{\gamma_n(s)}(u^{\bs\phi})) \\
      &\qquad - \big(d\phi_{\gamma_n(s)}^{-1}\cI(\gamma(s))d\phi_{\gamma_n(s)}\big)(u^{\bs\phi})(\partial_tu^{\bs\phi} - d\phi_{\gamma_n(s)}^{-1}X_{\cG(\gamma(s))}(\phi_{\gamma_n(s)}(u^{\bs\phi}))\Big)\,.
  \end{align*}
  Now, for a symplectomorphism $\psi$ and a Hamiltonian $F$ on $M$ we have
  $$X_{F\circ\psi}=X_{\psi^*F}=\psi^*X_F = d\psi^{-1} X_F\circ \psi\,,$$
  whence
  $$d\phi_{\gamma_n(s)}^{-1}X_{H_{\gamma_n(s)}}(\phi_{\gamma_n(s)}(u^{\bs\phi})) = X_{H_{\gamma_n(s)}\circ\phi_{\gamma_n(s)}}(u^{\bs\phi})\quad\text{and}$$
  $$d\phi_{\gamma_n(s)}^{-1}X_{\cG(\gamma(s))}(\phi_{\gamma_n(s)}(u^{\bs\phi})) = X_{\cG(\gamma(s))\circ\phi_{\gamma_n(s)}}(u^{\bs\phi})\,.$$
  Also, $d\phi_{\gamma_n(s)}^{-1}\cI(\gamma(s))d\phi_{\gamma_n(s)} = \phi_{\gamma_n(s)}^*\cI(\gamma(s))$, and
  $$\gamma_n'(s) = -\nabla_n\rho_n(\gamma(s))\,,$$
  thus the above equation for $u^{\bs\phi}$ becomes
  $$\partial_su^{\bs\phi} - X_{\nabla_n\rho_n(\gamma(s))H_{\gamma_n(s)}\circ\phi_{\gamma_n(s)}}(u^{\bs\phi}) - \phi_{\gamma_n(s)}^*\cI(\gamma(s))(u^{\bs\phi})\big(\partial_tu^{\bs\phi} - X_{\cG(\gamma(s))\circ\phi_{\gamma_n(s)}}(u^{\bs\phi})\big) = 0\,.$$
  This is precisely the equation
  $$\big(du^{\bs\phi} - X_{\cG^{\bs\phi}\circ\gamma}(u^{\bs\phi})\big)^{1,0}_{\cI^{\bs\phi}\circ\gamma}=0$$
  for
  $$\cG^{\bs\phi}(x) = \nabla_n\rho(x)H_{x_n}\circ\phi(x_n) \otimes ds + \cG(x)\circ\phi_{x_n}\otimes dt\,.$$
  Thus, $(\gamma,u)\mapsto (\gamma,u^{\bs\phi})$ indeed maps solutions of \eqref{eq:Floer_eq_along_gamma} to solutions of \eqref{eq:Floer_eq_conjugated}. It is also clearly bijective, the reverse bijection being given by $v\mapsto v^{\bs\phi^{-1}}$, where $\bs\phi^{-1}$ is the isotopy $\sigma\mapsto \phi_\sigma^{-1}$. It also clearly is a homeomorphism between the corresponding moduli spaces.
\end{proof}

\subsubsection{From cubes to slits}\label{sss:cubes_to_slits}

In this final preparatory step we show how a $\bs\phi$-parametrized Hamiltonian cube can be transformed into an $n$-slit of Hamiltonian $1$-forms with special properties.  Let $\cG$ be any $\bs\phi$-parametrized $n$-cube, and let $\epsilon>0$ be as in Definition \ref{def:phi_parametrized_cubes}. Fix a family of isotopies $(\bs\phi_\sigma)_{\sigma\in[0,1]}$ such that
\begin{itemize}
  \item For each $\sigma$, $\bs\phi_\sigma(t)$ is the identity for $t\in[0,\epsilon]$ and equals $\phi_\sigma$ for $t\in[1-\epsilon,1]$;
  \item For each $\sigma$, $\bs\phi_\sigma$ is a time-reparametrization of the isotopy $(\phi_t)_{t\in[0,\sigma]}$.
\end{itemize}
Note that in particular the time-$1$ map of $\bs\phi_\sigma$ is precisely $\phi_\sigma$, $\bs\phi_0$ is the constant isotopy at the identity, and $\bs\phi_1$ is a reparametrization of $\bs\phi$.
\begin{rem}
  If $H$ is a generating Hamiltonian for $(\phi_\sigma)_\sigma$, such a family of isotopies can be produced as follows. Fix a nondecreasing smooth function $p\fc[0,1]\to[0,1]$ with $p|_{[0,\epsilon]}\equiv 0$ and $p|_{[1-\epsilon,1]}\equiv 1$. Define $H^\sigma \fc [0,1]\times M\to \R$ by $H^\sigma_t=\sigma p'(t) H_{\sigma p(t)}$ and set $\phi_\sigma(t)=\phi^t_{H^\sigma}$.
\end{rem}
\noindent We can now define a cube of Hamiltonian $1$-forms $\cG^{(\bs\phi_\sigma)_\sigma}$ by setting
$$\cG^{(\bs\phi_\sigma)_\sigma}:=\big(\cG^{(\bs\phi_\sigma)_\sigma}|_r\big)_{r\in[0,1]}\quad\text{where}\quad \cG^{(\bs\phi_\sigma)_\sigma}|_r:=\big(\cG|_r\big)^{\bs\phi_r}\,.$$
Note that this is indeed a smooth cube of $1$-forms. Here is the chief property of this cube, where points in $[0,1]^n$ are denoted by $(x,z,r)$ with $x\in[0,1]^{n-2}$, $z,r\in[0,1]$:
\begin{lemma}\label{lem:cube_of_forms_conjugated_paramd_indep_z_r}
  The cube $\cG^{(\bs\phi_\sigma)_\sigma}$ is independent of $z$ and $r$ for $z\in[0,\epsilon] \cup [1-\epsilon,1]$.
\end{lemma}
\begin{proof}
  Let $H^r=(H^r_z)_{z\in[0,1]}$ be the normalized generating Hamiltonian for the isotopy $\bs\phi_r = (\bs\phi_r(z))_{z\in[0,1]}$, that is $\bs\phi_r(z)=\phi^z_{H^r}$. Since $\bs\phi_r$ is stationary on $[0,\epsilon]\cup[1-\epsilon,1]$, it follows that $H^r_t\equiv 0$ for $t\in[0,\epsilon]\cup[1-\epsilon,1]$. Recall that for $(x,z)\in[0,1]^{n-1}$
  $$(\cG|_r)^{\bs\phi_r}(x,z) = \nabla_{n-1}\rho_{n-1}(x,z)H^r_{z}\circ\bs\phi_r(z) \otimes ds + \cG(x,z,r)\circ\bs\phi_r(z)\otimes dt\,.$$
  For $z\in[0,\epsilon]$ we have $H^r_z\equiv 0$ and $\bs\phi_r(z)\equiv \id$, therefore
  $$(\cG|_r)^{\bs\phi_r}(x,z) = \cG(x,z,r)\circ\id\otimes dt = \cG(x,z,r)\otimes dt\,.$$
  Since $\cG$ is $\bs\phi$-parametrized, it follows from the definition that $\cG(x,z,r)$ is independent of $z$, $r$ for $z\in[0,\epsilon]$, that is $\cG(x,z,r)=\cG(x,0,0)$, meaning $(\cG|_r)^{\bs\phi_r}(x,z)=\cG(x,0,0)\otimes dt$.

  For $z\in[1-\epsilon,1]$ we have $H^r_z\equiv 0$ and $\bs\phi_r(z)\equiv \phi_r$, therefore
  $$(\cG|_r)^{\bs\phi_r}(x,z) = \cG(x,z,r)\circ\phi_r\otimes dt\,.$$
  Again, since $\cG$ is $\bs\phi$-parametrized, for $z\in[1-\epsilon,1]$ we have $\cG(x,z,r) = \cG(x,1,0)\circ\phi_r^{-1}$, and $\cG(x,z,0)$ is also independent of $z$. It follows that $(\cG|_r)^{\bs\phi_r}(x,z) = \cG(x,z,r)\circ\phi_r\otimes dt = \cG(x,1,0)\otimes dt$ is independent of $z,r$ there.
\end{proof}

Next, we will introduce a certain smooth map $\Xi_n \fc [0,1]^n\to\Slit^n$, whose existence and properties are the subject of Corollary \ref{cor:existence_properties_Xi_n}. Its goal is as follows. Suppose we are given an $n$-cube $\cK$ of Hamiltonian $1$-forms which is independent of $z,r$ for $z\in[0,\epsilon]\cup[1-\epsilon,1]$, where, like above, we denote points in $[0,1]^n$ by $(x,z,r)$ with $z,r\in[0,1]$. Then Corollary \ref{cor:descent_forms_Xi} states that there exists a unique $n$-slit $\wt\cK$ of Hamiltonian $1$-forms satisfying $\wt\cK\circ\Xi_n=\cK$.

Let us describe this in detail. First, we treat the case $n=2$. In the formulation we write $\Phi_t$ for the time-$t$ map of the negative gradient flow of $\rho_1$. Since $\rho_n=\rho_1\oplus\dots\oplus\rho_1$, it follows that the time-$t$ map of the negative gradient flow of $\rho_n$ is $\Phi_t^{\times n}=\Phi_t\times\dots\times\Phi_t$.
\begin{prop}\label{prop:existence_properties_Xi_2}
  There exists a unique continuous map $\Xi_2 \fc [0,1]^2 \to \Slit^2$ satisfying
  \begin{itemize}
    \item $\Xi_2(\tfrac12,r) = (\tfrac{r+1}3,\tfrac{2-r}3)$ for $r\in[0,1]$;
    \item $\Xi_2$ intertwines the negative gradient flows of the functions $\rho_1\oplus0$ on $[0,1]^2$ and $\rho_2$ on $\Slit^2$, that is $\Xi_2\circ(\Phi_t\times\id) = \Phi_t^{\times 2}\circ\Xi_2$.
  \end{itemize}
  Moreover, $\Xi_2$ is smooth and in addition satisfies:
  \begin{itemize}
    \item $\Xi_2(0,r) = (0,0)$, $\Xi_2(1,r) = (1,1)$ for $r\in[0,1]$;
    \item $\Xi_2$ maps $(0,1)\times[0,1]$ diffeomorphically onto $\Slit^2\setminus\{(0,0),(1,1)\}$;
    \item for each $r\in[0,1]$, $\Xi_2|_{[0,1]\times\{r\}}$ is a diffeomorphism from $[0,1]\cong[0,1]\times\{r\}$ onto the image $\Xi_2([0,1]\times\{r\})$, which is the closure of the negative gradient trajectory of $\rho_2$ passing through $(\tfrac{r+1}3,\tfrac{2-r}3)$.
  \end{itemize}
\end{prop}
\begin{figure}
	\centering
	\includegraphics[width=0.4\linewidth]{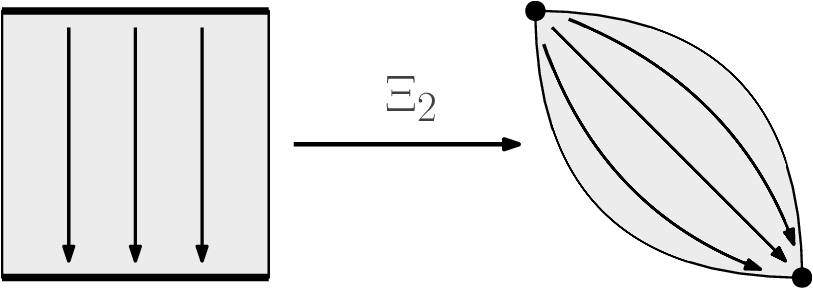}
	\caption{The map $\Xi_2$. On the left, the coordinate $z$ goes down, the coordinate $r$ goes right. On the right, the origin $0\in[0,1]^2$ is in the top left}
	\label{fig:Xi2}
\end{figure}
\noindent See Figure \ref{fig:Xi2} for an illustration. Proposition \ref{prop:existence_properties_Xi_2} is proved below, after the proof of Lemma \ref{lem:unique_descent_forms_by_Xi}. In what follows, we will need a higher-dimensional version of $\Xi_2$, namely we introduce
\begin{notation}
  For $n\geq 2$ we define $\Xi_n \fc [0,1]^n \to \Slit^n = [0,1]^{n-2}\times\Slit^2$ by
  $$\Xi_n:=\id_{[0,1]^{n-2}}\times\Xi_2\,.$$
\end{notation}
\noindent The following is an immediate consequence of Proposition \ref{prop:existence_properties_Xi_2}.
\begin{coroll}\label{cor:existence_properties_Xi_n}
  The map $\Xi_n \fc [0,1]^n\to\Slit^n$ is smooth and satisfies:
  \begin{itemize}
    \item $\Xi_n$ maps $[0,1]^{n-2}\times(0,1)\times[0,1]$ onto $[0,1]^{n-2}\times(\Slit^2\setminus\{(0,0),(1,1)\})$ diffeomorphically;
    \item $\Xi_n(x,0,r)=(x,0,0)$, $\Xi_n(x,1,r)=(x,1,1)$ for all $x\in[0,1]^{n-2}$, $r\in [0,1]$;
    \item $\Xi_n$ intertwines the negative gradient flows of $\rho_{n-1}\oplus 0$ on $[0,1]^n$ and $\rho_n$ on $\Slit^n$, that is $\Xi_n\circ (\Phi_t^{\times(n-1)}\times\id_{[0,1]})=\Phi_t^{\times n}\circ\Xi_n$;
    \item for each $r\in[0,1]$, $\Xi_n|_{[0,1]^{n-1}\times\{r\}}$ is a diffeomorphism from $[0,1]^{n-1}\cong[0,1]^{n-1}\times\{r\}$ onto its image $\Xi_n([0,1]^{n-1}\times\{r\})$.      \qed
  \end{itemize}
\end{coroll}
\noindent Let us see what Proposition \ref{prop:existence_properties_Xi_2} implies in terms of going from cubes to slits. In the following lemma, the manifold $P$ plays the role of parameters appearing in a function space. See Corollary \ref{cor:descent_forms_Xi} and its proof for the relevant application.
\begin{lemma}\label{lem:unique_descent_forms_by_Xi}
  Let $P,Q$ be smooth manifolds, and let $\cK \fc [0,1]^n\times P \to Q$ be a smooth map such that there exists $\epsilon \in (0,\frac12)$ with the property that for all $x\in[0,1]^{n-2}$, $r\in [0,1]$, and $p \in P$ we have: if $z \in [0,\epsilon]$, then $\cK(x,z,r;p)=\cK(x,0,0;p)$ and if $z \in [1-\epsilon,1]$ then $\cK(x,z,r;p)=\cK(x,1,0;p)$. Then there is a unique smooth $\wt\cK \fc \Slit^n\times P \to Q$ satisfying $\wt\cK \circ (\Xi_n\times \id_P)=\cK$.
\end{lemma}
\noindent The proof appears below. For now, let us derive from this the desired conclusion regarding the descent of cubes of Hamiltonian $1$-forms to slits.
\begin{coroll}\label{cor:descent_forms_Xi}
  Let $\cK \fc [0,1]^n \to \Omega^1(\R\times S^1,C^\infty(M))$ be an $n$-cube of Hamiltonian $1$-forms. Let $(z,r)$ denote the last two coordinates of $[0,1]^n=[0,1]^{n-2}\times[0,1]\times[0,1]$ and assume there is $\epsilon\in(0,\frac12)$ such that for $z\in[0,\epsilon]\cup[1-\epsilon,1]$, $\cK$ is independent of both $z$ and $r$. Then there is a unique smooth $n$-slit $\wt\cK$ of Hamiltonian $1$-forms satisfying $\wt\cK\circ\Xi_n=\cK$. Moreover, if $\cK$ takes values in $s$-independent forms, then so does $\wt\cK$.
\end{coroll}
\begin{proof}
  A Hamiltonian $1$-form on $\R\times S^1$ given by $F\otimes ds + G\otimes dt$ can be identified with the pair of functions $(F,G)$, or equivalently, with a smooth map $(F,G)\fc \R\times S^1\times M \to \R^2$. It follows that a smooth $n$-cube $\cK$ of Hamiltonian $1$-forms is identified with a smooth map $\cK' \fc [0,1]^n\times\R\times S^1\times M \to \R^2$. Moreover, by assumption, $\cK'$ is independent of $z,r$ for $z\in [0,\epsilon]\cup[1-\epsilon,1]$. Thus we are in the situation of Lemma \ref{lem:unique_descent_forms_by_Xi} with $P=\R\times S^1\times M$ and $Q=\R^2$, which yields a unique smooth map $\wt\cK' \fc \Slit^n\times\R\times S^1 \times M \to \R^2$ satisfying $\wt\cK'\circ(\Xi_n\times\id)=\cK'$. Unraveling the identifications, we see that $\wt\cK'$ is equivalent to a smooth $n$-slit of Hamiltonian $1$-forms $\wt\cK$ which satisfies $\wt\cK\circ\Xi_n=\cK$. The last assertion is obvious.
\end{proof}
\begin{proof}[Proof of Lemma \ref{lem:unique_descent_forms_by_Xi}]
  The uniqueness is obvious since $\Xi_n$ is surjective. Let us show that $\wt\cK$ is well-defined. Indeed, $\Xi_n$ is a bijection between $[0,1]^{n-2}\times(0,1)\times[0,1]$ and $[0,1]^{n-2}\times(\Slit^2\setminus\{(0,0),(1,1)\})$, meaning on the set $[0,1]^{n-2}\times(\Slit^2\setminus\{(0,0),(1,1)\})\times P$, $\wt\cK$ is given by the explicit formula
  $$\wt\cK|_{[0,1]^{n-2}\times(\Slit^2\setminus\{(0,0),(1,1)\})\times P} = \cK\circ\big(\Xi_n|_{[0,1]^{n-2}\times(0,1)\times[0,1]}\times \id_P\big)^{-1}\,.$$
  Now, for $x\in[0,1]^{n-2}$ and any $r\in[0,1]$ we have $\Xi_n(x,0,r)=(x,0,0)$, meaning for all $p\in P$ we must have
  $$\cK(x,0,r;p) = \wt\cK((\Xi_n\times\id_P)(x,0,r;p)) = \wt\cK(x,0,0;p)\,.$$
  However, $\cK(x,0,r;p)=\cK(x,0,0;p)$ for all $x,p$, meaning the right-hand side is indeed well-defined. Similarly,
  $$\wt\cK(x,1,1;p) = \cK(x,1,r;p)\quad\text{is well-defined.}$$
  It remains to prove that $\wt\cK$ is smooth. Define
  $$U_0=\Xi_2\big([0,\epsilon)\times[0,1]\big)\,,\quad U_1 = \Xi_2\big((1-\epsilon,1]\times[0,1]\big)\,.$$
  Since $\Xi_2$ is a continuous surjection from a compact space to a Hausdorff space, it is a quotient map. Moreover,
  $$\Xi_2^{-1}(U_0) = [0,\epsilon)\times[0,1]\,,\qquad\Xi_2^{-1}(U_1) = (1-\epsilon,1]\times[0,1]\,,$$
  both of which are open in $[0,1]^2$, implying that $U_0,U_1$ are open in $\Slit^2$. Thus we have an open cover $\Slit^2=\Slit^2\setminus\{(0,0),(1,1)\}\cup U_0\cup U_1$, which yields the following open cover of $\Slit^n\times P$:
  \begin{multline*}
    \Slit^n\times P = \big([0,1]^{n-2}\times (\Slit^2\setminus\{(0,0),(1,1)\})\times P\big)\\\cup \big([0,1]^{n-2}\times U_0\times P\big)\cup \big([0,1]^{n-2}\times U_1\times P\big)\,.
  \end{multline*}
  It suffices to prove that the restriction of $\wt\cK$ to each one of the sets in this open cover is smooth. On the first piece we have, thanks to the above:
  $$\wt\cK|_{[0,1]^{n-2}\times (\Slit^2\setminus\{(0,0),(1,1)\})\times P} = \cK\circ \big(\Xi_n|_{[0,1]^{n-2}\times(0,1)\times[0,1]}\times\id_P\big)^{-1}\,.$$
  This is smooth because $\cK$ is, and because $\Xi_n$ is a diffeomorphism $[0,1]^{n-2}\times(0,1)\times[0,1] \to [0,1]^{n-2}\times (\Slit^2\setminus\{(0,0),(1,1)\})$. To show smoothness on $[0,1]^{n-2}\times U_0\times P$, let us introduce the projection $\pr \fc \Slit^n=[0,1]^{n-2}\times\Slit^2\to[0,1]^{n-2}$ and the map $\cK_0\fc[0,1]^{n-2}\times P \to Q$ defined by $\cK_0(x;p) = \cK(x,0,0;p)$.
  For $u \in U_0$ let $(z,r) \in [0,1]^2$ be such that $\Xi_2(z,r) = u$. It follows that $z\in[0,\epsilon)$. From the definitions we then have for all $x\in[0,1]^{n-2}$, $p\in P$:
  \begin{multline*}
    \cK_0(x;p) = \cK(x,0,0;p) = \cK(x,z,r;p)=\wt\cK(\Xi_n(x,z,r);p)\\=\wt\cK(x,\Xi_2(z,r);p) = \wt\cK(x,u;p)\,.
  \end{multline*}
  It follows that
  $$\wt\cK|_{[0,1]^{n-2}\times U_0\times P} = \cK_0\circ(\pr|_{[0,1]^{n-2}\times U_0}\times\id_P)\,,$$
  which is manifestly smooth. A similar argument shows that $\wt\cK$ is smooth on $[0,1]^{n-2}\times U_1\times P$. The proof is complete.
\end{proof}

\begin{proof}[Proof of Proposition \ref{prop:existence_properties_Xi_2}]
  The uniqueness can be seen as follows. Thanks to the existence and uniqueness of solutions of ODEs, and the asymptotics of the negative gradient flow of $\rho_1$, for any $z\in(0,1)$ there is a unique $t\in\R$ such that $z=\Phi_t(\frac12)$. It follows that
  \begin{multline*}
    \Xi_2(z,r)=\Xi_2(\Phi_t(\tfrac12),r)=(\Xi_2\circ(\Phi_t\times\id))(\tfrac12,r)=(\Phi_t^{\times2}\circ\Xi_2)(\tfrac12,r)\\=\Phi_t^{\times2}(\Xi_2(\tfrac12,r))= \Phi_t^{\times2}(\tfrac{r+1}3,\tfrac{2-r}3)\,,
  \end{multline*}
  meaning $\Xi_2$ is uniquely determined on $(0,1)\times[0,1]$. The continuity of $\Xi_2$ then implies that $\Xi_2(0,r)=(0,0)$ and $\Xi_2(1,r)=(1,1)$ for all $r \in [0,1]$. Indeed, if $(0,1)\times[0,1]\ni(a,b)\to(0,r)$, then $b\to r$ and the unique $t(a)\in\R$ satisfying $\Phi_{t(a)}(\frac12)=a$ satisfies $t(a)\to-\infty$. It follows that
  $$\Xi_2(a,b)=\Xi_2(\Phi_{t(a)}(\tfrac12),b)=\Phi_{t(a)}^{\times2}(\tfrac{b+1}3,\tfrac{2-b}3) \to (0,0)\,,$$
  again thanks to the asymptotics of $\Phi_t$ as $t\to-\infty$. It follows that $\lim_{(a,b)\to(0,r)}\Xi_2(a,b)$ indeed exists and equals $(0,0)$. A similar argument shows that $\Xi_2(1,r)=(1,1)$. This proves both that $\Xi_2$ is uniquely determined and that it is continuous.

  As for existence, it can be inferred from properties of solutions of ODEs, however there is an explicit formula for $\Xi_2$. First, let us consider the map $\Psi \fc [0,1]\times(0,1) \to [0,1]^2$ given by
  $$\Psi(z,r) = \frac1\pi\left(\arccos\frac{\cos\pi z + \cos \pi r}{1+\cos\pi z\cos \pi r}\,,\arccos\frac{\cos\pi z - \cos \pi r}{1-\cos\pi z\cos \pi r}\right)\,.$$
  An elementary, but not entirely trivial check shows that $\Psi$ is the unique continuous map satisfying $\Psi(\tfrac12,r)=(r,1-r)$ for $r\in(0,1)$ and $\Psi\circ(\Phi_t\times\id)=\Phi_t^{\times2}\circ\Psi$, where the uniqueness is proved similarly to that of $\Xi_2$. This calculation is based on the fact that for $\rho_1(x)=\cos\pi x$ the negative gradient line $\gamma \fc \R\to (0,1)$ such that $\gamma(0)=\tfrac12$ is given by
  $$\gamma(t) = \frac1\pi\arccos(-\tanh \pi^2 t)\,,$$
  as well as on the addition formula for the hyperbolic tangent. It then follows that $\Xi_2$ must be given by
  $$\Xi_2(z,r) = \Psi(z,\tfrac{r+1}3)\,.$$
  The smoothness of $\Psi$ is proved by inspecting its first partial derivatives:
  \begin{equation}\label{eq:derivative_z_Psi}
    \partial_z\Psi = \left(\frac{\sin\pi r}{1+\cos\pi z\cos\pi r}\,,\frac{\sin\pi r}{1-\cos\pi z\cos\pi r}\right)\,,
  \end{equation}
  $$\partial_r\Psi = \left(\frac{\sin\pi z}{1+\cos\pi z\cos\pi r}\,,\frac{-\sin\pi z}{1-\cos\pi z\cos\pi r}\right)\,,$$
  which are obviously smooth. The smoothness of $\Xi_2$ follows from that of $\Psi$.

  Of the three claimed properties of $\Xi_2$ we have already established the first one. The assertion that $\Xi_2$ maps $(0,1)\times[0,1]$ diffeomorphically onto $\Slit^2\setminus\{(0,0),(1,1)\}$ follows from the fact that $\Psi$ maps $(0,1)^2$ diffeomorphically onto itself. This can probably be proved based on properties of solutions of ODEs, however there is an explicit smooth inverse for $\Psi|_{(0,1)^2}$:
  $$\Psi^{-1}(s,t) = \frac1\pi\left(\arccos\frac{\cos\pi\frac{s+t}2}{\cos\pi\frac{s-t}2}\,,\arccos\frac{\sin\pi\frac{t-s}2}{\sin\pi\frac{s+t}2}\right)\,.$$

  To prove the third claimed property, namely that $\Xi_2|_{[0,1]\times\{r\}}$ is a diffeomorphism onto the image, it suffices to show that $\Psi$ has the same property, that is for each $r\in(0,1)$, $\Psi|_{[0,1]\times\{r\}}$ is a diffeomorphism onto $\Psi([0,1]\times\{r\})$. Since $\Psi|_{[0,1]\times\{r\}} \fc [0,1]\times\{r\} \to [0,1]^2$ is a function of one variable, in order to prove it is a diffeomorphism onto the image, it suffices to show that it is a homeomorphism onto the image and that it is an immersion. To prove it is an immersion, consider its derivative, given by \eqref{eq:derivative_z_Psi}; it is clearly nowhere zero if $r \in (0,1)$. To show it is a homeomorphism, note that it is a injective continuous map from a compact space onto a Hausdorff space. To prove injectivity, recall that $\Psi|_{(0,1)\times\{r\}}$, as a restriction of the diffeomorphism $\Psi|_{(0,1)^2}$, is itself injective. Finally, the endpoints $(0,r)$, $(1,r)$ are mapped to $(0,0)$, $(1,1)$, respectively, which are two distinct points, and which do not lie on $\Psi((0,1)\times\{r\})\subset(0,1)^2$. This concludes the proof of the proposition.
\end{proof}

\subsubsection{Completing the proof of Theorem \ref{thm: invariance} (invariance)}\label{sss:proof_invariance}

\begin{figure}
	\centering
	\includegraphics[width=0.3\linewidth]{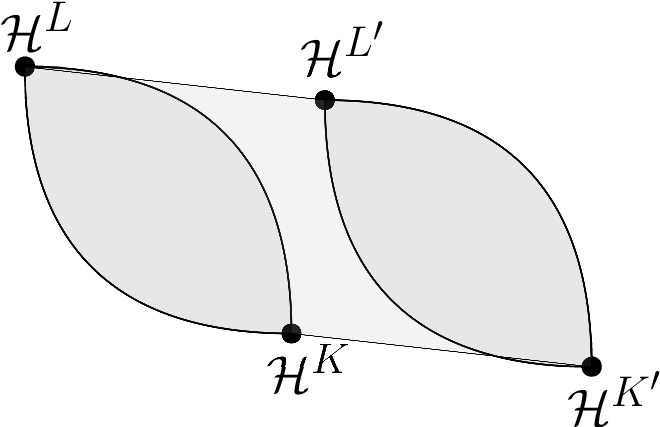}
	\caption{The slit-like $4$-ray of Hamiltonian forms $\wt\cK_\bullet$}
    \label{fig:Ham4ray}
\end{figure}

Let us abbreviate $\Xi=\Xi_4$ until the end of the proof. Fix a $\bs\phi$-parametrized $4$-ray $\cG_\bullet = (\cG_i)_{i\in \N}$ as in Theorem \ref{thm:exists_parametrized_ray}. Just like in the beginning of Section \ref{sss:cubes_to_slits}, pick a family of isotopies $(\bs\phi_\sigma)_\sigma$ and construct the corresponding $4$-cube of $1$-forms $\cK_i = \cG_i^{(\bs\phi_\sigma)_\sigma}$. Form the corresponding $4$-ray $\cK_\bullet=(\cK_i)_{i\in\N}$. The cubes $\cK_i$ are independent of $z,r$ for $z$ near $0,1$, thanks to Lemma \ref{lem:cube_of_forms_conjugated_paramd_indep_z_r}. Corollary \ref{cor:descent_forms_Xi} then implies that for each $i$ there is a unique smooth $4$-slit of Hamiltonian $1$-forms $\wt\cK_i$ with $\wt\cK_i\circ\Xi=\cK_i$. Thus $\wt\cK_\bullet=(\wt\cK_i)_{i\in\N}$ is a smooth slit-like $4$-ray of Hamiltonian $1$-forms. See Figure \ref{fig:Ham4ray} for a schematic depiction.
Here the acceleration $1$-rays $\cH^{L,L',K,K'}$ for $L,L',K,K'$, which are the corner $1$-rays of $\cG|_0$, appear as the vertices, and we think of the infinite ray coordinate as being hidden in each vertex.

Let $\wt\cJ_\bullet$ be any slit-like $4$-ray of almost complex structures, and let $\cJ_\bullet=(\wt\cJ_i\circ\Xi)_{i\in\N}$. Applying the Floer functor to the pair $(\wt\cK_\bullet,\wt\cJ_\bullet)$ yields a slit-like algebraic $4$-ray of the form
$$\xymatrix{CF^*(\cH^L)\ar@{=}[rr]\ar[dd]\ar[rd] && CF^*(\cH^L)\ar[dd]|\hole\ar[rd] \\
            & CF^*(\cH^{L'})\ar@{=}[rr]\ar[dd] && CF^*(\cH^{L'})\ar[dd] \\
            CF^*(\cH^K)\ar@{=}[rr]|\hole\ar[rd] && CF^*(\cH^K)\ar[rd] \\
            & CF^*(\cH^{K'})\ar@{=}[rr] && CF^*(\cH^{K'}) \\}$$
where we use Notation \ref{not:Floer_1_ray} for simplicity. Since the restrictions of $\wt\cK_\bullet$ to the four corner rays of the slit-like $4$-ray are the acceleration $1$-rays $\cH^{L,L',K,K'}$, it follows that upon taking the telescope of the above slit-like algebraic $4$-ray and completing, we obtain a slit-like algebraic $3$-cube of the form
$$\xymatrix{SC^*(L)\ar@{=}[rr]\ar[dd]\ar[rd] && SC^*(L)\ar[dd]|\hole\ar[rd] \\
            & SC^*(L')\ar@{=}[rr]\ar[dd] && SC^*(L')\ar[dd] \\
            SC^*(K)\ar@{=}[rr]|\hole\ar[rd] && SC^*(K)\ar[rd] \\
            & SC^*(K')\ar@{=}[rr] && SC^*(K') \\}$$
At this point we note that the $\wt\cK_i$, as functions defined on $\Slit^4=[0,1]^2\times\Slit^2$, are independent of the last two coordinates in a neighborhood of $[0,1]^2\times\{(0,0),(1,1)\}$, and that the restrictions of $\wt\cK_i$ to $[0,1]^2\times\{(0,0)\}$ and to $[0,1]^2\times\{(1,1)\}$ coincide with acceleration $2$-rays for the pairs $(L,L')$ and $(K,K')$, respectively. It follows that the oblique arrows in the last diagram are usual restriction maps (see Remark \ref{rem:Ham_forms_acc_data_induce_restrictions} for a more detailed argument), thus its cocone in the oblique direction is a slit-like $2$-cube
\begin{equation}\label{eq:diagram_res_01}
  \xymatrix{SC^*(L,L') \ar[r]^{\id} \ar[d]_{\res_0} & SC^*(L,L') \ar[d]^{\res_1} \\ SC^*(K,K') \ar[r]^{\id} & SC^*(K,K')}
\end{equation}
where $\res_0,\res_1$ are certain maps, which are identified in the following lemma.
\begin{lemma}\label{lem:computing_res_01}
  For $j=0,1$ we have $\res_j = \res(\cK_\bullet|_j,\cJ_\bullet|_j)$ (see Notation \ref{not:restriction_Ham_1_forms}).
\end{lemma}
\begin{rem}
  \emph{A priori}, despite the notation, we do not know if $\res_0,\res_1$ are usual restrictions, since we have used Hamiltonian forms in their definition.
\end{rem}
\noindent The lemma is proved below. Accepting it for a moment, let us see how this allows us to complete the proof of the theorem. Now, since $\cK_\bullet|_0=\cG_\bullet|_0$ is a usual monotone Hamiltonian $3$-ray, in fact an acceleration $3$-ray for the inclusion of pairs $(K,K')\subset(L,L')$, the map $\res(\cK_\bullet|_0,\cJ_\bullet|_0)$ is a usual restriction map. Let us consider the other map $\res(\cK_\bullet|_1,\cJ_\bullet|_1)$. From the definitions it follows that $\cK_\bullet|_1 = \big(\cG_\bullet|_1\big)^{\bs\phi_1}$. Setting $\psi_\sigma:=\bs\phi_1(\sigma)^{-1}$ and $\bs\psi:=(\psi_\sigma)_\sigma$, we see that $\cJ_\bullet|_1=\big((\cJ_\bullet|_1)^{\bs\psi}\big)^{\bs\phi_1}$, that is
$$\res(\cK_\bullet|_1,\cJ_\bullet|_1) = \res\Big((\cG_\bullet|_1)^{\bs\phi_1},\big((\cJ_\bullet|_1)^{\bs\psi}\big)^{\bs\phi_1}\Big)\,,$$
meaning the pair $(\cK_\bullet|_1,\cJ_\bullet|_1)$ is obtained from the pair $\big(\cG_\bullet|_1,(\cJ_\bullet|_1)^{\bs\psi}\big)$ by conjugation with respect to $\bs\phi_1$. Note that $\cG_\bullet|_1$ is an acceleration $3$-ray for the inclusion of pairs $(\phi(K),\phi(K'))\subset (L,L')$ (see Remark \ref{rem:phi_paramd_ray_slices_acc_rays}), therefore
$$\res\big(\cG_\bullet|_1,(\cJ_\bullet|_1)^{\bs\psi}\big) \fc SC^*(L,L') \to SC^*(\phi(K),\phi(K'))\quad\text{is a usual restriction}.$$
It follows from Proposition \ref{prop:conjugated_cube_isomorphic} that the maps $\res(\cK_\bullet|_1,\cJ_\bullet|_1)$, $\res\big(\cG_\bullet|_1,(\cJ_\bullet|_1)^{\bs\psi}\big)$ are isomorphic, as $1$-cubes, through the diagram
\begin{equation}\label{eq:right_half_diagram}
  \xymatrix{SC^*(L,L') \ar[d]_{\res(\cK_\bullet|_1,\cJ_\bullet|_1)} \ar[r]^{\id} & SC^*(L,L') \ar[d]^{\res\big(\cG_\bullet|_1,(\cJ_\bullet|_1)^{\bs\psi}\big)} \\ SC^*(K,K') \ar[r]^-{\phi_*} & SC^*(\phi(K),\phi(K'))}
\end{equation}
We have already remarked that the right vertical arrow is the usual restriction map. Combining the diagrams \eqref{eq:diagram_res_01}, \eqref{eq:right_half_diagram} and passing to cohomology, we obtain the following commutative diagram, where $\res$ denotes the usual restriction maps:
$$\xymatrix{SH^*(L,L')\ar[r]^{\id} \ar[d]_{\res} & SH^*(L,L')\ar[r]^{\id} \ar[d]_{\res(\cK_\bullet|_1,\cJ_\bullet|_1)} & SH^*(L,L')\ar[d]^{\res} \\ SH^*(K,K')\ar[r]^{\id} & SH^*(K,K')\ar[r]^-{\phi_*} & SH^*(\phi(K),\phi(K'))}$$
Going around the diagram yields the desired identity \eqref{eq:identity_invariance}. This completes the proof.

\begin{proof}[Proof of Lemma \ref{lem:computing_res_01}]
  Consider the bounding ``walls'' of $\Slit^4$, $W_j:=\Xi([0,1]^3\times\{j\})$, $j=0,1$. The maps $\res_j$, $j=0,1$, come, via the Floer functor, from the data $(\wt\cK_\bullet,\wt\cJ_\bullet)$ restricted to these walls, that is the map $\res_j$ comes from $(\wt\cK_\bullet|_{W_j},\wt\cJ_\bullet|_{W_j})$, where $\wt\cK_\bullet|_{W_j} = (\wt\cK_i|_{W_j})_{i\in\N}$, $\wt\cJ_\bullet|_{W_j} = (\wt\cJ_i|_{W_j})_{i\in\N}$. The lemma states that $\res_j$ coincides with the map $\res(\cK_\bullet|_j,\cJ_\bullet|_j)$. The maps $\res_j$ and $\res(\cK_\bullet|_j,\cJ_\bullet|_j)$ are induced by virtual counts of compactified moduli spaces consisting of pairs $(\gamma,u)$ where $\gamma$ is a negative gradient trajectory and $u$ is a Floer solution, possibly broken and with holomorphic spheres attached. The maps $\res_j$ are induced by such counts where $\gamma$ is a negative gradient trajectory of $\rho_4$ contained in $W_j$ while $u$ solves
  \begin{equation}\label{eq:tilde_K_gamma}
    \big(du-X_{\wt\cK_i\circ\gamma}(u)\big)^{1,0}_{\wt\cJ_i\circ\gamma}=0\,,
  \end{equation}
  where $i$ is the number of the constituent cube under consideration. On the other hand, the map $\res(\cK_\bullet|_j,\cJ_\bullet|_j)$ is induced by such counts with $\gamma$ a negative gradient trajectory of $\rho_3\oplus 0$ on $[0,1]^3\times\{j\}$, while $u$ solves
  \begin{equation}\label{eq:K_delta}
    \big(du-X_{\cK_i\circ\gamma}(u)\big)^{1,0}_{\cJ_i\circ\gamma}=0\,.
  \end{equation}
  The intuitive idea behind the equality $\res_j=\res(\cK_\bullet|_j,\cJ_\bullet|_j)$ is that $\Xi$ induces diffeomorphisms $[0,1]^3\times\{j\} \to W_j$, which intertwine the gradient flows of $\rho_3\oplus0$, $\rho_4$, respectively, and thus identifies the spaces of negative gradient trajectories on both sides (see Corollary \ref{cor:existence_properties_Xi_n}). The identification is given by $\gamma \mapsto \Xi\circ\gamma$ for a negative gradient trajectory $\gamma$ in $[0,1]^3\times\{j\}$. Also, $\Xi$ identifies the data on both sides: $\wt\cK_i|_{W_j}\circ\Xi=\cK_i|_j$, $\wt\cJ_i|_{W_j}\circ\Xi=\cJ_i|_j$.

  In more detail, let us show that $(\gamma,u)\mapsto(\Xi\circ\gamma,u)$ maps solutions of \eqref{eq:K_delta} to solutions of \eqref{eq:tilde_K_gamma}. Indeed, if $(\gamma,u)$ solves \eqref{eq:K_delta}, then
  $$\big(du-X_{\wt\cK_i\circ(\Xi\circ\gamma)}(u)\big)^{1,0}_{\wt\cJ_i\circ(\Xi\circ\gamma)}= \big(du-X_{(\wt\cK_i\circ\Xi)\circ\gamma}(u)\big)^{1,0}_{(\wt\cJ_i\circ\Xi)\circ\gamma}=
  \big(du-X_{\cK_i\circ\gamma}(u)\big)^{1,0}_{\cJ_i\circ\gamma}=0\,,$$
  as claimed. Thus the map $(\gamma,u)\mapsto(\Xi\circ\gamma,u)$ induces a homeomorphism between the compactified moduli spaces on both sides, and this includes broken trajectories and holomorphic spheres, since the PDE solutions are preserved by the map. This map also clearly preserves the linearized operators of the above parametrized equations, and thus can be used to identify stabilizations of these linearizations. The same holds for the sphere components, which are again untouched by the map. It follows that the map identifies the implicit atlases on both sides, and can also be used to identify auxiliary data, which results in an equality of the virtual counts. What we get from this is an identification of the algebraic $3$-rays coming from $(\wt\cK_\bullet|_{W_j},\wt\cJ_\bullet|_{W_j})$, on the one hand, and from $\res(\cK_\bullet|_j,\cJ_\bullet|_j)$, on the other hand. The claim that $\res_j=\res(\cK_\bullet|_j,\cJ_\bullet|_j)$ follows.
\end{proof}

\subsubsection{Proof of Theorem \ref{thm:exists_parametrized_ray}}\label{sss:proof_exists_parametrized_ray}

Recall that the theorem states that if $(K,K')\subset(L,L')$ are compact pairs and $\bs\phi=(\phi_t)_{t\in[0,1]}$ is a Hamiltonian isotopy such that $(\phi_t(K),\phi_t(K'))\subset(L,L')$ for all $t$, then there exists a monotone $\bs\phi$-parametrized Hamiltonian $4$-ray $\cG_\bullet=(\cG_i)_{i\in\N}$ such that $\cG_\bullet|_0$ is an acceleration $3$-ray for the inclusion of pairs $(K,K')\subset (L,L')$.
\begin{notation}
  For Hamiltonian $n$-cubes $\cG,\cG'$ we will denote $\cG\leq \cG'$ to mean that $\cG(x)\leq\cG'(x)$ as functions, for each $x$. If $\cG_\bullet=(\cG_i)_{i\in\N}$, $\cG'_\bullet=(\cG_i')_{i\in\N}$ are Hamiltonian $n$-rays, we write $\cG_\bullet\leq \cG_\bullet'$ to mean $\cG_i\leq\cG_i'$ for all $i$.
\end{notation}

The proof proceeds in stages:
\begin{itemize}
  \item First, we recursively construct acceleration $1$-rays $\cH_\bullet^A$ for $A\in\{L,L',K,K'\}$ such that $\cH_\bullet^L\leq \cH_\bullet^{L'}$, $\cH_\bullet^K \leq \cH_\bullet^{K'}$ and $\cH_\bullet^L\leq\cH_\bullet^K\circ\phi_t^{-1}$, $\cH_\bullet^{L'}\leq \cH_\bullet^{K'}\circ\phi_t^{-1}$ for all $t$.
  \item We then use interpolation to produce monotone $2$-rays $\cH^0,\cH^1$ such that the corner $1$-rays of $\cH^0$ are $\cH^L$, $\cH^{L'}$, while those of $\cH^1$ are $\cH^K$, $\cH^{K'}$, such that monotonicity with respect to the isotopy $\bs\phi$ is preserved in the sense that $\cH^0 \leq \cH^1\circ\phi_t^{-1}$ for all $t$;
  \item Finally, we construct the required monotone $\bs\phi$-parametrized $4$-ray using another interpolation step.
\end{itemize}
\noindent Let us describe all of this in detail.

\begin{lemma}\label{lemma: uniform_majorant}
    Let $A\subset M$ be compact, and let $\mathcal P\subset C_{A\subset M}^\infty$ be a nonempty family which is compact in the $C^0$-topology. Then there exists a nondegenerate Hamiltonian $G\in C^\infty_{A\subset M}$ such that for every $F\in\mathcal P$ we have $F<G$.
\end{lemma}
\begin{proof}
Let us define a function $\max \cP$ on $S^1\times M$ as follows:
$$\max\cP(t,m)=\sup_{F \in \cP}F(t,m)\,.$$
We claim that $\max\cP$ is a well-defined, real-valued, continuous function, which is everywhere negative on $S^1\times A$, and which satisfies $F\leq \max\cP$ for all $F\in \cP$. First, consider the evaluation $\ev_{(t,m)} \fc \cP \to \R$, $\ev_{(t,m)}(F)=F(t,m)$. Since $\ev_{(t,m)}$ is continuous with respect to the $C^0$-topology and $\cP$ is compact and nonempty, the image $\ev_{(t,m)}(\cP)$ is a compact nonempty subset of $\R$ and thus has a  maximum, which equals $\max\cP(t,m)$ by definition. Thus, $\max\cP$ is a well-defined real-valued function. The continuity of $\max\cP$ is a consequence of the uniform equicontinuity of $\cP$ which holds by the Arzel\`a--Ascoli theorem. Finally, if $(t,m)\in S^1\times A$, then $\max\cP(t,m) = \max \ev_{(t,m)}(\cP)$, where $\ev_{(t,m)}(\cP)\subset(-\infty,0)$ by assumption. Thus the maximum is negative. The last property, namely that $F\leq\max\cP$, holds trivially.

Since $A$ is compact, $S^1\times A$ is compact, and since $\max\cP$ is continuous and negative on $S^1\times A$, we deduce that there exists $\delta>0$ such that $\max\cP|_{S^1\times A}\leq -\delta$. Choose a smooth function $g_0$ satisfying $\|g_0-\max\cP\|_{C^0}<\delta/4$ and put $g_1=g_0+\delta/2$. Then $g_1-\max\cP>\delta/4$ and $g_1|_{S^1\times A} <-\delta/4$. Perturbing $g_1$, we can find a nondegenerate Hamiltonian $G$ satisfying $\|G-g_1\|_{C^0}<\delta/8$. Then $G-\max\cP>\delta/8$ and $G|_{S^1\times A}<-\delta/8$. Consequently, $G\in C^\infty_{A\subset M}$ and $F<G$ for every $F\in\mathcal P$.
\end{proof}

\begin{prop}\label{prop:adapted_acceleration_data}
    Let $(K,K'),(L,L')\in\cC\cP(M)$ and let $\phi=(\phi_t)_{t\in[0,1]}$ be a Hamiltonian isotopy such that for every $t\in [0,1]$ we have $(\phi_t(K),\phi_t(K'))\subset (L,L')$.
    Then there exist acceleration $1$-rays $\cH_\bullet^A=(\cH_i^A)_{i\in\N}$ for $A \in\{L,L',K,K'\}$, such that
    \[\cH_\bullet^L\leq\cH_\bullet^{L'},\;\cH_\bullet^{K}\leq \cH_\bullet^{K'},\;\text{and for all }t\in[0,1]:\;\cH_\bullet^L \leq \cH_\bullet^K\circ\phi_t^{-1},\;\cH_\bullet^{L'} \leq \cH_\bullet^{K'}\circ\phi_t^{-1}.\]
\end{prop}

\begin{proof}
    For each $A\in\{K,K',L,L'\}$, fix an arbitrary acceleration datum $\cR^A=(R_i^A)_{i\in\N}$ for $A$. First, we will construct by recursion strictly pointwise increasing sequences $(H_i^A \in C^\infty_{A\subset M})_{i\in\N}$ of nondegenerate Hamiltonians for $A \in\{L,L',K,K'\}$, such that:
    \begin{itemize}
        \item For every $A\in\{L,L',K,K'\}$ and $i\geq 2$, $R_i^A<H_i^A$;
        \item For every $i\in\N$, $H_i^L<H_i^{L'}$ and $H_i^{K}<H_i^{K'}$;
        \item For every $t\in [0,1]$ and $i \in \N$, we have $H_i^L\leq H_i^K\circ \phi_t^{-1}$ and $H_i^{L'}\leq H_i^{K'}\circ \phi_t^{-1}$.
    \end{itemize}
    Once such a sequence has been constructed, the required acceleration $1$-rays can be defined by interpolation. Let us describe the construction of the above sequences.

    \textbf{Base of the recursion:} For $i=1$, take a $C^2$-small Morse function $H$ satisfying $-1<H<0$. Put
    \[H_1^L=H-2, \qquad H_1^{K}=H_1^{L'}=H-1\qquad\text{and}\qquad H_1^{K'}=H.\]
    These functions satisfy all the conditions above.

    \textbf{Step of the recursion:} Let $i\geq 2$, and assume we have already constructed finite sequences $(H_j^L)_{j=1}^{i-1},(H_j^{L'})_{j=1}^{i-1},(H_j^{K})_{j=1}^{i-1},(H_j^{K'})_{j=1}^{i-1}$ of nondegenerate Hamiltonians such that for every $1\leq j\leq i-1$ they satisfy the above conditions. We will construct Hamiltonians $H_i^L,H_i^{L'},H_i^{K},H_i^{K'}$ such that the above conditions still hold. Denote $\cP_i^L=\{R_i^L,H_{i-1}^L\}\subset C^\infty_{L\subset M}$. Since $L$ is compact and $\cP_i^L$ is nonempty and compact with respect to the $C^0$-topology, Lemma \ref{lemma: uniform_majorant}, applied to $L$ and $\cP_i^L$, yields a nondegenerate Hamiltonian $H_i^L\in C^\infty_{L\subset M}$ satisfying $H_{i-1}^L<H_{i}^L$, and $R_i^L<H_i^L$.

    Since $L'\subset L$, we have $H_i^L\in C^\infty_{L'\subset M}$. Let $\cP_i^{L'}=\{R_i^{L'},H_{i-1}^{L'},H_i^L\}\subset C^\infty_{L'\subset M}$; it is a nonempty $C^0$-compact set, and thus Lemma \ref{lemma: uniform_majorant}, applied to $L'$ and $\cP^{L'}_i$, produces a nondegenerate Hamiltonian $H_i^{L'}\in C^\infty_{L'\subset M}$ such that $H_{i-1}^{L'}<H_{i}^{L'}$, $R_i^{L'}<H_i^{L'}$, and $H_i^L<H_i^{L'}$.

    Since $\phi_t(K)\subset L$ for all $t$, it follows that $H^L_i\circ\phi_t|_{S^1\times K}<0$. Since the set
    \[\cP_i^{K}=\{R_i^{K},H_{i-1}^{K}\}\cup\{H_i^L\circ\phi_t \mid t\in[0,1]\}\subset C^\infty_{K\subset M}\]
    is $C^0$-compact, Lemma \ref{lemma: uniform_majorant}, applied to $K$ and $\cP^K_i$, furnishes a nondegenerate Hamiltonian $H_i^{K}\in C^\infty_{K\subset M}$ with $H_{i-1}^{K}<H_{i}^{K}$, $R_i^{K}<H_i^{K}$, and such that for every $t\in [0,1]$ we have $H_i^L\circ \phi_t<H_i^{K}$, that is, $H_i^L<H_i^{K}\circ \phi_t^{-1}$.

    Similarly, consider the $C^0$-compact set
    \[\cP_i^{K'}=\{R_i^{K'},H_{i-1}^{K'},H_i^{K}\}\cup\{H_i^{L'}\circ\phi_t \mid t\in[0,1]\}\subset C^\infty_{K'\subset M}\,.\]
    An application of Lemma \ref{lemma: uniform_majorant} to this set and to $K'$ yields a nondegenerate Hamiltonian $H_i^{K'}\in C^\infty_{K'\subset M}$ satisfying $H_{i-1}^{K'}<H_{i}^{K'}$, $R_i^{K'}<H_i^{K'}$, $H_i^{K}<H_i^{K'}$, and such that for every $t\in [0,1]$ we have $H_i^{L'}\circ \phi_t<H_i^{K'}$, that is, $H_i^{L'}<H_i^{K'}\circ \phi_t^{-1}$. This completes the recursive construction.

    We now use interpolation to produce the requisite acceleration $1$-rays. Fix a non-decreasing smooth function $\chi \fc [0,1]\to[0,1]$ with $\chi \equiv0$ near $0$ and $\chi\equiv 1$ near $1$.

    Let $A \in\{L,L',K,K'\}$. Define $\cH_\bullet^A=(\cH_i^A)_{i\in \N}$ as follows:
    \[\cH_i^A (s)=(1-\chi(s)) H_{i}^A+\chi(s) H_{i+1}^A\quad \text{for every }s\in[0,1]\,, i \in \N\,.\]
    This $\cH_\bullet^A$ is clearly a monotone $1$-ray. Note that for every $i\in \N$ and $s,t\in [0,1]$, we have
    $$\cH_i^L(s)=(1-\chi(s))H_{i}^L+\chi(s)H_{i+1}^L<(1-\chi(s))H_{i}^{K}\circ\phi_t^{-1}+\chi(s)H_{i+1}^{K}\circ\phi_t^{-1}=\cH_i^{K}(s)\circ\phi_t^{-1}\,,$$
    and similarly $\cH_i^{L'}(s)<\cH_i^{K'}(s)\circ\phi_t^{-1}$. This shows that $\cH_\bullet^L\leq \cH_\bullet^K\circ\phi_t^{-1}$ and $\cH_\bullet^{L'}\leq \cH_\bullet^{K'}\circ\phi_t^{-1}$ for all $t$. It remains to show that each $\cH_\bullet^A$ is an acceleration ray for $A$. This follows immediately from the fact that $R^A_i<\cH^A_i(0)=H_i^A$ for $i\geq2$ and the fact that the $R^A_i$ have the correct pointwise limits.
\end{proof}

The following lemma shows that interpolation between monotone cubes again yields monotone cubes.
\begin{lemma}\label{lemma:monotone filling}
    Let $\cH_0,\cH_1\colon [0,1]^n\to C^\infty(S^1\times M)$ be monotone Hamiltonian $n$-cubes. Let $\chi\colon [0,1]\to [0,1]$ be a nondecreasing smooth function satisfying $\chi\equiv 0$ near $0$ and $\chi\equiv 1$ near $1$. Let $\cH \fc [0,1]^{n+1} = [0,1]^n\times[0,1]\to C^\infty(S^1\times M)$ be given by
    \[\cH(x,s)=(1-\chi(s))\cH_0(x)+\chi(s)\cH_1(x),\quad \text{for }x\in[0,1]^n\,,s\in[0,1]\,.\]
    If $\cH_0\leq\cH_1$ then $\cH$ is a monotone Hamiltonian $(n+1)$-cube. Furthermore, if $\cH_0,\cH_1$ are nondegenerate at the vertices, then so is $\cH$.
\end{lemma}
\begin{proof}
    First, note that $\cH$ is constant on a neighborhood of each vertex, and equal to a nondegenerate Hamiltonian there, since $\cH_0,\cH_1$ are. It remains to verify that $\cH$ is nondecreasing along negative gradient lines of $\rho_{n+1}$, or equivalently that $d\cH(\nabla\rho_{n+1}) \leq 0$. We have $\nabla\rho_{n+1}=\nabla\rho_n\oplus\rho_1'\partial_s$, and thus for $x \in [0,1]^n$, $s\in [0,1]$ we have
    \begin{align*}
      d_{(x,s)}\cH(\nabla\rho_{n+1}) &= (1-\chi(s))d_x\cH_0(\nabla\rho_n) + \chi(s)d_x\cH_1(\nabla\rho_n)\\
                                      &\qquad\qquad -\chi'(s)\rho_1'(s)\cH_0(x) + \chi'(s)\rho_1'(s)\cH_1(x)\\
                                   &= (1-\chi(s))d_x\cH_0(\nabla\rho_n) + \chi(s)d_x\cH_1(\nabla\rho_n)\\
                                      &\qquad\qquad +\rho_1'(s)\chi'(s)(\cH_1(x)-\cH_0(x))\,.
    \end{align*}
    Since $\cH_0$, $\cH_1$ are monotone, we have $d\cH_0(\nabla\rho_n) \leq 0$, $d\cH_1(\nabla\rho_n) \leq 0$, whence the first two summands after the last equality sign are nonpositive. We also have $\chi'(s) \geq 0$, $\rho_1'(s) \leq 0$, and finally $\cH_1(x)-\cH_0(x) \geq 0$ by assumption. Thus the third summand is nonpositive as well.
\end{proof}

We can now complete the proof of Theorem \ref{thm:exists_parametrized_ray}.
\begin{proof}
    By Proposition \ref{prop:adapted_acceleration_data}, we have acceleration $1$-rays $\cH_\bullet^A$ for $A\in\{L,L',K,K'\}$, such that $\cH_\bullet^L\leq \cH_\bullet^{L'}$, $\cH_\bullet^K\leq \cH_\bullet^{K'}$ and such that for all $t$, $\cH_\bullet^L\leq \cH_\bullet^K\circ\phi_t^{-1}$ and $\cH_\bullet^{L'}\leq \cH_\bullet^{K'}\circ\phi_t^{-1}$.

    Pick a smooth nondecreasing $\chi \fc[0,1]\to[0,1]$ with $\chi\equiv 0$ near $0$ and $\chi\equiv 1$ near $1$. For all $i\in \N$, define Hamiltonian $2$-cubes $\cH^0_i,\cH^1_i$ by
    $$\cH^0_i(x) = (1-\chi(x_2))\cH^L_i(x_1) + \chi(x_2)\cH_i^{L'}(x_1)\,,\;\cH^1_i(x) = (1-\chi(x_2))\cH^K_i(x_1) + \chi(x_2)\cH_i^{K'}(x_1)\,$$
    where $x=(x_1,x_2) \in [0,1]^2$. It can be checked that $\cH_\bullet^j:=(\cH^j_i)_{i\in\N}$, $j=0,1$, is a Hamiltonian $2$-ray, meaning the constituent cubes satisfy the required compatibility at the boundary walls, that the $x_2=0$ wall of $\cH_\bullet^0$ is the ray $\cH_\bullet^L$, the $x_2=1$ wall is $\cH_\bullet^{L'}$, and that the $x_2=0$ wall of $\cH_\bullet^1$ is $\cH_\bullet^K$ while the $x_2=1$ wall is $\cH_\bullet^{K'}$. Furthermore, it is clear from the construction that for all $t$ we have $\cH_\bullet^0\leq \cH_\bullet^1\circ\phi_t^{-1}$. It also follows from Lemma \ref{lemma:monotone filling} that $\cH_\bullet^0$, $\cH_\bullet^1$ are monotone.

    For each $i\in\N$ define a Hamiltonian $4$-cube $\cG_i$ by
    $$\cG_i(x,z,r) = (1-\chi(z))\cH^0_i(x) + \chi(z)\cH^1_i(x)\circ\phi_r^{-1}\,,\quad\text{where }x\in[0,1]^2,z,r\in[0,1]\,.$$
    Let us show that $\cG_i$ is a $\bs\phi$-parametrized cube. For $z$ near $0$ we have $\chi(z)\equiv 0$, that is $\cG_i(x,z,r)=\cH^0_i(x)$, in particular it is independent of both $z$ and $r$. For $z$ near $1$ we have $\chi(z)\equiv 1$, meaning $\cG_i(x,z,r)=\cH^1_i(x)\circ\phi_r^{-1} = \cG_i|_0(x,1)\circ\phi_r^{-1}$ there, that is $\cG_i|_r(\cdot,z)$ is independent of $z$ while the dependence on $r$ has the specific form required by the definition. Since $\cH^j_i$ is monotone and for all $t$ we have $\cH^0_i\leq\cH^1_i\circ\phi_t^{-1}$, it follows from Lemma \ref{lemma:monotone filling} that each $\cG_i|_r$ is monotone. Thus, $\cG_i$ is a monotone $\bs\phi$-parametrized $4$-cube. It is also clear from the construction that the cubes $\cG_i$ are compatible at the $x_1=0,1$ walls, and therefore that $\cG_\bullet=(\cG_i)_{i\in\N}$ is a $4$-ray. Thus, $\cG_\bullet$ is a $\bs\phi$-parametrized monotone $4$-ray. Finally, $\cG_\bullet|_0$ has as its four corner $1$-rays the original acceleration $1$-rays $\cH_\bullet^A$, $A\in\{L,L',K,K'\}$, and thus it is an acceleration $3$-ray for the inclusion of pairs $(K,K')\subset(L,L')$. This completes the proof.
\end{proof}

\bibliographystyle{abbrv}
\bibliography{Bib_Symplectophagy}

\begin{thebibliography}{10}

\bibitem{AGV_2024_framed_E2_structures}
M.~Abouzaid, Y.~Groman, and U.~Varolgunes.
\newblock Framed {E}2 structures in {F}loer theory.
\newblock {\em Advances in Mathematics}, 450:109755, 2024.

\bibitem{DG_2026_relSH_in_CPn}
A.~Dickstein and Y.~Ganor.
\newblock Relative symplectic cohomology in complex projective spaces.
\newblock {\em arXiv preprint arXiv:2606.17370}, 2026.

\bibitem{DGPZ_2024_Symp_top_and_IVMs}
A.~Dickstein, Y.~Ganor, L.~Polterovich, and F.~Zapolsky.
\newblock Symplectic topology and ideal-valued measures.
\newblock {\em Selecta Mathematica}, 30(5):88, 2024.

\bibitem{DGZ_Var_conj}
A.~Dickstein, Y.~Ganor, and F.~Zapolsky.
\newblock Work in progress.

\bibitem{DZ_2026_constraints_on_symp_qs}
A.~Dickstein and F.~Zapolsky.
\newblock Constraints on symplectic quasi-states.
\newblock {\em International Journal of Mathematics}, 37(07):2650029, 2026.

\bibitem{Gromov_1985_J-hol_curves}
M.~Gromov.
\newblock Pseudo holomorphic curves in symplectic manifolds.
\newblock {\em Inventiones mathematicae}, 82(2):307--347, 1985.

\bibitem{Eliashberg_Gromov_1991_convex_symp_mfds}
M.~Gromov and Y.~Eliashberg.
\newblock Convex symplectic manifolds.
\newblock In {\em Proceedings of Symposia in Pure Mathematics}, pages 135--162.
  American Mathematical Society, 1991.

\bibitem{Hind_2016_Sympl_isotopies_dim_greater_four}
R.~Hind.
\newblock Symplectic isotopies in dimension greater than four.
\newblock {\em J. Symplectic Geom.}, 14(4):1033--1057, 2016.

\bibitem{Hind_Kerman_2026_Constraints_Lag_embeddings}
R.~Hind and E.~Kerman.
\newblock New constraints on {L}agrangian embeddings and the shape invariant.
\newblock {\em arXiv preprint arXiv: 2602.10291}, 2026.

\bibitem{Lerman_1995_symplectic_cut}
E.~Lerman.
\newblock Symplectic cuts.
\newblock {\em Mathematical Research Letters}, 2(3):247--258, 1995.

\bibitem{MSV_2024_heavy_sets_and_SH}
C.~Y. Mak, Y.~Sun, and U.~Varolgunes.
\newblock A characterization of heaviness in terms of relative symplectic
  cohomology.
\newblock {\em Journal of Topology}, 17(1):e12327, 2024.

\bibitem{MS_1995_Intro}
D.~McDuff and D.~Salamon.
\newblock Introduction to symplectic topology.
\newblock 1995.

\bibitem{McDuff_Traynor_1993_camel}
D.~McDuff and L.~Traynor.
\newblock The 4-dimensional symplectic camel and related results.
\newblock {\em Symplectic geometry}, 192:169--182, 1993.

\bibitem{Pardon_2016_alg_vir_counting}
J.~Pardon.
\newblock An algebraic approach to virtual fundamental cycles on moduli spaces
  of pseudo-holomorphic curves.
\newblock {\em Geometry \& Topology}, 20(2):779--1034, 2016.

\bibitem{Pardon_2019_contact}
J.~Pardon.
\newblock Contact homology and virtual fundamental cycles.
\newblock {\em Journal of the American Mathematical Society}, 32(3):825--919,
  2019.

\bibitem{Schlenk_book}
F.~Schlenk.
\newblock {\em Embedding problems in symplectic geometry}, volume~40.
\newblock Walter de Gruyter, 2005.

\bibitem{Spanier_1966_Alg_Top}
E.~H. Spanier.
\newblock {\em Algebraic topology}.
\newblock McGraw-Hill Book Co., New York-Toronto-London, 1966.

\bibitem{Tonkonog_Varolgunes_2023_Super_rigidity_of_certain_skeleta_using_relSH}
D.~Tonkonog and U.~Varolgunes.
\newblock Super-rigidity of certain skeleta using relative symplectic
  cohomology.
\newblock {\em Journal of Topology and Analysis}, 15(01):57--105, 2023.

\bibitem{Varolgunes_2018_PhD}
U.~Varolgunes.
\newblock Mayer--{V}ietoris property for relative symplectic cohomology (2018).

\bibitem{Varolgunes_2021_MV_and_relSH}
U.~Varolgunes.
\newblock Mayer--{V}ietoris property for relative symplectic cohomology.
\newblock {\em Geometry \& Topology}, 25(2):547--642, 2021.

\bibitem{Viterbo_1992_GenFun}
C.~Viterbo.
\newblock Symplectic topology as the geometry of generating functions.
\newblock {\em Mathematische Annalen}, 292(1):685--710, 1992.

\end{thebibliography}

\noindent
\begin{tabular}{l}
	{\bf Adi Dickstein} \\
	D\'epartement de Math\'ematiques et de Statistique\\
	Universit\'e de Montr\'eal\\
	CP 6128 succ Centre-Ville\\
	Montr\'eal QC H3C 3J7, Canada\\
	{\em E-mail:}  \texttt{adi.dickstein@gmail.com}
\end{tabular}

\medskip

\noindent
\begin{tabular}{l}
	{\bf Yaniv Ganor} \\
	School of Mathematical Sciences,\\
	Faculty of Sciences,\\
	HIT - Holon Institute of Technology \\
	
	Golomb St. 52, Holon 5810201, Israel\\
	{\em E-mail:}  \texttt{ganory@gmail.com}
\end{tabular}

\medskip

\noindent
\begin{tabular}{l}
	{\bf Frol Zapolsky} \\
	University of Haifa \\
	Department of Mathematics \\
	Faculty of Natural Sciences \\
	3498838, Haifa, Israel \\
	\&\\
	MI SANU \\
	Kneza Mihaila 36 \\
	Belgrade 11001\\
	Serbia\\
	{\em E-mail:}  \texttt{frol.zapolsky@gmail.com}
\end{tabular}

\end{document}